\documentclass[11pt,twoside]{amsart}

\usepackage{geometry}
\usepackage{amsmath,amssymb,esint,amsthm,bbm}
\usepackage{hyperref}
\catcode`@=11 \@addtoreset{equation}{section} \catcode`@=12

\email{dario.monticelli@gmail.com (Dario Daniele Monticelli)}
\email{scott.rodney@gmail.com (Scott Rodney)}
\email{wheeden@math.rutgers.edu (Richard L. Wheeden)}

\newcommand{\field}[1]{\mathbf{#1}}
\newcommand{\R}{\field{R}}

\newcommand{\N}{\field{N}}

\newcommand{\wt}{\widetilde}
\newcommand{\norm}{\parallel}

\newcommand{\bea}{\begin{eqnarray}}
\newcommand{\cend} {\end{center}}

\newcommand{\eea}{\end{eqnarray}}
\newcommand{\rn}{\mathbb{R}^n}
\newcommand{\fracc}{\displaystyle\frac}

\newcommand{\disp}{\displaystyle}

\newcommand{\ra}{\rightarrow}

\newcommand{\beginc} {\begin{center}}

\newcommand{\lp} {L^p(\Omega)}
\newcommand{\lipl} {\text{Lip}_\text{loc}(\Omega)}
\newcommand{\lipo} {\text{Lip}_0(\Omega)}

\newcommand{\waQ} {W^{1,p}_Q(\Omega)}
\newcommand{\wbQ} {(W^{1,p}_Q)_0(\Omega)}

\newcommand{\scl} {\mathcal{L}}
\newcommand{\ep} {\varepsilon}

\newcommand{\avnorm}[2] {\| {#1} \|_{{#2};\overline{dx}}}
\newcommand{\avnormb}[3]{\| {#1} \|_{{#2} ,{#3};\overline{dx}}}

\newcommand{\Lnormb}[3]{\| {#1} \|_{{#2} ,{#3};dx}}

\newtheorem{thm}{\textbf{Theorem}}[section]
\newtheorem{lem}[thm]{\textbf{Lemma}}
\newtheorem{pro}[thm]{\textbf{Proposition}}
\newtheorem{rem}[thm]{\textbf{Remark}}
\newtheorem{cor}[thm]{\textbf{Corollary}}
\newtheorem{defn}[thm]{\textbf{Definition}}
\theoremstyle{remark}

\theoremstyle{definition}
\newtheoremstyle{Claim}{}{}{\itshape}{}{\itshape\bfseries}{:}{ }{#1}
\theoremstyle{Claim}

\begin{document}

\title [Boundedness of weak solutions]{Boundedness of weak solutions
of degenerate quasilinear equations with rough coefficients}
\vspace{-0,3cm}



\maketitle

\begin{center}
\textsc{\textmd{D. D. Monticelli\footnote{Dipartimento di Matematica
      "F. Enriques", Universit\`a degli Studi, Via C. Saldini 50,
      20133--Milano, Italy. Partially supported by GNAMPA, project
      ``Equazioni differenziali e sistemi dinamici''. Partially
      supported by MIUR, project ``Metodi variazionali e topologici
      nello studio di fenomeni nonlineari''.},
      S. Rodney\footnote{Department of Mathematics, Physics and
      Geology, Cape Breton University, P.O. Box 5300, Sydney, Nova
      Scotia B1P 6L2, Canada.} and R. L. Wheeden\footnote{Department of
      Mathematics, Rutgers University, Piscataway, NJ 08854, USA.}}}
\end{center}



\begin{abstract}
We derive local boundedness estimates for weak solutions of a large
class of second order quasilinear equations. The structural
assumptions imposed on an equation in the class allow vanishing of
the quadratic form associated with its principal part and require no
smoothness of its coefficients. The class includes second order linear
elliptic equations as studied in \cite{GT} and second order
subelliptic linear equations as in \cite[2]{SW1}.  Our results also
extend ones obtained by J. Serrin \cite{S} concerning local
boundedness of weak solutions of quasilinear elliptic equations.
\end{abstract}

\maketitle

\section{\textbf{Introduction}}\label{intro}

The main purpose of this paper is to prove local boundedness of weak
solutions $u$ of rough subelliptic quasilinear equations of the form
\begin{equation}\label{eqdiff}
\text{div}\big(A(x,u,\nabla u)\big) = B(x,u, \nabla u)
\end{equation}
in an open set $\Omega \subset \R^n$. Further regularity results
will be studied in  a sequel to this paper. We will assume that the
vector-valued function $A$ and the scalar function $B$ satisfy
specific structural restrictions on their size, but not on their
smoothness, relative to a symmetric nonnegative semidefinite matrix
$Q(x)$. Thus the quadratic form $Q(x,\xi)= \langle
Q(x)\xi,\xi\rangle,\,\xi \in \R^n,$ may vanish when $\xi \neq 0$.

More precisely, given $p$ and an $n\times n$ matrix $Q$ with
$1<p<\infty$ and $|Q| \in L_{loc}^{p/2}(\Omega)$, our weak solutions
are pairs $(u, \nabla u)$ belonging to an appropriate Banach space
${\mathcal{W}}^{1,p}_Q(\Omega)$.  As described in \cite{SW2},
$\mathcal{W}^{1,p}_Q(\Omega)$ is obtained via isomorphism from the
degenerate Sobolev space $W^{1,p}_Q(\Omega)$ defined to be the completion
with respect to the norm
\begin{equation}\label{Wnorm}
 ||u||_{W_Q^{1,p}(\Omega)} = \left(\int_\Omega |u|^p \,dx +
   \int_\Omega Q(x, \nabla u)^{\frac{p}{2}} dx\right)^{\frac{1}{p}}
\end{equation}
of the class of functions in $\text{Lip}_{\text{loc}}(\Omega)$ with
finite $W^{1,p}_Q(\Omega)$ norm. We will give some further discussion
about these Banach spaces below. The structural conditions which we
assume are that there exists a vector ${\tilde A}(x,z,\xi)$,
$(x,z,\xi)\in \Omega\times\R\times\R^n$, with values in $\R^n$, such
that for a.e. $x\in\Omega$ and all $(z,\xi)\in\R\times\R^n$,
\begin{equation}\label{struct}
\begin{cases}
(i)\quad A(x,z,\xi)=\sqrt{Q(x)}{\tilde A}(x,z,\xi),\\
(ii)\quad\xi \cdot A(x,z,\xi) \geq a^{-1}\Big|\sqrt{Q(x)}\, \xi\Big|^p
- h|z|^\gamma - g,\\
(iii)\quad\Big|\tilde{A}(x,z,\xi)\Big| \leq a\Big|\sqrt{Q(x)}\,
\xi\Big|^{p-1} + b|z|^{\gamma -1} + e,\\
(iv)\quad\Big|B(x,z,\xi)\Big|\leq c\Big|\sqrt{Q(x)}\,
\xi\Big|^{\psi-1} +d|z|^{\delta-1}+f,
\end{cases}
\end{equation}
where $a, \gamma, \psi, \delta> 1$ are constants,
and $b, c, d, e, f, g, h$ are nonnegative functions of $x$. In fact,
when dealing with a particular weak solution $(u, \nabla u)$,  it
will be enough to assume that parts (ii), (iii) and (iv) of
(\ref{struct}) hold with $z$ and $\xi$ replaced by $u(x)$ and
$\nabla u(x)$ respectively.

The sizes of $\gamma, \psi$ and $\delta$ will be further restricted in
terms of $p$ and a natural ``Sobolev gain factor'' $\sigma > 1$ to be
described below in (\ref{Sobolev}), while the functions $b, c, d, e,
f, g, h$ will be assumed to lie in appropriate Lebesgue or Morrey
classes related to $p, \sigma, \gamma, \psi$ and $\delta$. For the
classical Euclidean metric, non-degenerate $Q$ and $1<p<n$, the
Sobolev gain is $\sigma = n/(n-p).$ In general, we will always
restrict $\gamma, \psi, \delta$ to the ranges
\begin{equation}\label{ranges}
\gamma\in (1,\sigma(p-1)+1),\quad \psi\in (1,p+1-\sigma^{-1}),\quad
\delta \in (1,p\sigma).
\end{equation}
We will often refer to $a,b,c, d,e,f,g,h$ as {\it structural
coefficients}, or simply as {\it coefficients}. Except for $a$, which
is constant, the coefficients must always satisfy certain minimal local
integrability requirements: see \eqref{integrability}.

We remark that the set of structural properties (\ref{struct}) is
invariant under replacing the symmetric nonnegative semidefinite
matrix $Q(x)$ by another symmetric nonnegative semidefinite matrix
$M(x)$ which is equivalent to it, i.e., which satisfies
\begin{equation}\label{47}
\frac{1}{C}\langle Q(x)\xi,\xi\rangle\leq\langle
M(x)\xi,\xi\rangle\leq C\langle Q(x)\xi,\xi\rangle\qquad\text{for
all }\xi\in\R^n,\,\text{ a.e. }x\in\Omega;
\end{equation}
see Theorem \ref{32} in Appendix 2.

We also remark that the structural assumptions (\ref{struct}) are
equivalent to the following set of assumptions: there exists a
nonnegative function $\tilde{a}(x,z,\xi)$, $(x,z,\xi)\in
\Omega\times\R\times\R^n$, such that for a.e. $x\in\Omega$, all
$\eta\in\R^n$ and all $(z,\xi)\in\R\times\R^n$,
\begin{equation}\label{struct3}
\begin{cases}
(i)\quad \big|\eta\cdot A(x,z,\xi)\big| \le
\big|\sqrt{Q(x)}\eta\big|\, \tilde{a}(x,z,\xi),\\
(ii)\quad \xi \cdot A(x,z,\xi) \geq a^{-1}\Big|\sqrt{Q(x)}\,
\xi\Big|^p - h|z|^\gamma - g,\\
(iii)\quad \tilde{a}(x,z,\xi) \leq a\Big|\sqrt{Q(x)}\, \xi
\Big|^{p-1} + b|z|^{\gamma -1} + e\\
(iv)\quad \Big|B(x,z,\xi)\Big|\leq c\Big|\sqrt{Q(x)}\, \xi
\Big|^{\psi-1}+d|z|^{\delta-1}+f;
\end{cases}
\end{equation}
see Theorem \ref{36} in Appendix 2.

Historically speaking, in the classical elliptic case ($Q(x) =$
Identity), structural conditions more restrictive than (\ref{struct})
were considered by J. Serrin \cite{S}, who derived a broad class of
regularity results for weak solutions  of (\ref{eqdiff}). Our ranges
of the parameters $\gamma, \psi, \delta$ are however wider than those
studied in \cite{S}, where these parameters are all equal to $p$. In
case $p=2$, our ranges correspond more closely to those in
\cite[p. 176]{G} and \cite{GM}, although we miss some endpoint
values. These latter papers impose continuity conditions on
coefficients which we do not assume, but which lead to stronger
regularity and also to results for systems. N. Trudinger
\cite{T} also derived regularity results in the elliptic case,
relaxing some structural conditions under the assumption of local
boundedness of weak solutions, but generally for the same choices as
in \cite{S}. We note in passing that the equation for the
$p$-Laplacian, namely $\text{div} (\nabla u\, |\nabla u|^{p-2})=0$,
as well as Yamabe type equations $\Delta u - Ru + {\bar R} u^{q-1}=0$
for $q<2n/(n-2)$ are included in the case $Q(x) =$ Id (with $p=2$ for
the Yamabe type equations).

In the subelliptic case, by which we mean the case when $Q(x)$ may
be singular, regularity results including local boundedness of weak
solutions are derived in \cite[2]{SW1} for {\it linear} equations
with rough coefficients and nonhomogeneous terms. The form of the
linear equations studied there is
\begin{equation}\label{linear}
\mbox{div }\big(M(x)\nabla u\big) +H(x)R(x)u + S(x)'G(x)u +F(x)u =
F_1(x) + T(x)'G_1(x),
\end{equation}
where $M(x)$ is a symmetric matrix whose quadratic form $M(x,\xi)$
satisfies
\[
c_1 Q(x,\xi) \le M(x,\xi) \le c_2 Q(x,\xi)
\]
for some positive constants $c_1, c_2$, and where $H, G, F, F_1,
G_1$ are functions on $\Omega$ and $R, S, T$ are vector fields
on $\Omega$ that are subunit with respect to $Q(x)$. Here we say
that a vector field $V(x)$ is {\it subunit with respect to $Q$} if
$V(x)u(x) = v(x)\cdot \nabla u(x)$ and $|v(x)\cdot \xi|^2 \le
Q(x,\xi)$ for almost all $x\in \Omega$, all $\xi \in \R^n$ and all
Lipschitz continuous functions $u$ on $\Omega$. By direct
computation (see Theorems \ref{35} and \ref{32} and relations
(\ref{48}) in Appendix 2), such a linear subelliptic equation
satisfies the structural conditions (\ref{struct}) with $p = \gamma
= \psi = \delta = 2$. Our principal result includes the local
boundedness estimates in \cite[2]{SW1} for solutions of (\ref{linear})
and can be viewed as an extension of them to solutions of quasilinear
equations.

The regularity results in \cite[2]{SW1} for equations of type
(\ref{linear}) were derived in an axiomatic setting which assumes
the existence of appropriate Sobolev-Poincar\'e estimates in a space
of homogeneous type, as well as the existence of sequences of
Lipschitz cutoff functions. We will derive our estimates for weak
solutions of (\ref{eqdiff}) in a quasimetric setting, but our axioms
are generally less restrictive than those in \cite[2]{SW1}. For
example, we do not need the assumption in \cite[2]{SW1} that Lebesgue
measure satisfies the doubling property relative to quasimetric
balls. In fact, our main result Theorem \ref{main_thm}
requires no hypothesis at all about doubling, and Corollaries
\ref{allp}--\ref{second>p} use only the condition $D_{q^*}$ listed
in Definition \ref{weakdoub}. Also, our main Sobolev-Poincar\'e
assumption will be one of Sobolev type for compactly supported
functions. Unlike \cite[2]{SW1}, where not only local boundedness but
also H\"older continuity of weak solutions is obtained, we will not
require any Poincar\'e estimate for non-compactly supported
functions. However, depending on the order of integrability of $|Q|$,
we sometimes assume that functions $w$ in $W^{1,p}_Q(\Omega)$ satisfy
higher local integrabililty than order $p$.

In order to state our main theorem, we now briefly describe the
axiomatic framework. A fuller discussion can be found in
\cite[2]{SW1}. Facts about degenerate Sobolev spaces
$W_Q^{1,p}(\Omega$) are given in \cite{SW2}, as well as in \cite{R}
when $p=2$, and we now recall some of them. Let $Lip_{Q,p}(\Omega)$
denote the class of locally Lipschitz functions with finite
$W^{1,p}_Q(\Omega)$ norm. By definition, $W^{1,p}_Q(\Omega)$ is the
Banach space of equivalence classes of Cauchy sequences in
$Lip_{Q,p}(\Omega)$. Consider the
form-weighted space consisting of all measurable
$\R^n$-valued functions $\mathbf{f}(x), x\in \Omega$, satisfying
\begin{equation}\label{calLp}
||\mathbf{f}||_{\mathcal{L}^p(\Omega,Q)} = \left\{\int_\Omega Q(x,
  \mathbf{f}(x))^{\frac{p}{2}} dx \right\}^{\frac{1}{p}} < \infty.
\end{equation}
Identifying measurable $\R^n$-valued functions $\mathbf{f},
\mathbf{g}$ which satisfy
$||\mathbf{f}-\mathbf{g}||_{\mathcal{L}^p(\Omega,Q)} =0$,
(\ref{calLp}) defines a norm on the resulting vector space
of equivalence classes; we define $\mathcal{L}^p(\Omega, Q)$ as the
space consisting of these
equivalence classes. When $p=2$, $\mathcal{L}^2(\Omega, Q)$ is shown
to be a Hilbert space with inner product $\langle \mathbf{f}, \mathbf{g}
\rangle = \int_\Omega \mathbf{f}(x)'Q(x) \mathbf{g}(x)\,dx$ in Theorem
4 of \cite{SW2}, and the arguments in the proof there show that
$\mathcal{L}^p(\Omega, Q)$ is a Banach space with norm (\ref{calLp})
for $1\le p <\infty$.   If $\{w_k\}_{k=1}^\infty\in
W^{1,p}_Q(\Omega)$, i.e., if $\{w_k\}_{k=1}^\infty$ is a Cauchy
sequence in $W^{1,p}_Q(\Omega)$ norm of $Lip_{Q,p}(\Omega)$
functions, there is a unique pair $(w, \mathbf{v})\in L^p(\Omega)\times
\mathcal{L}^p(\Omega,Q)$ such that $w_k\rightarrow w$ in $L^p(\Omega)$
and $\nabla w_k \rightarrow \mathbf{v}$ in $\mathcal{L}^p(\Omega,
Q)$. The pair $(w, \mathbf{v})$ represents the equivalence class in
$W^{1,p}_Q(\Omega)$ which contains the Cauchy sequence $\{w_k\}$. Any
pair $(w,\mathbf{v})$ representing an equivalence class in
$W^{1,p}_Q(\Omega)$ is said to belong to the space
${\mathcal{W}}^{1,p}_Q(\Omega)$.  Thus, $\mathcal{W}^{1,p}_Q(\Omega)$
is the image of the isomorphism $\mathcal{J}:W^{1,p}_Q(\Omega)\ra
L^p(\Omega)\times \mathcal{L}^p(\Omega,Q)$ defined by
\bea \mathcal{J}([\{w_k\}]) = (w,\mathbf{v}),\nonumber
\eea
where $[\{w_k\}]$ denotes the equivalence class in
$W^{1,p}_Q(\Omega)$ containing the Cauchy sequence $\{w_k\}$.
Therefore, $\mathcal{W}^{1,p}_Q(\Omega)$ is a closed subspace of
$L^p(\Omega)\times \mathcal{L}^p(\Omega,Q)$ and hence a Banach space
as well.  As the spaces ${\mathcal{W}}^{1,p}_Q(\Omega)$ and
$W^{1,p}_Q(\Omega)$ are isomorphic, we will often refer to elements
$(w,\mathbf{v})$ of $\mathcal{W}^{1,p}_Q(\Omega)$ as elements of
$W^{1,p}_Q(\Omega)$, where the isomorphism is taken in context.  We
caution the reader that $\mathbf{v}$ is not generally uniquely
determined by $w$ for pairs $(w, \mathbf{v})$ in
$\mathcal{W}_Q^{1,p}(\Omega)$, i.e., the projection
\[
P: \mathcal{W}^{1,p}(\Omega) \rightarrow L^p(\Omega)
\]
obtained by mapping a pair onto its first component is not always an
injection, as shown by an example in \cite{FKS}.
Nevertheless, we will generally abuse notation and denote
representative pairs in $W^{1,p}_Q(\Omega)$ by $(w, \nabla w)$
instead of $(w, \mathbf{v})$. Moreover, we will often abuse notation
even further by simply writing $w$ instead of the pair $(w,\nabla
w)$. Some additional facts about degenerate Sobolev spaces are listed
in Section 2 and Appendix 1 below.

In \cite{SW2}, the notion of the {\it regular} gradient
$\nabla_{reg}w$ of an element $w$ in $W^{1,p}_Q(\Omega)$ is introduced
and used to derive results related to regularity of linear subelliptic
equations. However, in the present paper, we have been able to avoid
this technical device; see the comment which follows Corollary
\ref{goodlevelsets2}.

By a quasimetric $\rho$ on $\Omega$, we mean a finite nonnegative
function on $\Omega\times \Omega$ such that for some constant
$\kappa \ge 1$,
\[
\rho(x,y) = 0 \mbox{ iff } x=y
\]
\[
\rho(x,y) \le \kappa [\rho(x,z)+ \rho(y,z)] \mbox{ if } x,y,z\in
\Omega.
\]
For simplicity, we will also assume that $\rho$ is symmetric, i.e.,
that $\rho(x,y)= \rho(y,x)$ if $x,y \in \Omega$. For $x\in\Omega$
and $r>0$, define the sets
\begin{eqnarray}
\nonumber B(x,r)\!&\!=\!&\!\{y\in\Omega: \rho(x,y)<r\},\\
\nonumber D(x,r)\!&\!=\!&\!\{y\in\Omega:  |x-y|<r\},
\end{eqnarray}
and assume that $B(x,r)$ is Lebesgue measurable for every $r>0$,
$x\in\Omega$. We call $B(x,r)$ the quasimetric ball (or $\rho$-ball)
with center $x$ and radius $r$, and we sometimes write $B_r(x)$
or simply $B_r$ instead of $B(x,r)$. Throughout the paper we will
assume that
\begin{eqnarray}
\label{Cond1}&& \text{for all } x\in\Omega, \, |x-y|\rightarrow 0
\text{ if } \rho(x,y)\rightarrow 0,
\end{eqnarray}
and in some of our results we will also assume that
\begin{eqnarray}
\label{Cond0}&& \text{for all } x\in\Omega,\, \rho(x,y)\rightarrow 0
\text{ if } |x-y|\rightarrow 0.
\end{eqnarray}
We remark that condition \eqref{Cond1} is equivalent to requiring
that
\begin{equation}\label{Cond1.1}
\begin{array}{l}
\text{for every }x\in\Omega\text{ and every }\epsilon>0\text{ there
exists }\delta>0,\text{ depending on }x\text{ and }\epsilon,\\
\text{such that }B(x,\delta)\subset D(x,\epsilon),
\end{array}
\end{equation}
while condition \eqref{Cond0} is equivalent to
\begin{equation}\label{Cond0.1}
\begin{array}{l}
\text{for every }x\in\Omega\text{ and every }r>0\text{ there
exists }s>0,\text{ depending on }x\text{ and }r,\\
\text{such that }D(x,s)\subset B(x,r).
\end{array}
\end{equation}
Then condition \eqref{Cond0}, or equivalently condition
\eqref{Cond0.1}, implies that $|B(x,r)|>0$ for every
$\rho$-ball with $r>0$. By Lemma \ref{L0} in Section \ref{section2},
condition \eqref{Cond1} implies that for every $x\in\Omega$, one has
$\overline{B(x,r)}\subset\Omega$ if $r$ is smaller than a suitable
$r_0=r_0(x)>0$; here $\overline{E}$ denotes the Euclidean closure of a
set $E\subset\Omega$.

Given $p$, $1<p < \infty$, and a nonnegative semidefinite quadratic
form $Q(x,\xi) = \langle Q(x)\xi,\xi\rangle$, where $Q(x)$ is a
symmetric matrix for each $x\in \Omega$ and $|Q|\in
L_{loc}^{p/2}(\Omega)$, we need the following Sobolev estimate:
there exist $\sigma >1$ and $C>0$ such that for every $\rho$-ball
$B_r= B(y,r)$ with $y\in \Omega$ and $0<r<r_1(y)$ for a suitable
$r_1(y)>0$,
\begin{equation}\label{Sobolev}
\left(\frac{1}{|B_r|}\int_{B_r} |w|^{p\sigma}
dx\right)^{\frac{1}{p\sigma}}\le C \left[
r \left(\frac{1}{|B_r|}\int_{B_r} Q(x, \nabla w)^{\frac{p}{2}} dx
\right)^{\frac{1}{p}}+\left(\frac{1}{|B_r|}\int_{B_r} |w|^p
dx\right)^{\frac{1}{p}}\right]
\end{equation}
for all pairs $(w, \nabla w) \in (W_Q^{1,p})_0(B_r)$. Here
$(W_Q^{1,p})_0(B_r)$ denotes the analogue of the space
$W_Q^{1,p}(B_r)$ defined earlier but now the completion with
respect to (\ref{Wnorm}), with $\Omega$ now replaced by $B_r$, is
formed by using Lipschitz functions with compact support in
$B_r$. Even though $\nabla w$ may not be determined uniquely by $w$,
it follows that (\ref{Sobolev}) holds for all $(w,\nabla w) \in
(W_Q^{1,p})_0(B_r)$ provided it holds for all Lipschitz functions with
compact support in $B_r$. We also note that since $Q(x,\xi) =
|\sqrt{Q(x)}\, \xi|^2$, (\ref{Sobolev}) can be rewritten as
\[
\left(\frac{1}{|B_r|}\int_{B_r} |w|^{p\sigma}
dx\right)^{\frac{1}{p\sigma}}\le C\left[
r \left(\frac{1}{|B_r|}\int_{B_r} \big|\sqrt{Q} \nabla w\big|^p dx
\right)^{\frac{1}{p}}+\left(\frac{1}{|B_r|}\int_{B_r} |w|^p
dx\right)^{\frac{1}{p}}\right].
\]
The number $\sigma$ is a factor which measures the ``Sobolev gain''
in integrability of $w$, from $L^p(B_r)$ to $L^{p\sigma}(B_r)$
independently of $B_r$; $\sigma$ plays a crucial role in our
results.

We will always assume that $r_1(y)\leq r_0(y)$ for every $y\in\Omega$,
where $r_1(y)$ is as in \eqref{Sobolev} and $r_0(y)$ is as in Lemma
\ref{L0} In particular, it then follows that the closure of any
ball $B(y,r)$ with $r<r_1(y)$ lies in $\Omega$.

We also require the existence of appropriate sequences of Lipschitz
cutoff functions (called ``accumulating sequences of Lipschitz
cutoff functions'' in \cite{SW1}), namely, we require that for some
exponent $s^*$,  $p\sigma^\prime<s^* \le \infty$, there are positive
constants $\tau, N$ and $C_{s^*}$, with $\tau <1$, such that for
every ball $B(y,r)$ with $0< r< r_1(y)$, there
is a sequence of Lipschitz functions $\{\eta_j\}_{j=1}^\infty$ with
the properties
\begin{equation}\label{cutoff}
\begin{cases}
\text{supp}\,\eta_1\subset B(y,r)\\
0\leq\eta_j\leq 1\quad \text{ for all }j\geq 1\\
B(y,\tau r)\subset\{x\in B(y,r):\eta_j(x)=1\}\quad\text{for all $j\ge 1$}\\
\text{supp}\,\eta_{j+1}\subset\{x\in
B(y,r):\eta_j(x)=1\}\quad\text{for all $j\ge 1$}\\
\displaystyle{\left(\frac{1}{|B(y,r)|}\int_{B(y,r)}\big|\sqrt{Q}\nabla\eta_j
\big|^{s^*} dx\right)^{\frac{1}{s^*}}\leq
C_{s^*}\frac{N^j}{r}}\quad\text{for all $j\ge 1$.}
\end{cases}
\end{equation}
We remark that the above condition is slightly different from that
appearing in \cite{SW1}, and it is actually weaker.  Indeed, the key final
property in \eqref{cutoff} is weaker than its analogue in \cite{SW1}
where the exponential growth constant $N^j$ is replaced by
$j^N$. Further, it is assumed in \cite{SW1} that
$r_1(y)=\delta_0\text{dist}(y,\partial\Omega)$ for some $\delta_0>0$,
where ``$`\text{dist}$'' denotes the standard Euclidean distance in
$\R^n$. The second property, $0\leq \eta_j\leq 1$ for every $j$,
is not required in \cite{SW1}.  However, if $\{\eta_j\}_{j=1}^\infty$
is some collection which satisfies \eqref{cutoff} except for the
second part, simply define a new collection
$\{\tilde{\eta}_j\}_{j=1}^\infty$ by
\begin{equation*}
\tilde{\eta}_j(x) = \left\{
\begin{array}{cl}
\eta_j(x) & \quad\text{if } 0<\eta_j(x)<1,\\
0 & \quad\text{if } \eta_j(x)\leq 0,\\
1 & \quad\text{if } \eta_j(x)\geq 1.
\end{array} \right.
\end{equation*}
for each $j$.  This new collection then satisfies \eqref{cutoff} as
written. We also remark that since $s^*>p\sigma'$, we may choose a
number $s'>\sigma'$ so that $s^* = s'p$.  The exponent $s$ which is
dual to $s^\prime$, i.e., so that $1/s + 1/s' =1$, satisfies $1 \le
s <\sigma$ and plays an important role in our results.

Another assumption, generally simpler than \eqref{cutoff}, which we
will impose in our main theorem is that there exists $t$, $1\le t\le
\infty$, such that for every $\rho$-ball $B(y,r)$ with $0< r<
r_1(y)$ and every $\eta= \eta_j$ in the corresponding sequence
$\{\eta_j\}$ provided by \eqref{cutoff},
\begin{equation}\label{E3}
\left(\int_{B(y,r)}|\sqrt{Q}\nabla\eta|^{tp}\,dx\right)^{\frac{1}{tp}}
<\infty.
\end{equation}
In fact, by \eqref{cutoff}, condition \eqref{E3} is automatically
satisfied for every $t$ with $1\leq t\leq s^*/p$. On the other hand,
(\ref{E3}) might hold for larger values of $t$ independently of
\eqref{cutoff}; for example, if $Q(x)$ is bounded, then \eqref{E3}
holds with $t= \infty$ for all $B(y,r)$ with closure in $\Omega$ and for
every $\eta \in Lip_0(\Omega)$, even if \eqref{cutoff} is not valid
for any $s^*$. To derive some of the preliminary results in Section 2,
we will assume \eqref{E3} for more restricted classes of balls
$B(y,r)$ and functions $\eta$. In any case, \eqref{E3} as
well as \eqref{3.6-0} below are only qualitative conditions, in the
sense that the constants involved in both of them will not enter our
final estimates.

In our main theorem, \eqref{E3} will be paired with the following
assumption, where $t'$ is the usual dual index of $t$ given by
$1/t+1/t' =1$: for every $\rho$-ball $B(y,r)$ with $0 <r < r_1(y)$,
there is a constant $c_2=c_2(B(y,r))$ so that
for all $f\in Lip_{loc}(\Omega)$,
\begin{equation}\label{3.6-0}
\bigg(\int_{B(y,r)} | f|^{pt'}dx\bigg)^\frac{1}{pt'} \leq
c_2\,||f||_{W^{1,p}_Q(\Omega)} = c_2\bigg(\int_{\Omega}|\sqrt{Q}\nabla f|^pdx
+\int_{\Omega}|f|^pdx\bigg)^\frac{1}{p}.
\end{equation}
It is easy to see that condition \eqref{3.6-0} holds for all elements
of $W^{1,p}_Q(\Omega)$ and not just for functions in
$Lip_{loc}(\Omega)$. 

In Section 2, \eqref{E3} and \eqref{3.6-0} will
be used to derive a useful version of the product rule. They will also be
used to prove that functions in $W^{1,p}_Q(\Omega)$, which are
generally without compact support, have sufficiently high local
integrability in the presence of the Sobolev estimate
\eqref{Sobolev} for compactly supported ones. See Proposition
\ref{prop1} for an estimate of $||w||_{L^{p\sigma}(B(y,\tau r))},
0<\tau <1,$ in case $w \in W^{1,p}_Q(\Omega)$ and $B(y,r)$ is any
$\rho$-ball with $0<r< r_1(y)$.  As is true for \eqref{E3}, we
sometimes assume in Section 2 that \eqref{3.6-0} holds for a smaller
class of balls.

\begin{rem}
Note that condition \eqref{3.6-0} becomes weaker as $t'$ becomes
smaller. In particular, if \eqref{E3} holds with $t=\infty$ (e.g., if
$Q \in L^\infty_{\text{loc}}(\Omega)$ or if \eqref{cutoff} is valid
with $s^* = \infty$), then $t^\prime=1$ and \eqref{3.6-0} is
trivially true.\\
\indent When \eqref{cutoff} holds for some $s^*>p\sigma'$, then
\eqref{E3} is automatically true with $t= s^*/p$, and the
corresponding $t'$ in (\ref{3.6-0}) satisfies $1\leq t'<\sigma$.
In case $t'<\sigma$,  (\ref{3.6-0}) is considerably weaker than the
Sobolev inequality \eqref{Sobolev} when restricted to Lipschitz
functions $f$ with compact support in $B(y,r)$. On the other hand,
\eqref{3.6-0} is assumed to hold for any locally Lipschitz function
whether it is compactly supported in $B(y,r)$ or not. \\
\indent In case the Poincar\'e inequality
$$
\left(\frac{1}{|B_r|}\int_{B_r} |f- f_{B_r}|^{pt'}
dx\right)^{\frac{1}{pt'}}\le C r \left(\frac{1}{|B_r|}\int_{B_r} Q(x,
\nabla f)^{\frac{p}{2}} dx \right)^{\frac{1}{p}}, \quad f_{B_r}
=\frac{1}{|B_r|}\int_{B_r} f dx,
$$
holds with $B_r = B(y,r)$ for all $f\in Lip_{loc}(\Omega)$, then
\eqref{3.6-0} clearly holds as well.
\end{rem}

In many cases of interest, conditions
(\ref{Sobolev}), (\ref{cutoff}), (\ref{E3}) and \eqref{3.6-0}
automatically hold. An enormous related literature exists, and we
refer to \cite{SW1} for an introduction to it. In particular,
(\ref{cutoff}) is known to hold with $s^*=\infty$ for the {\it
subunit} balls $\mathcal{K}(x,r)$ associated with a quadratic form
$Q(x,\xi)$ that is continuous in $x$, provided the Fefferman-Phong
condition \cite{FP} holds, i.e., provided there are positive constants
$c_0, \epsilon$ such that for every $\mathcal{K}(x,r)$ with closure in
$\Omega$, there is a Euclidean ball $D(x,r)$ satisfying
\[
D(x,r) \subset \mathcal{K}(x,c_0r^\epsilon).
\]
Notice that this condition in particular implies \eqref{Cond0.1},
i.e., condition \eqref{Cond0}, for the subunit balls
$\mathcal{K}(x,r)$. In order to elaborate, we extend (as in
\cite{SW1}) the notion of
subunit metric to a nonnegative continuous quadratic form $Q(x,\xi)$
on $\Omega$ by defining
\begin{equation}\label{subunitmetric}
\delta(x,y) = \inf \{r>0: \gamma(0)=x, \gamma(r) = y, \gamma \mbox{
is a Lipschitz subunit curve in $\Omega$}\},
\end{equation}
where a Lipschitz curve $\gamma: [0,r] \rightarrow \Omega$ is said
to be {\it subunit} (with respect to $Q(x,\xi))$ if
\[
(\gamma'(t)\cdot \xi)^2 \le Q(\gamma(t), \xi)
\]
for a.e. $t\in [0,r]$ and all $\xi\in \R^n$. Then $\delta(x,y)$ is a
symmetric metric on $\Omega$, although possibly infinite if $Q$ is
degenerate. If $\delta(x,y)$ is finite for all $x,y\in \Omega$, the
subunit balls $\mathcal{K}(x,r)$ are defined by
\begin{equation}\label{subunitballs}
\mathcal{K}(x,r) = \{y\in \Omega: \delta(x,y) <r\},\quad x\in \Omega,\,
0<r<\infty.
\end{equation}
Assuming that $Q$ is continuous, that $\delta(x,y)$ is finite, and
that the Fefferman-Phong containment condition holds, it is shown in
\cite{SW1} (and, under more restrictive assumptions, in the related
references listed there) that (\ref{cutoff}) holds with $s^*=\infty$
for the balls $\mathcal{K}(x,r)$.

We say that a pair $(u, \nabla u) \in W_Q^{1,p}(\Omega)$ is a {\it
weak solution} of (\ref{eqdiff}) if
\begin{equation}\label{weaksolution}
\int_\Omega \big[\nabla\varphi\cdot A(x,u,\nabla u) +\varphi
B(x,u,\nabla u)\big]\,dx =0 \quad\mbox{ for all $\varphi \in
\text{Lip}_0(\Omega),$}
\end{equation}
where $\text{Lip}_0(\Omega)$ denotes the class of Lipschitz
functions with compact support in $\Omega$.

The main results of this paper are the following theorem and
corollaries, in which we will use the notation
\begin{equation*}
\Lnormb{f}{\alpha}{E}=\left(\int_E|f(x)|^\alpha\,dx
\right)^\frac{1}{\alpha},\qquad
\avnormb{f}{\alpha}{E}:=\left(\frac{1}{|E|}\int_E|f(x)|^\alpha\,
dx\right)^\frac{1}{\alpha}
\end{equation*}
whenever $E\subset\Omega$ is Lebesgue measurable, $f$ is a Lebesgue
measurable function on $E$, and $\alpha >0$.

\begin{thm}\label{main_thm}
Let $\Omega$ be an open set in $\R^n$, $1<p<\infty$, and $Q(x,\xi) =
\langle Q(x)\xi,\xi\rangle$ be a symmetric nonnegative semidefinite
quadratic form on $\Omega$ with $|Q|\in L_{loc}^{p/2}(\Omega)$.
Suppose that $(\Omega, \rho)$ is a quasimetric space, that condition
\eqref{Cond1} holds, and that there exists $\sigma >1$
such that the Sobolev estimate \eqref{Sobolev} holds for all $(w,
\nabla w) \in \big(W_Q^{1,p}\big)_0(B)$ for all $\rho$-balls
$B=B(y,r)$ with $0<r< r_1(y)$. Let $A(x,z,\xi)$ and $B(x,z,\xi)$
satisfy the structural conditions (\ref{struct}) with
\begin{equation}\label{ranges1}
\gamma=\delta=p, \quad \psi \in [p, p+1-\sigma^{-1}).
\end{equation}
Suppose that condition \eqref{cutoff} about Lipschitz cutoff
functions holds for some $\tau\in(0,1)$ and $s^*>p\sigma^\prime$,
and that conditions \eqref{E3} and \eqref{3.6-0} hold with
$1/t+1/t'=1$ for some $t\geq1$ and all $\rho$-balls $B$ as above.
Let $(u, \nabla u)\in W_Q^{1,p}(\Omega)$ be a weak solution of
\eqref{eqdiff} in $\Omega$ and let $B(y,r)$ be a $\rho$-ball
with $0<r<\tau r_1(y)$.

Furthermore, given $k>0$ and
$\epsilon_1,\epsilon_2,\epsilon_3\in(0,1]$, let
\begin{equation*}
\begin{array}{rclcrcl}
\bar{u}&=&|u|+k,&\qquad&\bar{b}&=&b+k^{1-p}e,\\
\bar{h}&=&h+k^{-p}g,&\qquad&\bar{d}&=&d+k^{1-p}f,
\end{array}
\end{equation*}
and define
\begin{eqnarray} \label{Zetabar}
{\bar Z}\!&\!=\!&\! 1+r^{p-1}\avnormb{\bar{b}}{{p'\sigma'}}{B(y,r)}+
\left(r^p \avnormb{c^\frac{p}{p+1-\psi}
\bar{u}^\frac{p(\psi-p)}{p+1-\psi}}{\frac{p
\sigma'}{p-\epsilon_1}}{B(y,r)}\right)^{\frac{1}{\epsilon_1}}\\
\nonumber\!&\!\!&\!\hspace{2cm} +\left(r^p
\avnormb{\bar{h}}{\frac{p\sigma'}{p-\epsilon_2}}{B(y,r)}
\right)^\frac{1}{\epsilon_2}+ \left(r^p
\avnormb{\bar{d}}{\frac{p\sigma'}{p-\epsilon_3}}{B(y,r)}
\right)^\frac{1}{\epsilon_3}.
\end{eqnarray}

Then
\begin{equation}\label{Linfinityestimate}
\|\bar{u}\|_{L^\infty(B(y,\tau r))}\leq C{\bar
Z}^{\Psi_0}\avnormb{\bar{u}}{sp}{B(y,r)},
\end{equation}
where $s$ is the dual exponent of the number $s^\prime$ which
satisfies $s^*=s'p$, $C$ is a constant independent of $u,k, B(y,r),
b,c,d,e,f,g$ and $h$, and where $\Psi_0=\frac{s}{\sigma-s}$.
\end{thm}

\begin{rem}\label{useCond0}
Under the hypothesis of Theorem \ref{main_thm}, Proposition
\ref{prop1} guarantees that the factor $\avnormb{\bar{u}}{sp}{B(y,r)}$
in \eqref{Linfinityestimate} is finite. If $\bar{Z}$ is finite for all
$B(y,r)$ as above and condition \eqref{Cond0} is satisfied, then
Theorem \ref{main_thm} gives local boundedness of weak solutions of
equation \eqref{eqdiff} in $\Omega$, i.e., weak solutions in
$W^{1,p}_Q(\Omega)$ are bounded in every compact subset of $\Omega$;
see Section 2 for the simple proof.
\end{rem}

\begin{rem}\label{otherserrinestimate}
The proof of Theorem \ref{main_thm} also provides an $L^p$ estimate for
the size of $\sqrt{Q}\nabla u$ when $(u,\nabla u)$ is a weak
solution. In fact, under the same assumptions as in Theorem
\ref{main_thm},
\begin{equation}\label{ose}
||\sqrt{Q}\nabla u||_{L^p(B(y,\tau r))} \le C \bar{Z}\left(\frac{1}{r}
  ||\bar{u}||_{L^p(B(y,r))} + ||\bar{u}||_{L^{t'p}(B(y,r))}\right),
\end{equation}
where the norms are now unnormalized. This estimate is an analogue of
one obtained in \cite{S} in the nondegenerate case. It follows from
\eqref{ne16} below by choosing $q=1$ and $\eta =\eta_1$ there, and by
applying \eqref{3.6-0} to the first term on the right in \eqref{ne16}.
\end{rem}

\begin{rem}\label{weakerstruct}
As mentioned earlier, if we are
  dealing with a particular weak solution $(u,\nabla u)$, then parts
  (ii), (iii) and (iv) of the structural assumptions \eqref{struct}
  required in Theorem \ref{main_thm} (where $\gamma=\delta =p$) can be
  weakened without affecting the conclusion. In particular, it is
  enough to assume that they hold when the general variable $z\in  \R$
  is replaced by $u(x), x\in \Omega$, i.e., to assume that for
  a.e. $x\in\Omega$ and all $(z,\xi) \in \R\times\R^n$,
\begin{equation}\label{relaxedstruct}
\begin{cases}
(i)\quad A(x,z,\xi)=\sqrt{Q(x)}{\tilde A}(x,z,\xi),\\
(ii)\quad\xi \cdot A(x,u(x),\xi) \geq a^{-1}\Big|\sqrt{Q(x)}\, \xi\Big|^p
- h(x)|u(x)|^p - g(x),\\
(iii)\quad\Big|\tilde{A}(x,u(x),\xi)\Big| \leq a\Big|\sqrt{Q(x)}\,
\xi\Big|^{p-1} + b(x)|u(x)|^{p-1} + e(x),\\
(iv)\quad\Big|B(x,u(x),\xi)\Big|\leq c\Big|\sqrt{Q(x)}\,
\xi\Big|^{\psi-1} +d(x)|u(x)|^{p-1}+f(x).
\end{cases}
\end{equation}
When $\gamma$ or $\delta$ exceeds $p$ and we assume the structural
conditions (\ref{struct}), this fact will be used in some of the
corollaries below to deduce boundedness results from the case
$\gamma=\delta=p$ considered in Theorem \ref{main_thm}; it will allow
us to bundle some powers of $|u(x)|$ together with the coefficients.
\end{rem}


We now turn to the question of estimating the expression
$\bar{Z}$ defined in \eqref{Zetabar}, and in particular of determining
when it is finite. 

In case $e=f=g=0$, we can let $k$ tend to $0$ in
\eqref{Linfinityestimate} to obtain \eqref{Linfinityestimate} for $u$
instead of $\bar{u}$. In our applications of Theorem \ref{main_thm},
provided $e,f,g$ are not all identically $0$ in $B_r=B(y,r)$, we will
choose the constant $k$ to be
\begin{eqnarray}
\nonumber \qquad\quad k(r)\!\!&\!\!=\!\!&\!\! k(y,r)\\
\label{k}\!\!&\!\!=\!\!&\!\!\left(r^{p-1}\avnormb{e}{p'\sigma'}{B_r}
\right)^{\frac{1}{p-1}}\!\! +\! \left(r^p
\avnormb{g}{\frac{p\sigma'}{p-\epsilon_2}}{B_r}\right)^{\frac{1}{p}}
\!\!+\!\left(r^p \avnormb{f}{\frac{p\sigma'}{p-\epsilon_3}}{B_r}
\right)^{\frac{1}{p-1}}\!\!.
\end{eqnarray}
As above, in case $k=0$ in (\ref{k}), then in order to be able to
apply Theorem \ref{main_thm}, we will instead choose any positive
number for $k$ and then let this number tend to $0$. In any case, it
follows from (\ref{k}) that the three terms of (\ref{Zetabar})
corresponding to ${\bar b}, {\bar h}$ and ${\bar d}$ satisfy
\begin{eqnarray}
\nonumber\!\!\! &\!\!\!\!&\!\!\!r^{p-1}
\avnormb{\bar{b}}{p'\sigma'}{B_r} +\left(r^p
\avnormb{\bar{h}}{\frac{p\sigma'}{p-\epsilon_2}}{B_r}
\right)^\frac{1}{\epsilon_2} + \left(r^p
\avnormb{\bar{d}}{\frac{p\sigma'}{p-\epsilon_3}}{B_r}
\right)^\frac{1}{\epsilon_3}\leq\\
\nonumber\!\!\! &\!\!\!\!&\!\!\!\hspace{0.6cm}
r^{p-1}\avnormb{b}{p'\sigma'}{B_r}+k^{1-p}r^{p-1}\avnormb{e}{p'\sigma'}{B_r}\\
\nonumber\!\!\! &\!\!\!\!&\!\!\!\hspace{1.3cm} + \left(r^p
\avnormb{h}{\frac{p\sigma'}{p-\epsilon_2}}{B_r} + k^{-p}
r^p\avnormb{g}{\frac{p\sigma'}{p-\epsilon_2}}{B_r}
\right)^\frac{1}{\epsilon_2}\\
\nonumber\!\!\! &\!\!\!\!&\!\!\!\hspace{2cm} +\left(r^p
\avnormb{d}{\frac{p\sigma'}{p-\epsilon_3}}{B_r} + k^{1-p}
r^p\avnormb{f}{\frac{p\sigma'}{p-\epsilon_3}}{B_r}
\right)^\frac{1}{\epsilon_3}\leq\\
\nonumber\!\!\! &\!\!\!\!&\,\,\,\,1+2^\frac{1}{\epsilon_2}\!
+2^\frac{1}{\epsilon_3}\! + r^{p-1} \avnormb{b}{p'\sigma'}{B_r} +\!
2^\frac{1}{\epsilon_2}\!\left(r^p
\avnormb{h}{\frac{p\sigma'}{p-\epsilon_2}}{B_r}
\right)^\frac{1}{\epsilon_2}\!\! +\!2^\frac{1}{\epsilon_3}\!
\left(r^p \avnormb{d}{\frac{p\sigma'}{p- \epsilon_3}}{B_r}
\right)^\frac{1}{\epsilon_3}\!\!\!.
\end{eqnarray}
Consequently, with $k$ defined by (\ref{k}), we may replace ${\bar
Z}$ in (\ref{Zetabar}) and (\ref{Linfinityestimate}) by the following
analogous expression in which ${\bar b}, {\bar h}, {\bar d}$ are
replaced respectively by $b, h, d$:
\begin{eqnarray} \label{Zeta}
\qquad Z &=& 1+r^{p-1}\avnormb{b}{p'\sigma'}{B_r}+ \left(r^p
\avnormb{c^\frac{p}{p+1-\psi}
\bar{u}^\frac{p(\psi-p)}{p+1-\psi}}{\frac{p\sigma'}{p-
    \epsilon_1}}{B_r}\right)^{\frac{1}{\epsilon_1}}\\
\nonumber&&\hspace{3,2cm} +\left(r^p
\avnormb{h}{\frac{p\sigma'}{p-\epsilon_2}}{B_r}
\right)^\frac{1}{\epsilon_2}+ \left(r^p
\avnormb{d}{\frac{p\sigma'}{p-\epsilon_3}}{B_r}
\right)^\frac{1}{\epsilon_3}.
\end{eqnarray}
Strictly speaking, the additive constant $1$ in (\ref{Zeta}) should be
replaced by $1+2^{1/\epsilon_2} + 2^{1/\epsilon_3}$, but we can
incorporate extra constant factors depending on $\epsilon_2,
\epsilon_3$ in the constant $C$ in (\ref{Linfinityestimate}).

In order to better understand the expression $Z$ in (\ref{Zeta}), we
first note that its form leads naturally to the following definition
of spaces of Morrey type for quasimetric balls.

\begin{defn}\label{Morrey} Let $\alpha, \beta$ satisfy $0<\alpha
  <\infty$ and $0<\beta \le \infty$. We say that a measurable function
  $m(x)$ on $\Omega$ belongs to the Morrey class
  $M^\beta_\alpha(\Omega)$ if
\begin{equation}\label{morrey}
\|m\|_{M^\beta_\alpha(\Omega)} =\sup\left\{ r^\alpha
\left(\frac{1}{|B(y,r)|}\int_{B(y,r)} |m(x)|^\beta
dx\right)^{\frac{1}{\beta}}\right\}=\sup\left\{r^\alpha
\avnormb{m}{\beta}{B_r}\right\} <\infty,
\end{equation}
where the sup is taken over all $\rho$-balls $B_r=B(y,r)$ with
$r<\min\{1,r_1(y)\}$. We recall that the closure of any such ball is
contained in $\Omega$. In case $\beta = \infty$, (\ref{morrey}) means
that
$$\|m\|_{M^\infty_\alpha(\Omega)}=\sup \{r^\alpha\, \text{ess
sup}_{B(y,r)} |m|\} <\infty.$$
\end{defn}

Using this notation, the expression $Z$ in \eqref{Zeta} satisfies
\begin{equation}\label{newZeta}
Z \le1+ ||b||_{M^{p'\sigma'}_{p-1}(\Omega)}\! + ||c^{\frac{p}{p+1
-\psi}}{\bar u}^{\frac{p(\psi-p)}{p+1-
\psi}}||_{M^{\frac{p\sigma'}{p-\epsilon_1}}_{p}(\Omega)
}^{\frac{1}{\epsilon_1}}\! +||h||_{M^{\frac{p\sigma'}{p-
\epsilon_2}}_{p}(\Omega)}^{\frac{1}{\epsilon_2}}\! +
||d||_{M^{\frac{p\sigma'}{p-\epsilon_3}}_{p}(\Omega)}^{\frac{1}{\epsilon_3}}\!.
\end{equation}
However, since $Z$ involves only a single ball, it is more local
than the right-hand side of (\ref{newZeta}), and we will often take
further advantage of its local nature before using Morrey classes.

In general, there is no simple way to characterize Morrey classes in
terms of Lebesgue classes. However, it is possible to combine a
Lebesgue condition with an (algebaic) growth condition on $|B_r|$ in
order to estimate the size of $r^\alpha \avnormb{m}{\beta}{B_r}$ and
determine upper bounds for $Z$. To do this, we will use the
following simple observations.

For balls as in Definition \ref{Morrey}, note that
\[
r^\alpha \left(\frac{1}{|B(y,r)|}\int_{B(y,r)} |m|^\beta dx
\right)^{\frac{1}{\beta}} \le \left(\sup
\frac{r^\alpha}{|B(y,r)|^{1/\beta}}\right)
\,\Lnormb{m}{\beta}{B(y,r)},
\]
where the supremum is taken over all such balls.

\begin{defn} \label{weakdoub} If $q^*$ satisfies $0<q^* <\infty$ and
  there is a positive constant $c_0$ such that
\begin{equation}\label{D*}
|B(y,r)| \ge c_0 r^{q^*}
\end{equation}
for all $\rho$-balls $B(y,r)$ with $r<\min\{1,r_1(y)\}$, we will say
that condition $D_{q^*}$ holds.
\end{defn}

Condition \eqref{D*} is related to, but weaker than, the local doubling
condition
\begin{eqnarray}\label{predoubling}
|B(x,2r)| \le C|B(x,r)|,\quad x\in \Omega, 0<r<\bar{r}(x),
\end{eqnarray}
where $\bar{r}(x)<r_0(x)/2$. It is well-known that \eqref{predoubling}
implies there are positive constants $C, D^*$ such that
\begin{equation}\label{doubling}
|B(x,R)| \le C\left( \frac{R}{r}\right)^{D^*} |B(y,r)| \quad\mbox{if
$B(y,r) \subset B(x,R)$}
\end{equation}
and $r$ is sufficiently small. Note that \eqref{D*} follows from
(\ref{doubling}) with $D^*=q^*$ if $\rho$ is bounded in $\Omega$ since
then by choosing $R= \sup\{\rho(x,y): x,y\in \Omega\}$, we have for
all $x\in \Omega$ that $\Omega \subset B(x,R)$, and consequently
$B(x,R) = \Omega$. Moreover, even if $\rho$ is unbounded in $\Omega$,
\eqref{D*} follows from \eqref{doubling} if $\inf\{|B(x,1)|: x\in
\Omega\} >0$. 

Then for $\alpha, \beta$ as above,
\begin{equation} \label{D*plusLEB}
r^\alpha \left(\frac{1}{|B(y,r)|}\int_{B(y,r)} |m|^\beta dx
\right)^{\frac{1}{\beta}} \le
\begin{cases}
C \Lnormb{m}{\beta}{B(y,r)}\quad
  \text{if $\beta <\infty$ and $D_{\alpha\beta}$ holds} \\
||m||_{L^\infty(B(y,r))}
\end{cases}
\end{equation}
with $C= c_0^{-1/\beta}$, where to obtain the second option we have used
$\alpha>0$ and $r<1$. Thus $L^\beta(\Omega) \subset
M^\beta_\alpha(\Omega)$ if $\beta=  \infty$, or if $\beta<\infty$ and
$D_{\alpha\beta}$ holds.

If $m$ is a product, $m(x) = m_1(x)m_2(x)$, and if $\beta_1, \beta_2$
satisfy $0< \beta_1, \beta_2 \le \infty$ and $\frac{1}{\beta} =
\frac{1}{\beta_1} + \frac{1}{\beta_2}$, then H\"older's inequality
implies that for any $B_r= B(y,r)$,
\[
r^\alpha \left(\frac{1}{|B_r|}\int_{B_r} |m_1m_2|^\beta dx
\right)^{\frac{1}{\beta}} \le r^\alpha \left(\frac{1}{|B_r|}\int_{B_r}
|m_1|^{\beta_1} dx \right)^{\frac{1}{\beta_1}}  \left(\frac{1}{|B_r|}
\int_{B_r} |m_2|^{\beta_2} dx \right)^{\frac{1}{\beta_2}}.
\]
Combining this with (\ref{D*plusLEB}) gives (again we denote
$B_r=B(y,r)$)
\begin{equation}\label{product}
r^\alpha \left(\frac{1}{|B_r|}\int_{B_r} |m_1m_2|^\beta dx
\right)^{\frac{1}{\beta}} \le C\Lnormb{m_1}{\beta_1}{B(y,r)}
    \left(\frac{1}{|B_r|} \int_{B_r} |m_2|^{\beta_2}
    dx\right)^{\frac{1}{\beta_2}}\quad \text{provided}
\end{equation}
\[
\frac{1}{\beta} = \frac{1}{\beta_1} + \frac{1}{\beta_2}, \text{  and
either $\beta_1=\infty$, or $\beta_1<\infty$ and $D_{\alpha\beta_1}$
holds}.
\]

Estimates (\ref{D*plusLEB}) and (\ref{product}) serve as a basis for
finding different ways to majorize (\ref{Zeta}) by using $D_{q^*}$
conditions and Lebesgue classes, and lead to the corollaries below.
We emphasize that our corollaries cover only a few special cases and
do not exhaust all possibilities. Recall from (\ref{ranges}) that we
always assume $\gamma, \psi, \delta$ satisfy
\[
\gamma \in (1,\sigma(p-1)+1),\quad \psi\in (1,p+1-\sigma^{-1}), \quad
\delta \in (1,p\sigma).
\]
The fewest technicalities arise when $\gamma= \delta = \psi = p$, and
we begin with that case. The result we will state aims at making the
weakest possible integrability assumptions on the coefficients, and
consequently it makes a strong assumption about the order of the $D$
condition. As is apparent from (\ref{D*}) and (\ref{D*plusLEB}),
since $r<1$ and $\alpha\beta$ increases with $\beta$, a general
principle is that the better the coefficients are (i.e., the higher
their integrability becomes), then the weaker the required $D$ condition
becomes.

\begin{cor}\label{allp} Suppose the same hypotheses and notation as in
  Theorem \ref{main_thm} hold, but now also that $\psi =p$ (i.e.,
  $\gamma=   \delta = \psi =p$), that condition $D_{q^*}$ holds for
  some $q^*\le p\sigma'$, and that $b, e \in L^{p'\sigma'}(B(y,r))$, $c\in
  L^{p\sigma'(1+\epsilon)}(B(y,r))$ and $d, f, h, g \in L^{\sigma'(1 +
  \epsilon)}(B(y,r))$ for some $\epsilon >0$. Then
\begin{equation}\label{Linfinityestimate1}
\|u\|_{L^\infty(B(y,\tau r))}\leq
C\left\{\left(\frac{1}{|B(y,r)|}\int_{B(y,r)} |u|^{sp} dx
\right)^\frac{1}{sp} + K(y,r)\right\},\quad\text{where}
\end{equation}
\begin{eqnarray}
\nonumber
K(y,r)\!\!&\!\!=\!\!&\!\!r^{1-\frac{q^*}{p\sigma'}}
\Big(\Lnormb{e}{p'\sigma'}{B(y,r)}^{\frac{1}{p-1}}
+r^{\{1-\frac{q^*}{p\sigma'} +
  \frac{q^*}{\sigma'}\frac{\epsilon}{1+\epsilon}\}\frac{1}{p-1}}
\Lnormb{f}{\sigma'(1+\epsilon)}{B(y,r)}^{\frac{1}{p-1}}\\
\nonumber\!\!&\!\!\!\!&\!\!\hspace{7cm}+
r^{\frac{q^*}{p\sigma'}\frac{\epsilon}{1+\epsilon}}
\Lnormb{g}{\sigma'(1+\epsilon)}{B(y,r)}^{\frac{1}{p}}\Big)
\end{eqnarray}
and $C$ depends on all relevant parameters including $\epsilon$, the
constant in the $D_{q^*}$ condition and the sum of the corresponding
norms of $b, c, d, h$ over $B(y,r)$, but does not depend on $u,
B(y,r), e, f$ or $g$. In particular, if $s=1$, i.e., if the cutoff condition
(\ref{cutoff}) holds in the $L^\infty$ sense, then
\begin{equation}\label{Linfinityestimate1*}
\|u\|_{L^\infty(B(y,\tau r))}\leq C\left\{\left(
\frac{1}{|B(y,r)|}\int_{B(y,r)} |u|^p dx\right)^\frac{1}{p}+ K(y,r)\right\}
\end{equation}
with $K(y,r)$ and $C$ as above.
\end{cor}

Note that for the case of the standard Euclidean structure,
$p\sigma' =n$ and condition $D_n$ automatically holds. Hence,
since (\ref{cutoff}) is true with $s^*=\infty$ in this situation,
estimate (\ref{Linfinityestimate1*}) then applies and includes the
local boundedness result of \cite[Theorem 1, p. 555]{S}.

If $e,f,g$ vanish identically in $B(y,r)$, then $K(y,r)=0$ in
Corollary \ref{allp}. Also, if $q^*<p\sigma'$ and the corresponding
norms of $e,f,g$ over all of $\Omega$ are finite, note that $K(y,r) \le
c r^\eta$ for some $\eta >0$ which depends on $q^*$.

The proof of Corollary \ref{allp} is an application of Theorem
\ref{main_thm} and follows from \eqref{Zeta} and \eqref{D*plusLEB},
without needing to use \eqref{product}. We choose $k$ as in
\eqref{k} and drop $k$ on the left side of
\eqref{Linfinityestimate}, thereby replacing ${\bar u}$ by $u$ on
the left side. However, on the right side, we use the hypotheses to
estimate $Z$ and show that $k\le CK(y,r)$. As examples of the
required computations, let us briefly indicate how to estimate the
term of $Z$ which corresponds to $b$ and the term of $k$ which
corresponds to $g$. Denoting $B(y,r) =B$ and using the estimate
$|B|\ge c_0r^{q^*}$, we have
\begin{eqnarray}
\nonumber r^{p-1} \avnormb{b}{p'\sigma'}{B}\!&\!=\!&\!
\left(\frac{r}{|B|^{1/(p\sigma')}}\right)^{p-1}\Lnormb{b}{p'\sigma'}{B}
\,\,\,\le\,\,\,\left(c_0^{-\frac{1}{p\sigma'}}
r^{1-\frac{q^*}{p\sigma'}}\right)^{p-1}\Lnormb{b}{p\sigma'}{B}\\
\nonumber\!&\!\le\!&\!
c_0^{-\frac{1}{p'\sigma'}}\Lnormb{b}{p\sigma'}{B}
\end{eqnarray}
since $q^*\le p\sigma'$ and $r\le 1$. Similarly, choosing
$\epsilon_2 = \epsilon p/(1+\epsilon)$, we obtain $p-\epsilon_2 =
p/(1+\epsilon)$ and
\[
\left(r^p\avnormb{g}{\frac{p\sigma'}{p-\epsilon_2}}{B}\right)^{\frac{1}{p}}
\le\left(c_0^{-\frac{1}{p\sigma'(1+\epsilon)}}r^{\{1-
  \frac{q^*}{p\sigma'}\}+\frac{q^*}{p\sigma'}
  \frac{\epsilon}{1+\epsilon}}\right)
\Lnormb{g}{\sigma'(1+\epsilon)}{B}^{\frac{1}{p}}.
\]
Further details are left to the reader.

Our next corollary gives an estimate when all of $\gamma, \psi,
\delta$ are less than $p$. In this case, we can easily replace each of
the structural assumptions (\ref{struct})(ii), (iii) and (iv) by a
similar one involving only $p$ and modified coefficients. For example,
if $\gamma <p$, we can use the simple estimate
\[
b|z|^{\gamma-1}+e \le b|z|^{p-1} + (b+e)
\]
to see that (\ref{struct})(iii) implies
\[
\Big|\tilde{A}(x,z,\xi)\Big| \leq a\Big|\sqrt{Q(x)}\,
\xi\Big|^{p-1} + b|z|^{p -1} + (b+e).
\]
Similarly, an analogue of (\ref{struct})(ii) holds with
$-h|z|^\gamma -g$ replaced by $-h|z|^p - (h+g)$, and if $\psi, \delta
<p$ then (\ref{struct})(iv) gives
\[
\Big|B(x,z,\xi)\Big|\leq c\Big|\sqrt{Q(x)}\, \xi
\Big|^{p-1}+d|z|^{p-1}+ (f+c+d).
\]
It follows that when $\gamma, \psi, \delta$ are less than $p$,
(\ref{struct}) implies its analogue with $\gamma, \psi, \delta$ all
replaced by $p$ and with $e, g, f$ replaced by $e+b, g+h, f+c+d$
respectively. Hence we immediately obtain the next corollary from the
previous one.

\begin{cor}\label{all<p} Suppose the same hypotheses and
  notation as in Theorem \ref{main_thm} hold with these exceptions:
  the structural assumptions (\ref{struct}) hold for some $\gamma,
  \psi, \delta <p$; condition $D_{q^*}$ holds for some $q^*\le p\sigma'$;
  and $b, e \in  L^{p'\sigma'}(B(y,r))$, $c\in L^{p\sigma'(1+
  \epsilon)}(B(y,r))$ and $d, f, h, g \in L^{\sigma'(1+
  \epsilon)}(B(y,r))$ for some $\epsilon >0$. Then
  (\ref{Linfinityestimate1}) is true with
\begin{eqnarray}
\nonumber
K(y,r)\!\!&\!\!=\!\!&\!\!r^{1-\frac{q^*}{p\sigma'}}\Big(\Lnormb{e+
  b}{p'\sigma'}{B(y,r)}^{\frac{1}{p-1}}
+r^{\frac{q^*}{p\sigma'}\frac{\epsilon}{1+\epsilon}}
\Lnormb{g+h}{\sigma'(1+\epsilon)}{B(y,r)}^{\frac{1}{p}}\\
\nonumber\!\!&\!\!\!\!&\!\!\hspace{3.5cm}+
r^{\{1-\frac{q^*}{p\sigma'}+\frac{q^*}{\sigma'}\frac{\epsilon}{1
    +\epsilon}\}\frac{1}{p-1}}
\Lnormb{f+c+d}{\sigma'(1+\epsilon)}{B(y,r)}^{\frac{1}{p-1}}\Big)
\end{eqnarray}
and with $C$ depending on all relevant parameters including
$\epsilon$, the constant in the $D_{q^*}$ condition and the norms of
$b, c, d, h$ over $B(y,r)$, but not on $u, B(y,r), e, f$ or $g$.
\end{cor}

Next we list corollaries in case all of $\gamma, \psi, \delta$ in
(\ref{struct}) exceed $p$. In this situation, we will use the
observation in Remark \ref{weakerstruct}; in fact, the last three
assumptions in either (\ref{struct}) or its weaker analogue when $z$
is replaced by $u(x)$ (for a fixed weak solution $u$) yield
\begin{equation}\label{bundledstruct}
\begin{cases}
(i)\quad A(x,z,\xi)=\sqrt{Q(x)}{\tilde A}(x,z,\xi),\\
(ii)\quad\xi \cdot A(x,u(x),\xi) \geq a^{-1}\Big|\sqrt{Q(x)}\, \xi\Big|^p
- \Big(h(x)|u(x)|^{\gamma-p}\Big) |u(x)|^p - g(x),\\
(iii)\quad\Big|\tilde{A}(x,u(x),\xi)\Big| \leq a\Big|\sqrt{Q(x)}\,
\xi\Big|^{p-1} + \Big(b(x)|u(x)|^{\gamma-p}\Big)|u(x)|^{p-1} + e(x),\\
(iv)\quad\Big|B(x,u(x),\xi)\Big|\leq c\Big|\sqrt{Q(x)}\,
\xi\Big|^{\psi-1} +\Big(d(x)|u(x)|^{\delta-p}\Big)|u(x)|^{p-1}+f(x)
\end{cases}
\end{equation}
for the same $\gamma, \psi, \delta$ and $a, b, c, d, e, f, g, h$ as in
(\ref{struct}). Consequently, denoting
\[
b^*= b|u|^{\gamma-p},\quad d^* = d|u|^{\delta-p},\quad  h^* =
h|u|^{\gamma-p},
\]
we obtain the conditions (\ref{relaxedstruct}) with $b, d, h$ there
replaced respectively by $b^*, d^*, h^*$. Using Remark
\ref{weakerstruct} as well as (\ref{Zeta}), where the form of
the constant $k$ in the formula $\bar{u} = |u|+k$ is still the same
as in (\ref{k}), we may apply Theorem \ref{main_thm} with ${\bar Z}$
replaced by
\begin{eqnarray}
\nonumber Z^*\!&\!=\!&\!
1+r^{p-1}\avnormb{b|u|^{\gamma-p}}{p'\sigma'}{B_r} + \left(r^p
\avnormb{c^\frac{p}{p+1-\psi}
\bar{u}^\frac{p(\psi-p)}{p+1-\psi}}{\frac{p\sigma'}{p-
    \epsilon_1}}{B_r}\right)^\frac{1}{\epsilon_1}\\
\label{Zeta*}\!&\!\!&\!\hspace{0.5cm} +\left(r^p
\avnormb{h|u|^{\gamma-p}}{\frac{p\sigma'}{p-\epsilon_2}}{B_r}
\right)^\frac{1}{\epsilon_2} + \left(r^p
\avnormb{d|u|^{\delta-p}}{\frac{p\sigma'}{p-\epsilon_3}}{B_r}
\right)^\frac{1}{\epsilon_3}.
\end{eqnarray}

The terms in (\ref{Zeta*}) can be treated by applying (\ref{product}),
and we obtain the following corollaries. In the first one, we make the
strongest possible assumption on the coefficients, namely that they
are all bounded. In this case, we require no $D$ condition at
all. By using (\ref{Zeta*}) together with (\ref{product}) for $\beta_1=
\infty$, we obtain
\begin{eqnarray}
\nonumber \qquad Z^*\!\!&\!\!\le\!\!&\!\!1+
\|b\|_{L^\infty(B_r)}\avnormb{|u|^{\gamma-p}}{p'\sigma'}{B_r} +
\left(\|c\|_{L^\infty(B_r)}^\frac{p}{p+1-\psi}
\avnormb{\bar{u}^\frac{p(\psi-p)}{p+1-\psi}}{\frac{p\sigma'}{p-
    \epsilon_1}}{B_r}\right)^\frac{1}{\epsilon_1}\\
\label{Zeta*strong}\!&\!\!&\!+\left(\|h\|_{L^\infty(B_r)}
\avnormb{|u|^{\gamma-p}}{\frac{p\sigma'}{p-\epsilon_2}}{B_r}
\right)^\frac{1}{\epsilon_2}\!\!
+\! \left(\|d\|_{L^\infty(B_r)}
\avnormb{|u|^{\delta-p}}{\frac{p\sigma'}{p-\epsilon_3}}{B_r}
\right)^\frac{1}{\epsilon_3}\!\!.
\end{eqnarray}

In the third term on the right side of the estimate for $Z^*$, we use
\eqref{k} and the fact that $r\le 1$ to obtain
$$
\bar{u}= |u| + k \le |u| +
\left(\avnormb{e}{p'\sigma'}{B_r}\right)^{\frac{1}{p-1}} +
\left(\avnormb{g}{\frac{p\sigma'}{p-\epsilon_2}}{B_r}\right)^{\frac{1}{p}}
+ \left( \avnormb{f}{\frac{p\sigma'}{p-\epsilon_3}}{B_r}
\right)^{\frac{1}{p-1}}.
$$
Choosing $\epsilon_1, \epsilon_2, \epsilon_3$ small, we
then easily obtain from (\ref{Zeta*strong}) that for any sufficiently
small $\epsilon>0$, depending on $\gamma, \delta, \psi, p$ and $\sigma$,
there are constants $\theta, C_1, C_2, L$ satisfying
\begin{equation}\label{theta}
\theta = \max\left\{(\gamma-p)p'\sigma',
 \frac{(\psi-p)p\sigma'(1+\epsilon)}{p+1-\psi},
 (\delta-p)\sigma'(1+\epsilon) \right\},
\end{equation}
\begin{equation}\label{C1}
C_1 = C_1\big(\epsilon,p,\psi, \|c\|_{L^\infty(B(y,r))},
||e||_{L^\infty(B(y,r))}, ||f||_{L^\infty(B(y,r))},
||g||_{L^\infty(B(y, r))}\big)
\end{equation}
\[
\text{with } C_1 =0 \quad\text{when $c\equiv 0$ in $B(y,r)$ or when
  $e,  f, g \equiv 0$ in $B(y,r)$},
\]
\begin{equation}\label{C2}
C_2= C_2\big(\epsilon, p, \psi, ||b||_{L^\infty(B(y,r))},
||c||_{L^\infty(B(y,r))},
||d||_{L^\infty(B(y,r))}, ||h||_{L^\infty(B(y,r))}\big)
\end{equation}
\[
\text{with } C_2 = 0 \quad\text{when $b, c, d, h \equiv 0$ in $B(y,r)$,}
\]
and
\begin{equation}\label{L}
L = L(\epsilon, p, \gamma, \psi, \delta)
\end{equation}
 such that
$$
Z^* \le 1+ C_1+ C_2\left[1 + \left(\frac{1}{|B(y,r)|}
  \int_{B(y,r)} |u|^\theta dx \right)^\frac{1}{\theta}\right]^L.
$$
Moreover, for small $\epsilon$, the restrictions (\ref{ranges}) imply
that $\theta <p\sigma$, and consequently that $u\in L^\theta(B(y,r))$.
Thus we obtain the following estimate.

\begin{cor} \label{first>p} Suppose the same hypotheses and
notation as in Theorem \ref{main_thm} hold with these exceptions:
the structural assumptions (\ref{struct}) hold for some $\gamma,
\psi, \delta >p$ which satisfy (\ref{ranges}), and the coefficients
$b, c, d, e, f, g, h$ are bounded in $\Omega$. For small
$\epsilon>0$, define $\theta, C_1, C_2$ and $L$ as in (\ref{theta}),
(\ref{C1}), (\ref{C2}) and (\ref{L}), respectively. Then for any
$\rho$-ball $B(y,r)$ with $0<r< \tau r_1(y)$,
\begin{equation}\label{Linfinityestimate2}
\begin{array}{l}
\displaystyle\|u\|_{L^\infty(B(y,\tau r)} \le C \left\{1+ C_1+
C_2\left[ 1
  + \left(\frac{1}{|B(y,r)|} \int_{B(y,r)} |u|^\theta dx
  \right)^\frac{1}{\theta}\right]^L  \right\}^{\Psi_0} \times\\
\displaystyle\hspace{5cm} \left\{ \left(\frac{1}{|B(y,r)|}
\int_{B(y,r)} |u|^{sp} dx\right)^\frac{1}{sp} + K(y,r)\right\},
\end{array}
\end{equation}
\[\text{where}\quad
K(y,r) = \left(r\|e\|_{L^\infty(B(y,
    r))}^{\frac{1}{p-1}} + r^{p'}\|f\|_{L^\infty(B(y,
    r))}^{\frac{1}{p-1}} +r \|g\|_{L^\infty(B(y,
    r))}^{\frac{1}{p}}\right)
\]
and $C$ is as in (\ref{Linfinityestimate}). In particular, if
    (\ref{cutoff}) holds in the $L^\infty$ sense, then
(\ref{Linfinityestimate2}) holds with $s=1$. In this corollary, no
    $D$ condition is needed.
\end{cor}

As noted above, when $\epsilon$ is small, the value of $\theta$ in
(\ref{Linfinityestimate2}) satisfies $\theta <p\sigma$. The largest
power of $|u|$ which is a priori locally integrable is $p\sigma$,
and our next corollary gives an estimate when $\theta$ is replaced by
$p\sigma$, still asuming that all of $\gamma, \psi, \delta$ exceed
$p$. In this situation, the conditions required of the coefficients are
weaker than boundedness, but an appropriate restriction in terms of a
$D$ condition is required.

\begin{cor}\label{second>p} Suppose the same hypotheses and
notation as in Theorem \ref{main_thm} hold with these exceptions:
the structural assumptions (\ref{struct}) hold for some $\gamma,
\psi, \delta >p$ which satisfy (\ref{ranges});  the $D_{q^*}$
condition holds for some $q^*>0$ as described below; for a given
$\rho$-ball $B(y,r)$ with $r<\tau r_1(y)$, the coefficients satisfy
\[
b\in L^{\mathcal B}(B(y,r)),\, {\mathcal B}\ge
  \frac{p\sigma}{\sigma(p-1)+1- \gamma};  \quad c\in  L^{\mathcal
  C}(B(y,r)),\, {\mathcal C} > \frac{p\sigma}{\sigma(p +1-\psi)-1};
\]
\[
d \in L^{\mathcal D}(B(y,r)),\, {\mathcal D}> \frac{p\sigma}{p\sigma-
\delta}; \quad e  \in L^{\mathcal
  E}(B(y,r)),\, {\mathcal E}\ge p'\sigma';
\]
\[
 f \in L^{\mathcal F}(B(y,r)),\, {\mathcal F}>\sigma';\quad  g
 \in L^{\mathcal G}(B(y,r)),\, {\mathcal G} > \sigma'; \quad h\in
 L^{\mathcal H}(B(y,r)),\, {\mathcal H}
>\frac{p\sigma}{p\sigma- \gamma};
\] and
\[
q^*\le \min \big\{{\mathcal B}(p-1),{\mathcal E}(p-1),{\mathcal
  C}(p+1- \psi),   {\mathcal D}p, {\mathcal F}p,
  {\mathcal   G}p, {\mathcal H}p\big\}.
\]
Then
\begin{equation}\label{Linfinityestimate3}
\begin{array}{l}
\displaystyle\|u\|_{L^\infty(B(y,\tau r))}\!\le\! C\! \left\{1+C_1'
+
  C_2' \left[ 1 + \left(\frac{1}{|B(y,r)|} \int_{B(y,r)}
  |u|^{p\sigma} dx \right)^\frac{1}{p\sigma}\right]^L
  \right\}^{\Psi_0}\!\!\!\times\\
\displaystyle\hspace{5,5cm}\left\{\left(\frac{1}{|B(y,r)|}
\int_{B(y,r)}
  |u|^{sp} dx \right)^\frac{1}{sp} + K(y,r)\right\}
\end{array}
\end{equation}
where
\[
K(y,r) = \left(r^{1-\frac{q^*}{(p-1)\mathcal{E}}}\Lnormb{e}{\mathcal
E}{B(y,r)}^{\frac{1}{p-1}} +
r^{p'\{1-\frac{q^*}{p\mathcal{F}}\}}\Lnormb{f}{\mathcal
F}{B(y,r)}^{\frac{1}{p-1}} +r^{1-\frac{q^*}{p\mathcal{G}}}
\Lnormb{g}{\mathcal G}{B(y,r)}^{\frac{1}{p}}\right),
\]
$C$ is as in (\ref{Linfinityestimate}), $L>0$ is a constant
depending on $\gamma, \psi, \delta, p, \sigma, {\mathcal B},
{\mathcal C}, {\mathcal D}, {\mathcal H}$, and the constants $C_1',
C_2'$ are analogous to the constants $C_1, C_2$ in (\ref{C1}),
(\ref{C2}) but with the $L^\infty(\Omega)$ norms of $b, c, d, e, f, g,
h$ replaced by their corresponding norms listed above. Also, $C_1',
C_2'$ have the same vanishing properties as $C_1, C_2$ but now depend
on the constants in the $D_{q^*}$ condition as well as $\gamma, \psi,
\delta, p, \sigma, {\mathcal C}, {\mathcal D}, {\mathcal
H}$ and the norm of $c$; $C_1'$ depends furthermore on the norms of
$e,f,g$, and $C_2'$ on ${\mathcal B}$ and the norms of $b,d,h$.
Again, in case (\ref{cutoff}) holds in the $L^\infty$ sense,
then (\ref{Linfinityestimate3}) holds with $s=1$.
\end{cor}

The proof of Corollary \ref{second>p} is computational but
straightforward; it is based on (\ref{product}) and
(\ref{D*plusLEB}). Details are left to the reader. Further
computations show that if $q^* = p\sigma'$, then all conditions
involving $q^*$ in Corollary \ref{second>p} are satisfied.  This is
the case in the classical Euclidean setting if $p<n$ since then
$\sigma' = n/p$, giving $p\sigma' = n$.

\section{\textbf{Preliminary Definitions and
Lemmas}}\label{section2}

Let $\Omega\subset\rn$ be a bounded open set and  $|Q|\in
L_\text{loc}^\frac{p}{2}(\Omega)$. For $w\in\lipl$, recall that
\bea
 \label{W^1p_Qnorm}\|w\|_{\waQ}& =&\left(\int_\Omega|w|^p\,dx +
 \int_\Omega\langle\nabla w,Q\nabla w\rangle^\frac{p}{2}\,
   dx\right)^\frac{1}{p}\\
 \nonumber&=&\left(\int_\Omega|w|^p\,dx+\int_\Omega|\sqrt{Q}\nabla w|^p\,
   dx\right)^\frac{1}{p}\\
 \nonumber&=&\left(\Lnormb{w}{p}{\Omega}^p+\Lnormb{\,|\sqrt{Q}\nabla
   w|\,}{p}{\Omega}^p\right)^\frac{1}{p}.
\eea

By definition, $\waQ$ is the completion with respect to
$\|\cdot\|_{\waQ}$ of those functions in $\lipl$ with finite
$\waQ$-norm. Also, $\wbQ$ denotes the completion with respect to
$\|\cdot\|_{\waQ}$ of $\lipo$.

A sequence $\{w_i\}_{i\in\N}\subset\lipl$ such that
$\|w_i\|_{\waQ}<\infty$ for every $i$ and which is Cauchy with
respect to  $\|\cdot\|_{\waQ}$ identifies an element of $\waQ$. Then
$\{w_i\}$ is a Cauchy sequence in $\lp$ and $\{\sqrt{Q}\nabla
w_i\}$ is a Cauchy sequence in $[\lp]^n$. Hence, up to subsequences,
as $i\rightarrow\infty$,
\bea
 \label{Lp-part}w_i\!&\!\longrightarrow\!&\!w\qquad\qquad\qquad\text{
   in }\lp \text{  and a.e. in }\Omega,\\
 \label{scriptL-part}\nabla w_i\!&\!\longrightarrow\!&\!\mathbf{v}:=
\nabla  w\qquad\quad\text{ in   }\scl^p(\Omega)\text{ and}\\
 \label{grad-part}\sqrt{Q}\nabla w_i\!&\!\longrightarrow\!&\!
 \sqrt{Q}\nabla w\,\,\,\,\,\,\text{ in
  }[\lp]^n\text{ and a.e. in }\Omega.
\eea

We adopt the abuses of notation mentioned in the Introduction. Thus,
we will not generally distinguish between $W^{1,p}_Q(\Omega)$ and its
isomorphic copy $\mathcal{W}^{1,p}_Q(\Omega)$ defined to be the
collection of pairs $(w, \nabla w)$ which arise as in \eqref{Lp-part},
\eqref{scriptL-part}, \eqref{grad-part}. We will often write
simply $w$ instead of $(w, \nabla w)$ even though $\nabla w$ may not
be uniquely determined by $w$.

It follows from \eqref{Lp-part}, \eqref{scriptL-part} and
\eqref{grad-part} that \eqref{W^1p_Qnorm} also holds for a generic
element $w\in\waQ$. Similarly, it follows by passing to the limit that
conditions like \eqref{3.6-0} hold for all functions in
$W^{1,p}_Q(\Omega)$ instead of just for $Lip_{loc}(\Omega)$. In order
to deal with the left side of such inequalities when passing to the
limit,  we generally use Fatou's lemma.

The role of condition \eqref{Cond1} is illustrated in the next simple
lemma.

\begin{lem}\label{L0} If \eqref{Cond1} holds, then for every
$y\in\Omega$ there exists $r_0=r_0(y)>0$ such that
$\overline{B(y,r)}\subset\Omega$ for all $r\in(0,r_0]$.
\end{lem}

\textbf{Proof:} Let $y\in \Omega$. Since $\Omega$ is open, there
exists $\epsilon >0$ such that $D(y,\epsilon) \subset \Omega$. By
\eqref{Cond1.1}, there is $r_0>0$ for which $B(y,r_0) \subset D(y,
\epsilon/2)$, and it follows that the closure of $B(y,r_0)$ lies in
$\Omega$. $\qquad\Box$

We now derive a useful version of the product rule. See the comments
after Lemma \ref{L1} for a more global version.

\begin{pro}\label{C4}
Let \eqref{Cond1} be true. Suppose that \eqref{E3} holds for a
particular $\rho$-ball $B$ with closure in $\Omega$ and a particular
function $\eta \in Lip_0(B)$. Suppose also that, for $t\geq1$ as in
\eqref{E3}, condition \eqref{3.6-0} holds for $B$ with $t'$ given by
$1/t+1/t'=1$. If $\theta\geq1$ and $w\in W^{1,p}_Q(\Omega),$ then
$\eta^\theta w\in (W^{1,p}_{Q})_0(B)$ and
$$\sqrt{Q}\nabla(\eta^\theta w)=\theta\eta^{\theta-1}w\sqrt{Q}
\nabla\eta+\eta^\theta\sqrt{Q}\nabla w\qquad\text{a.e. in }\Omega.$$
\end{pro}

\textbf{Proof:}  Since $w\in W^{1,p}_Q(\Omega)$, there is a
sequence $\{w_i\}\subset Lip_{loc}(\Omega)$ representing $w$ which is
Cauchy in $W^{1,p}_Q(\Omega)$. Taking a subsequence, we may assume
that \eqref{Lp-part}, \eqref{scriptL-part} and \eqref{grad-part} hold.
Now fix $B$ and $\eta$ as in the hypotheses, and consider the sequence
$\varphi_i=\eta^\theta w_i$. Clearly,
$\varphi_i\in Lip_0(B)$.

\textit{Claim 1}: $\varphi_i\longrightarrow\varphi:=\eta^\theta w$
a.e. in $\Omega$ and in $L^p(\Omega)$.\\
In fact,  $\varphi_i\longrightarrow\eta^\theta w$ a.e. in $\Omega$ by
our assumptions on the sequence $\{w_i\}$. Also,
$$|\varphi_i-\eta^\theta w|^p= |\eta|^{\theta p}|w_i-w|^p\leq
C(\theta,\eta)^p\,|w_i-w|^p,$$ and hence
$\Lnormb{\varphi_i-\eta^\theta w}{p}{\Omega}\leq
C(\theta,\eta)\Lnormb{w_i-w}{p}{\Omega}\longrightarrow0$, which
proves the claim.

\textit{Claim 2:} For a.e. $x\in \Omega$ and in $[L^p(\Omega)]^n$ norm,
\begin{equation}\label{50}
\sqrt{Q}\nabla\varphi_i\longrightarrow
\theta\eta^{\theta-1}w\sqrt{Q}\nabla\eta+\eta^\theta\sqrt{Q}\nabla
w.
\end{equation}
By the product rule for Lipschitz functions, $\sqrt{Q}\nabla\varphi_i
=\theta\eta^{\theta-1} w_i\sqrt{Q}\nabla\eta+\eta^\theta
\sqrt{Q}\nabla w_i$ a.e. in $\Omega$. Then \eqref{50} holds for
a.e. $x\in \Omega$ by the convergence properties of $w_i$ and
$\sqrt{Q} \nabla w_i$. Moreover
\begin{equation*}
\begin{array}{l}
 \displaystyle\int_\Omega\left|\sqrt{Q}\nabla\varphi_i-\theta
    \eta^{\theta-1}w\sqrt{Q}\nabla\eta-\eta^\theta\sqrt{Q}\nabla
    w\right|^pdx\\
 \displaystyle\hspace{2cm}\leq2^p\int_\Omega \left[\theta
    |\eta|^{(\theta-1)p}|w_i-w|^p|\sqrt{Q}\nabla\eta|^p+|\eta|^{\theta
    p}|\sqrt{Q}\nabla w_i-\sqrt{Q}\nabla w|^p\right]\,dx\\
 \displaystyle\hspace{2cm}\leq C(\theta,\eta)^p \int_\Omega \left[|w_i
    -w|^p|\sqrt{Q}\nabla\eta|^p+|\sqrt{Q}\nabla w_i-\sqrt{Q}\nabla
    w|^p\right]\,dx.
\end{array}
\end{equation*}
Now $$\int_\Omega |\sqrt{Q}\nabla w_i-\sqrt{Q}\nabla
    w|^p\,dx=\Lnormb{\sqrt{Q}\nabla w_i-\sqrt{Q}\nabla
    w}{p}{\Omega}^p\leq\|w_i-w\|_{W^{1,p}_Q(\Omega)}^p\longrightarrow0.$$
By H\"older's inequality and \eqref{3.6-0} (recall that \eqref{3.6-0}
    holds for general elements of $W^{1,p}_Q(\Omega)$),
\begin{eqnarray}
 \nonumber\int_\Omega |w_i-w|^p|\sqrt{Q}\nabla\eta|^p\,dx\!\!&\!\!=
 \!\!&\!\! \int_B |w_i-w|^p|\sqrt{Q}\nabla\eta|^p\,dx\\
\label{E4}\!\!&\!\!\leq\!\!&\!\!
     \left(\int_{B}|w_i-w|^{pt'}dx\right)^\frac{1}{t'}
     \left(\int_{B}|\sqrt{Q}\nabla\eta|^{pt}dx
     \right)^\frac{1}{t}\\
\nonumber\!\!&\!\!\leq\!\!&\!\! c_2^p  \|w_i-w\|^p_{W^{1,p}_Q(\Omega)}
 \Lnormb{\sqrt{Q}\nabla\eta}{pt}{B}^p,
\end{eqnarray}
and the last right-hand side tends to $0$ by \eqref{E3}. Hence,
$$\Lnormb{\sqrt{Q}\nabla\varphi_i-\big(\theta\eta^{\theta-1}w
  \sqrt{Q}\nabla\eta+ \eta^\theta\sqrt{Q}\nabla w\big)}{p}{\Omega}
  \longrightarrow0,$$
and in particular $\{\sqrt{Q}\nabla \varphi_i\}$ is a
Cauchy sequence in $[L^p(\Omega)]^n$. Claim 2 is thus proved.

It follows that the sequence $\{\varphi_i\}\subset Lip_0(B)$ identifies
an element of $(W^{1,p}_Q)_0(B)$, that $\varphi_i$ converges to
$\eta^\theta w$ a.e. in $\Omega$ and in $L^p(\Omega)$, and that
$$\sqrt{Q}\nabla(\eta^\theta w)=\theta\eta^{\theta-1}w\sqrt{Q}\nabla
w+\eta^\theta\sqrt{Q}\nabla w\qquad\text{a.e. in }\Omega.$$ The
proof is now complete.$\qquad\Box$

Next we will derive a result about higher local integrability of
functions in $W_Q^{1,p}(\Omega)$, whether or not they have compact
support.

\begin{pro}\label{prop1} Assume that \eqref{Cond1} and the Sobolev
inequality \eqref{Sobolev} hold. Let $0<\tau <1$ and $B=B(y,r)$ be a
$\rho$-ball with $r<r_1(y)$. Suppose that \eqref{E3} holds for $B$
and a function $\eta \in Lip_0(B)$ which equals $1$ on $B(y,\tau
r)$. With $t$ as in \eqref{E3}, assume that \eqref{3.6-0} holds for
$B$ with $t'$ given by $1/t+1/t'=1$.  Then $w\in
L^{p\sigma}(B(y,\tau r))$ for every $w\in W^{1,p}_Q(\Omega)$, and
$$\Lnormb{w}{p\sigma}{B(y,\tau r)}\leq C\|w\|_{W^{1,p}_Q(\Omega)},
$$
with $C>0$ depending on $p, \sigma, \max |\eta|, B$ and the constants
which arise in \eqref{Sobolev}, \eqref{E3} and \eqref{3.6-0}, but
independent of $w$.
\end{pro}

We note that under the hypotheses of Theorem \ref{main_thm}, by
using the functions $\eta_j$ in \eqref{cutoff}, the hypotheses of
Proposition \ref{prop1} are met for all $\rho$-balls $B=B(y,r)$ with
$r<r_1(y)$, and so the conclusion of Proposition \ref{prop1} holds
for all such $B$.

{\bf Proof:} Let $w \in W^{1,p}_Q(\Omega)$ and let $B, \eta$ satisfy
the hypotheses. Denote $\tau B= B(y,\tau r)$. Since $\eta =1$ on $\tau B$,
$$
\int_{\tau B} |w|^{p\sigma} dx \le \int_{B} |\eta
w|^{p\sigma} dx.
$$
By Proposition \ref{C4}, $\eta w \in (W^{1,p}_Q)_0(B)$ and satisfies
the product rule. Applying \eqref{Sobolev} with constant $c_1$,  we have
\begin{eqnarray}
\nonumber\frac{1}{|B|} \int_{B} |\eta w|^{p\sigma}dx &\leq&
(2c_1)^{p\sigma} \bigg[\frac{r^p}{|B|}\int_{B}|\sqrt{Q}
\nabla(\eta w)|^pdx+ \frac{1}{|B|}
\int_{B} |\eta w|^pdx \bigg]^\sigma \\
\nonumber&\leq& Cc_1^{p\sigma} \bigg[\frac{r^p}{|B|}
\int_{B}|\eta \sqrt{Q}\nabla w|^pdx +
\frac{r^p}{|B|} \int_{B}| w\sqrt{Q} \nabla\eta|^pdx\\
\label{3.6-2}&&\hspace{6.7cm}
+\frac{1}{|B|}\int_{B} |\eta w|^pdx \bigg]^\sigma \\
\nonumber&\leq& Cc_1^{p\sigma}\bigg[\frac{(r\,\max|\eta|)^p}{|B|}
  \int_{B}| \sqrt{Q}\nabla w|^pdx+ \frac{r^p}{|B|} \int_{B}|
w\sqrt{Q}\nabla\eta|^pdx\\
\nonumber&&\hspace{6cm}
+\frac{(\max|\eta|)^p}{|B|} \int_{B} |w|^pdx\bigg]^\sigma \\
\nonumber&\leq& C \left[\|w\|_{W^{1,p}_Q(B)}^{p} +\int_{B}|w
\sqrt{Q}\nabla \eta|^pdx\right]^\sigma,
\end{eqnarray}
where $C>0$ is a constant depending on $p,\sigma, B, \max |\eta|$ and
$c_1$. We will use \eqref{E3} and \eqref{3.6-0} to estimate the last
integral on the right; recall again that \eqref{3.6-0} holds for any
$w\in W^{1,p}_Q(\Omega)$. By H\"older's inequality (cf. \eqref{E4}),
\begin{eqnarray}
\nonumber\int_{B}|w\sqrt{Q}\nabla \eta|^pdx&\leq&
\bigg(\int_{B} |w|^{pt'}dx\bigg)^\frac{1}{t'}
\bigg(\int_{B}| \sqrt{Q}\nabla\eta|^{pt}dx \bigg)^\frac{1}{t}\\
\label{3.6-3}&\leq&
c_3^p\,c_2^p\,\|w\|_{W^{1,p}_Q(\Omega)}^p,\quad \text{where }c_3^{pt}
= \int_{B}| \sqrt{Q}\nabla\eta|^{pt}dx.
\end{eqnarray}
Combining estimates gives
$$\|\eta w\|_{p\sigma,\tau B;dx} \leq C \|w\|_{W^{1,p}_Q(\Omega)}
$$
with $C>0$ now also depending on $c_2$ and $c_3$. $\qquad\Box$

Condition \eqref{Cond0} provides a simple way to extend some of our
results proved for individual balls to general compact subsets of
$\Omega$. As an example, let us verify Remark \ref{useCond0} of the
Introduction. Let $\Omega'$ be a compact set in $\Omega$ and $u(x)$
be a function on $\Omega$ with the property that for all $B=B(y,r)$
with $r<r_1(y)$, $u$ is bounded on $\tau B$ for some $\tau \in
(0,1)$. For such $B$, by using \eqref{Cond0}, there is an open
concentric Euclidean ball $D\subset \tau B$. It follows from the
Heine-Borel Theorem that $\Omega'$ can be covered by a finite number
of such $D$, and so also by a finite number of balls $\tau B$ in
which $u$ is bounded. Consequently $u$ is bounded on $\Omega'$,
which verifies Remark \ref{useCond0}.

Similarly, \eqref{Cond0} leads to the following extension of
Proposition \ref{prop1}, whose proof we omit.

\begin{lem}\label{L1}  Assume that \eqref{Cond1} and
\eqref{Cond0} hold as well as the Sobolev inequality
\eqref{Sobolev}. Suppose that for each $y\in \Omega$, there is a
ball $B$ with center $y$ and radius $r<r_1(y)$ such that \eqref{E3}
holds for some $\eta \in Lip_0(B)$ which equals $1$ on $\tau B$ for
some $\tau \in (0,1)$ and some $t\ge 1$. Suppose also that
\eqref{3.6-0} holds for $B$ and $t'$ with $1/t+1/t'=1$. The values
of $\tau, t, t'$ may vary with $y$. Then for every compact subset
$\Omega'$ of $\Omega$, there is a constant $C$ depending on
$\Omega'$ so that
\begin{equation}\label{locglobsobolev}
\Lnormb{w}{p\sigma}{\Omega'}\leq C\|w\|_{\waQ}\quad\text{for all
$w\in W^{1,p}_Q(\Omega)$.}
\end{equation}
\end{lem}

In passing, we note that under the same hypotheses as in Lemma
\ref{L1}, the product rule in Proposition \ref{C4} extends to
Lipshitz functions $\eta$ supported in $\Omega$ (not just those
supported in a ball), provided $\eta$ satisfies the global condition
$$
\int_\Omega |\sqrt{Q}\nabla \eta|^{p\sigma'}\,dx < +\infty.
$$
The proof is similar to the one of Proposition \ref{C4}, using the
conclusion of Lemma \ref{L1} to modify the argument for \eqref{E4}. We
will not use this fact and so we omit the details of its proof.

\subsection{Weak Solutions}

As in the Introduction, we say that a pair $(u, \nabla u)
\in\waQ$ is a weak solution of equation \eqref{eqdiff} if
\begin{equation}\label{weaksol}
\int_\Omega\big[\nabla\varphi\cdot A(x,u,\nabla u)+\varphi
B(x,u,\nabla u)\big]=0\qquad\text{ for all }\varphi\in\lipo.
\end{equation}
If $(u, \nabla u)$ is a weak solution, we will sometimes simply say
that $u$ is a weak solution without explicitly mentioning $\nabla
u$. If $u$ is a weak solution, the class of functions
$\varphi$ for which \eqref{weaksol} holds can be enlarged from
$\lipo$; see Proposition \ref{Lambdaextd}. We shall refer to such
functions as {\it test functions}.

We start by showing that the notion of a weak solution is well-defined
and that the class of test functions can be enlarged from $\lipo$.

\begin{pro}\label{ABinLq} Assume that \eqref{struct} holds with
\begin{equation}\label{ranges2}
\gamma \in (1, \sigma(p-1)+1), \quad \psi \in (1,
p+1-\sigma^{-1}),\quad \delta \in (1, p\sigma),
\end{equation}
and that
\begin{equation}\label{integrability}
 \begin{array}{c}
  \begin{array}{ccc}
  c\in L^\frac{\sigma p}{\sigma p-1-\sigma(\psi-1)}_\text{loc}(\Omega),\,\,\,&\,\,\,e\in
   L^{p^\prime}_\text{loc}(\Omega),\,\,\,&\,\,\,f\in L^{(\sigma
    p)^\prime}_\text{loc}(\Omega),
  \end{array}\\
  \begin{array}{cc}
  b\in L^\frac{\sigma p}{\sigma\!(p-1)-\gamma+1}_\text{loc}(\Omega),\,\,\,&\,\,\,d\in
  L^\frac{p\sigma}{p\sigma-\delta}_\text{loc}(\Omega),
  \end{array}
 \end{array}
\end{equation}
where $p^\prime=\frac{p}{p-1}$, $(\sigma p)^\prime=\frac{\sigma
p}{\sigma p-1}$ are the conjugate exponents of $p$ and $\sigma p$
respectively. Let $0<\tau<1$ and $B= B(y,r)$ be a $\rho$-ball with
$r<\tau r_1(y)$, and the hypotheses of Proposition \ref{prop1} are
satisfied, but with $r$ there replaced by $r/\tau$.  Assume also
that \eqref{Sobolev} holds, and let ${\tilde A}(x,z,\xi)$ be defined
as in \eqref{struct}. Then for every $u\in\waQ$,
$$
\big|\widetilde{A}(\cdot,u,\nabla u)\big|\in L^{p^\prime}(B)
\quad\mbox{and}\quad B(\cdot,u,\nabla u)\in L^{(\sigma
p)^\prime}(B),
$$
with
\begin{equation*}
  \begin{array}{l}
    \Lnormb{\widetilde{A}(\cdot,u,\nabla u)}{p^\prime}{B}\leq\\
    \hspace{2.8cm} C_1\Big(p,a,\gamma,\Lnormb{\sqrt{Q}\nabla
    u}{p}{B},\Lnormb{e}{p}{B}, \Lnormb{b}{\frac{\sigma
    p}{\sigma(p-1)-\gamma+1}}{B},\Lnormb{u}{p\sigma}{B}\Big),\\
    \Lnormb{B(\cdot,u,\nabla u)}{(\sigma p)^\prime}{B}\leq\\
    \hspace{2.8cm}C_2\Big(p,\sigma,\delta,\psi,\Lnormb{c}{\frac{\sigma
    p}{\sigma p-1-\sigma(\psi-1)}}{B},
      \Lnormb{\sqrt{Q}\nabla u}{p}{B},\\
      \hspace{7.8cm}\Lnormb{d}{\frac{p\sigma}{p\sigma-
    \delta}}{B},\Lnormb{f}{(p\sigma)^\prime}{B},
      \Lnormb{u}{p\sigma}{B}\Big).
  \end{array}
\end{equation*}
\end{pro}

\textbf{Proof:} Let $B$ be a $\rho$-ball which satisfies the
hypotheses. By \eqref{struct},
\begin{eqnarray}
 \nonumber\int_{B}\big|\widetilde{A}(x,u(x),\nabla u(x))
   \big|^\frac{p}{p-1}\,dx&\leq&\int_{B} \left(a|\sqrt{Q}\nabla
   u|^{p-1}+b|u|^{\gamma-1}+e\right)^\frac{p}{p-1}\,dx\\
 \nonumber&\leq&\int_{B}3^\frac{p}{p-1}\left(a^\frac{p}{p-1}|\sqrt{Q}\nabla
    u|^p+b^\frac{p}{p-1}|u|^{p\frac{\gamma-1}{p-1}}+e^\frac{p}{p-1}\right)dx.
\end{eqnarray}
Using H\"{o}lder inequality with exponents $\sigma\frac{p-1}{\gamma
-1}$ and $\frac{\sigma (p-1)}{\sigma(p-1)-\gamma+1}$ on the second
term (note that $\sigma\frac{p-1}{\gamma-1} >1$ by \eqref{ranges2})
gives
\begin{eqnarray}
 \nonumber\int_{B}\big|\widetilde{A}(x,u(x),\nabla
    u(x))\big|^\frac{p}{p-1}\,dx&\leq&3^\frac{p}{p-1}
    \Big(a^\frac{p}{p-1}\Lnormb{\,|\sqrt{Q}\nabla u|\,}{p}{B}^p\\
 \nonumber&&\hspace{0,9cm}+\Lnormb{b}{\frac{\sigma p}{\sigma\!(p-1)
    -\gamma+1}}{B}^{p^\prime} \Lnormb{u}{p\sigma}{B}^{(\gamma-
    1)p^\prime}+\Lnormb{e}{p^\prime}{B}^{p^\prime}\Big).
\end{eqnarray}
This is finite by \eqref{grad-part}, \eqref{integrability} and
Proposition \ref{prop1}. In the same way,
\begin{eqnarray}
 \nonumber\int_{B}\big|B(x,u(x),\nabla u(x))\big|^\frac{\sigma
    p}{\sigma     p-1}\,dx&\leq&\int_{B}\left(c|\sqrt{Q}\nabla
    u|^{\psi-1}+d|u|^{\delta-1}+f\right)^\frac{\sigma p}{\sigma
    p-1}\,dx\\
 \nonumber&\leq&\int_{B}3^\frac{\sigma p}{\sigma p-1}
    \left(c^\frac{\sigma p}{\sigma p-1}|\sqrt{Q}\nabla
    u|^\frac{\sigma p(\psi-1)}{\sigma p-1}+d^\frac{\sigma p}{\sigma
    p-1}|u|^\frac{\sigma p(\delta-1)}{\sigma p-1} \right.\\
 \nonumber&&\hspace{6,1cm}\left.+f^\frac{\sigma p}{\sigma
 p-1}\right)\,dx.
\end{eqnarray}
Using H\"{o}lder inequality with exponents $\frac{\sigma
p-1}{\sigma(\psi-1)}$ and $\frac{\sigma p-1}{\sigma p-1-
\sigma(\psi-1)}$ on the first term and with $\frac{\sigma
p-1}{\delta-1}$ and $\frac{\sigma p -1}{\sigma p-\delta}$ on the
second one (note that $\frac{\sigma p-1}{\sigma(\psi-1)}, \frac{\sigma
p -1}{\sigma p-\delta} >1$ by \eqref{ranges2}) gives
\begin{eqnarray}
 \nonumber\int_{B}\big|B(x,u(x),\nabla u(x))\big|^\frac{\sigma
    p}{\sigma p-1}\,dx&\leq&   3^\frac{\sigma p}{\sigma p-1}
    \Big(\Lnormb{c}{\frac{\sigma p}{\sigma
    p-1-\sigma(\psi-1)}}{B}^{(p\sigma)^\prime}
    \Lnormb{\sqrt{Q}\nabla u}{p}{B}^\frac{p\sigma(\psi-
    1)}{\sigma p-1}\\
 \nonumber&&\hspace{1,1cm}+\Lnormb{d}{\frac{p\sigma}{p\sigma-
    \delta}}{B}^{(p\sigma)^\prime}
    \Lnormb{u}{p\sigma}{B}^\frac{\sigma p(\delta-1)}{\sigma p-1}
    +\Lnormb{f}{(p\sigma)^\prime}{B}^{(p\sigma)^\prime}\Big),
\end{eqnarray}
which is finite by \eqref{grad-part}, \eqref{integrability} and
Proposition \ref{prop1}. Indeed,
\bea
 \nonumber\int_{B}\!\!c^\frac{\sigma p}{\sigma p-1}|\sqrt{Q}\nabla
   u|^s\,dx\!\!&\!\!\leq\!\!&\!\!\left(\int_{B}\!\!c^\frac{\sigma p}{\sigma p-1-\sigma(\psi-1)}\,
   dx\!\right)^{\!\!\frac{\sigma p-1-\sigma(\psi-1)}{\sigma p-1}}\!\!\left(
   \int_{B}\!|\sqrt{Q}\nabla u|^p\,dx\!\right)^{\!\!\frac{\sigma(\psi-1)}{\sigma p-1}}\!<\!\infty,\\
 \nonumber\int_{B}d^{(\sigma p)^\prime}|u|^\frac{\sigma p(\delta-1)}{\sigma p-1}\,dx\!\!&\!\!\leq\!\!&\!\!
   \left(\int_{B}d^\frac{\sigma p}{\sigma p-\delta}\,dx\right)^\frac{\sigma p-\delta}{\sigma p-1}
   \left(\int_{B}|u|^{\sigma p}\,dx\right)^\frac{\delta-1}{\sigma
   p-1}<\infty.
\eea

Thus $\big|\widetilde{A}(\cdot,u,\nabla u)\big|\in
L^{p^\prime}(B)$ and $B(\cdot,u,\nabla u)\in L^{(\sigma
p)^\prime}(B)$, and the proposition is established.$\qquad\Box$

\begin{cor}\label{Lambdamap} Let the hypotheses of Proposition
\ref{ABinLq} hold and let $B = B(y,r)$ be a $\rho$-ball as described
there. Then for every $\varphi\in\text{Lip}_0(B)$ and every $u\in\waQ$,
$$
\int_\Omega\big[\nabla\varphi\cdot A(x,u,\nabla u)+\varphi
B(x,u,\nabla u) \big]\,dx < \infty.
$$
\end{cor}
\textbf{Proof:} Since no confusion should arise, we will use $B$ to
denote both the ball $B(y,r)$ and the function $B(x,u,\nabla u)$. Then
$$
\left|\int_\Omega\big[\nabla\varphi\cdot A(x,u,\nabla u)+\varphi
   B(x,u,\nabla u)\big]dx\right| \leq \int_{B} \big[|\nabla\varphi
   A(x,u,\nabla u)|+|\varphi B(x,u,\nabla u)|\big]dx
$$
$$
\leq \int_{B}|\sqrt{Q}\nabla\varphi|\,|\widetilde{A}(x,u,\nabla u)|\,dx
 +\int_{B}|\varphi|\,|B(x,u,\nabla u)|\,dx
$$
$$
\leq \Lnormb{\sqrt{Q}\nabla\varphi}{p}{B}
 \,\Lnormb{\widetilde{A}}{p^\prime}{B} +\Lnormb{\varphi}{\sigma p}{B}
 \,\Lnormb{B}{(\sigma p)^\prime}{B}.
$$

By the Sobolev inequality \eqref{Sobolev},
\begin{eqnarray}
 \nonumber\left|\int_\Omega\big[\nabla\varphi\cdot A(x,u,\nabla u)
   +\varphi B(x,u,\nabla u)\big]\,dx\right|\!&\!\leq\!&\!\|\varphi
   \|_{W^{1,p}_Q(B)}  \Lnormb{\widetilde{A}}{p^\prime}{B}\\
 \nonumber\!&\!\!&\!\hspace{0,6cm}+C(B,p)\|\varphi\|_{W^{1,p}_Q(B)}
   \Lnormb{B}{(\sigma p)^\prime}{B},
\end{eqnarray}
and hence
\begin{equation}\label{extendLambda}
  \begin{array}{l}
    \displaystyle\left|\int_\Omega\big[\nabla\varphi\cdot A(x,u,\nabla
    u)+\varphi B(x,u,\nabla u)\big] \,dx\right|\leq\\
    \hspace{4,4cm} \left(\Lnormb{\widetilde{A}}{p^\prime}{B}+
    C(B,p)\Lnormb{B}{(\sigma p)^\prime}{B}\right)
    \|\varphi\|_{W^{1,p}_Q(B)}.
  \end{array}
\end{equation}
The last quantity is finite because of Proposition \ref{ABinLq}, the
fact that $\varphi\in\text{Lip}_0(B)$, and our hypothesis that
$|Q|\in L^\frac{p}{2}_\text{loc}(\Omega)$.$\qquad\Box$

\begin{pro}\label{10} Under the hypotheses of Proposition \ref{ABinLq},
the map $\Lambda:\text{Lip}_0(B(y,r))\times\waQ\longrightarrow\rn$
defined by $$\Lambda(\varphi,u)=\int_\Omega\big[\nabla\varphi\cdot
A(x,u,\nabla u)+\varphi B(x,u,\nabla u)\big]\,dx$$ can be extended by
continuity so as to be defined on
$\big(W^{1,p}_{Q}\big)_0(B(y,r))\times\waQ$. Also, if
$\varphi\in\big(W^{1,p}_{Q}\big)_0(B(y,r))$ and $u\in\waQ$, then
\begin{equation}\label{Lambdaextd}
\Lambda(\varphi,u)=\int_\Omega\big[\sqrt{Q}\nabla\varphi\cdot
\widetilde{A}(x,u,\nabla u)+\varphi B(x,u,\nabla u)\big]\,dx.
\end{equation}
\end{pro}
\textbf{Proof:} We will again use $B$ to denote both $B(y,r)$
and $B(x,u,\nabla u)$. The map $\Lambda$ is well-defined
on $\text{Lip}_0(B)\times\waQ$ by Corollary \ref{Lambdamap}. For
fixed $u\in\waQ$, the map $\varphi\mapsto\Lambda(\varphi,u)$ is linear
in $\varphi\in \text{Lip}_0(B)$, and by \eqref{extendLambda},
$$|\Lambda(\varphi,u)|\leq C\|\varphi\|_{W^{1,p}_Q(B)},$$
with $C$ depending on $B$, $p$, $u$ and $\widetilde{A}$ but not on
$\varphi$. Then the linear map is continuous and can be extended by
continuity to $\big(W^{1,p}_{Q}\big)_0(B)$, since this is the
completion of $\text{Lip}_0(B)$.

In order to prove \eqref{Lambdaextd}, let $u\in\waQ$, $\varphi\in
\big(W^{1,p}_{Q}\big)_0(B)$, and $\{\varphi_i\}_{i\in\N}\subset
\text{Lip}_0(B)$ be a Cauchy sequence representing $\varphi$. Then \bea
 \label{4}\varphi_i&\longrightarrow&\varphi\quad\text{in
 }W^{1,p}_Q(B).
\eea
Moreover, by the previous estimates,
$$\left|\int_\Omega\big[\sqrt{Q}\nabla\varphi\cdot
\widetilde{A}(x,u,\nabla u)+\varphi B(x,u,\nabla
u)\big]\,dx\right|<\infty.$$ Then
\begin{equation*}
\begin{array}{l}
 \displaystyle\left|\int_\Omega\big[\sqrt{Q}\nabla\varphi\cdot
   \widetilde{A}(x,u,\nabla u)+\varphi B(x,u,\nabla
   u)\big]\,dx-\Lambda(\varphi_i,u)\right|\\\vspace{0,1cm}
 \hspace{2cm}\displaystyle=\left|\int_{B}\big[\sqrt{Q}\nabla\varphi\cdot
   \widetilde{A}+\varphi B-\sqrt{Q}\nabla\varphi_i\cdot\widetilde{A}-\varphi_i B\big]\,dx\right|\\\vspace{0,1cm}
 \hspace{2cm}\displaystyle\leq\int_{B}\big[|\sqrt{Q}\nabla\varphi
   -\sqrt{Q}\nabla\varphi_i|  |\widetilde{A}|+|\varphi-\varphi_i|
   |B|\big]\,dx\\\vspace{0,1cm}
 \hspace{2cm}\leq\displaystyle\Lnormb{\widetilde{A}}{p^\prime}{B}
   \Lnormb{\sqrt{Q}\nabla\varphi-\sqrt{Q}\nabla\varphi_i}{p}{B}
   +\Lnormb{B}{(\sigma p)^\prime}{B}\Lnormb{\varphi- \varphi_i}{\sigma
   p}{B}\\
 \hspace{2cm}\leq\left(\Lnormb{\widetilde{A}}{p^\prime}{B}+
   C\Lnormb{B}{(\sigma p)^\prime}{B}\right)\|\varphi-
   \varphi_i\|_{W^{1,p}_Q(B)},
\end{array}
\end{equation*}
where we used the Sobolev inequality \eqref{Sobolev} to obtain the
last inequality. Since $\|\varphi-\varphi_i\|_{W^{1,p}_Q(B)}\rightarrow0$
by \eqref{4}, we get
$$\Lambda(\varphi,u):=\lim_{i\rightarrow\infty}\Lambda(\varphi_i,
u)=\int_\Omega\big[\sqrt{Q}\nabla\varphi\cdot \widetilde{A}(x,u,
\nabla u)+\varphi B(x,u,\nabla u)\big]\,dx,$$
and \eqref{Lambdaextd} is established.$\qquad\Box$

\subsection{A useful test function}

Let $k,l,q,\mu, \beta\in\R$ with $q\geq1$, $l>k\geq0$, $\mu=p\sigma-1$
and $\beta =(\mu+1)q-\mu$. For any $t\in\R$, set
$\bar{t}=|t|+k$. Define
\bea
 \label{Ftest}F(\bar{t})&=&
    \begin{cases}
    \bar{t}^q\qquad\qquad\qquad\quad\,\qquad k\leq\bar{t}\leq l,\\
    q l^{q-1}\bar{t}-(q-1)l^q\quad\,\,\,\quad\bar{t}\geq l,
    \end{cases}\\\vspace{0.1cm}
 \label{Gtest}G(t)&=&\text{sign}(t)\left\{F(\bar{t})
F^\prime(\bar{t})^\mu-q^\mu k^\beta\right\}
    \qquad\text{ for }\,t\in\R.
\eea

As in the nondegenerate case studied in \cite{S}, we would like to use
the function
$$\phi(x)=\eta(x)^p G(u(x)),\qquad x\in\Omega,$$
as a test function in \eqref{weaksol}, where $\eta\in
\text{Lip}_0(\Omega)$ is any of the cutoff functions provided by
\eqref{cutoff}, and $u\in\waQ$ is a weak solution of the differential
equation \eqref{eqdiff}. In order to show that $\phi$ a feasible test
function, we begin by showing that there is a sequence
$\{l_i\}\subset\R^+$, $l_i\nearrow\infty$, such that if we choose
these $l$'s in definitions \eqref{Ftest} and \eqref{Gtest}, then
$G(u)\in\waQ$.

\begin{lem}\label{levelset}
Let $u\in L^\alpha(\Omega)$, $\alpha\in[1,\infty)$, $u\geq0$ a.e. in
$\Omega$. For any $l\in\R^+$, let $E_l=\{x\in\Omega:\,u(x)=l\}$.
Then the set
$$\Sigma =\left\{l\in\R^+:\,|E_l|>0\right\}$$
is countable.
\end{lem}
\textbf{Proof:} We claim that for every $\ep>0$, the set $\Sigma_\ep
=\{l>\ep:\,|E_l|>0\}$ is countable. For $j\in \N$, let
$\Sigma_{\ep,j}=\left\{l>\ep:\, |E_l|>\frac{1}{j}\right\}$. Then
$\Sigma_\ep=\bigcup_{j\in\N} \Sigma_{\ep,j},$ and it is enough to show
that each $\Sigma_{\ep,j}$ is countable. Fix $\epsilon, j$ and let
$\{l_i:i\in I\}$ be a sequence of distinct points in
$\Sigma_{\ep,j}$. Then
$$
\int_{\Omega}u^\alpha\,dx\geq\int_{\bigcup_{i\in I}
  E_{l_i}}u^\alpha\,dx = \sum_{i\in I} l_i^\alpha|E_{l_i}|\geq
\frac{\ep^\alpha}{j} \sum_{i\in I} 1.
$$
Since $u\in L^\alpha(\Omega)$, it follows that $\Sigma_{\ep,j}$ is actually
finite, and the claim follows. Since $\Sigma_\frac{1}{m}$ is countable
for every $m\in\N$, the set
$$\left\{l>0:\,|E_l|>0\right\}=\bigcup_{m\in\N}\Sigma_\frac{1}{m}$$
is countable too. Thus the set $\Sigma=\left\{l\geq0:\,|E_l|>0
\right\}$ is also countable, which proves the lemma.$\qquad\Box$

\begin{cor}\label{goodlevelsets} Given a sequence $\{u_i\}_{i\in\N}\subset
L^\alpha(\Omega)$, there is a sequence of positive numbers
$l_j\nearrow\infty$ such that
\begin{eqnarray}
\nonumber|E_{i,l_j}|&=&\left|\left\{x\in\Omega:\,|u_i(x)|=l_j
\right\}\right|=0\qquad\text{ for all }\,i,j \in\N.
\end{eqnarray}
\end{cor}

\textbf{Proof:} The sets $\Sigma_i:=\{l\in\R^+:\,|E_{i,l}|>0\}$ are
countable for every $i$ by Lemma \ref{levelset}, and hence the set
$\Sigma: =\bigcup_i \Sigma_i=\{l\in\R^+:\,|E_{i,l}|>0
\text{ for some }i\}$ is countable. Then $\R^+\setminus\Sigma$ is
uncountable, and in particular there is a sequence
$\{l_j\}_{j\in\N}\subset\R^+\setminus\Sigma$ such that
$l_j\nearrow\infty$. Since $l_j\in\R^+\setminus\Sigma$ for every $j$,
we have
$$|E_{i,l_j}|=0\qquad\text { for all }i, j \in\N.$$
The corollary is proved.$\qquad\Box$

The next fact can be proved in a similar way.
\begin{cor}\label{goodlevelsets2} Given a sequence
$\{u_i\}_{i\in\N}\subset L^\alpha(\Omega)$, there is a sequence of
positive numbers $\lambda_j\searrow0^+$ such that
\begin{eqnarray}
\nonumber|E_{i,\lambda_j}|&=&\left|\left\{x\in\Omega:\,|u_i(x)|=
\lambda_j\right\}\right|=0\qquad\text{ for all }i, j \in\N.
\end{eqnarray}
\end{cor}

Lemma \ref{levelset} and Corollaries \ref{goodlevelsets} and
\ref{goodlevelsets2} provide a means to avoid using the notion
of regular gradient introduced in \cite{SW2}.  This simplifies some
technical aspects in \cite{SW2} and leads to relatively short proofs
of results like Theorem \ref{Gtest-thm} and Lemma \ref{appendix1}.

\begin{thm}\label{Gtest-thm} Let $k\geq0$, $q\geq1$, $\mu=p\sigma-1$,
$\beta=(\mu+1)q-\mu$ and $\bar{t}=|t|+k$ for every $t\in\R$. Given a
function $u\in\waQ$, there exists a sequence $l_j\nearrow\infty$ such
that if we define
\bea
 \nonumber F_j(\bar{t})&=&
    \begin{cases}
    \bar{t}^q\qquad\qquad\qquad\quad\qquad k\leq\bar{t}\leq l_j,\\
    q l_j^{q-1}\bar{t}-(q-1)l_j^q\quad\,\,\,\quad\bar{t}\geq l_j,
    \end{cases}\\\vspace{0.1cm}
 \nonumber G_j(t)&=&\text{\emph{sign}}(t)\left\{F_j(\bar{t})
F_j^\prime(\bar{t})^\mu-q^\mu k^\beta\right\}
    \qquad\text{ for }t\in\R,
\eea then $G_j(u)\in\waQ$ for all $j$, and
\begin{equation}\label{G-Qgrad}
 \sqrt{Q}\nabla\left(G_j(u)\right)=G^\prime_j(u)\sqrt{Q}\nabla
    u\qquad\text{a.e. in }\Omega.
\end{equation}
\end{thm}

\begin{rem}\label{rem1}
In the proof of Theorem \ref{Gtest-thm}, we will use the following
facts for every $j\in\N$:
\begin{itemize}
\item[i)] $F_j\in\mathcal{C}^1\left([k,\infty)\right)$, with
\bea
 \nonumber F^\prime_j(\bar{t})&=&
    \begin{cases}
    q\bar{t}^{q-1}\qquad\qquad\quad k\leq\bar{t}\leq l_j,\\
    q l_j^{q-1}\,\qquad\qquad\quad\bar{t}\geq l_j,
    \end{cases}\\\vspace{0.1cm}
 \nonumber F^{\prime\prime}_j(\bar{t})&=&
    \begin{cases}
    q(q-1)\bar{t}^{q-2}\qquad k<\bar{t}<l_j,\\
    0\!\qquad\qquad\qquad\quad\bar{t}>l_j.
    \end{cases}
\eea
\item[ii)] $0\leq F_j(\bar{t})\leq ql_j^{q-1}\bar{t}+l^q_j$,\\
$0\leq F^\prime_j(\bar{t})\leq q
l_j^{q-1}$,\\
$F^{\prime\prime}_j(\bar{t})\geq0$ and $F^{\prime\prime}_j(\bar{t})$
is a bounded function away from $\bar{t}=0$.
\item[iii)] $G_j\in\mathcal{C}^0(\R)$ and $G_j$ is differentiable
everywhere except at $\pm(l_j-k)$ where it has ``corners''. Indeed \bea
 \nonumber
 G^\prime_j(t)&=& F^\prime_j(\bar{t})^{\mu+1}+
    \mu F_j(\bar{t})F^{\prime\prime}_j(\bar{t})
F^\prime_j(\bar{t})^{\mu-1}\\\vspace{0.1cm}
 \nonumber &=&\begin{cases}
    q^{\mu+1}\bar{t}^{(\mu+1)(q-1)}+\mu(q-1)q^\mu\bar{t}^{q+q-2+(q-1)(\mu-1)}\quad\quad 0<|t|<l_j-k,\\
    q^{\mu+1}l_j^{(\mu+1)(q-1)}\quad\qquad\qquad\qquad\qquad\qquad\qquad\quad\,\,\,|t|>l_j-k,
    \end{cases}\\\vspace{0.1cm}
\nonumber &=&
    \begin{cases}
    \beta q^\mu\bar{t}^{(\mu+1)(q-1)}\qquad\qquad 0<|t|<l_j-k,\\
    q^{\mu+1}l_j^{(\mu+1)(q-1)}\qquad\qquad\,\,\,\,\,\,\,\,|t|>l_j-k,
    \end{cases}\\\vspace{0.1cm}
 \nonumber &=&
    \begin{cases}
    q^{-1}\beta F^\prime_j(\bar{t})^{\mu+1}\qquad\quad\,\, 0<|t|<l_j-k,\\
    F^\prime_j(\bar{t})^{\mu+1}\qquad\qquad\quad\,\,\,|t|>l_j-k.
    \end{cases}
\eea \noindent The discontinuity of $G^\prime_j$ at $0$ is removable.
\item[iv)] $|G_j(t)|\leq
F_j(\bar{t}) F^\prime_j(\bar{t})^\mu\leq\left(q
  l_j^{q-1}\right)^\mu F_j(\bar{t})$,\\
$0\leq G_j^\prime(t)\leq
q^{-1}\beta F_j^\prime(\bar{t})^{\mu+1}\leq q^{-1}\beta
l_j^{(q-1)(\mu+1)}q^{\mu+1}=\beta q^\mu l_j^{(\mu+1)(q-1)}<\infty.$
\end{itemize}
\end{rem}
\textbf{Proof of Theorem \ref{Gtest-thm}}: Let $u\in\waQ$ and
$\{u_i\}\subset\lipl$ be a sequence representing $u$. Then
up to subsequences,
\bea
 \nonumber u_i\!&\!\longrightarrow\!&\!u\qquad\qquad\quad\text{ in
 }\waQ,\\
 \nonumber u_i\!&\!\longrightarrow\!&\!u\qquad\qquad\quad\text{ in }\lp \text{
   and a.e. in }\Omega\text{ and}\\
 \nonumber \sqrt{Q}\nabla u_i\!&\!\longrightarrow\!&\!\sqrt{Q}\nabla u\,\,\qquad\text{ in
  }[\lp]^n\text{ and a.e. in }\Omega.
\eea Use Corollary \ref{goodlevelsets} to choose a sequence $r_j
\nearrow\infty$ such that
\bea
 \nonumber|E_{i,r_j}|\!&\!=\!&\!\left|\{x\in\Omega:\,|u_i(x)|=r_j\}\right|=0
    \qquad\text{ for }i,j \in\N,\\
 \nonumber|E_{r_j}|\!&\!=\!&\!\left|\{x\in\Omega:\,|u(x)|=r_j\}\right|=0
    \qquad\,\,\text{ for }j\in\N.
\eea Define $l_j=r_j+k$ for $j\in\N$.

\textit{Claim 1:} $G_j(u_i)\in\lipl$ for every $i,j \in\N$.\\
By the fundamental theorem of calculus,
$$G_j(t)=\int_0^tG^\prime_j(s)\,ds+G(0)\qquad\text{ for all }t\in\R,j\in\N.$$
Then for $x,y\in\Omega$,
$$\left|G_j(u_i(x))-G_j(u_i(y))\right|=\left|\int_{u_i(y)}^{u_i(x)}G_j^\prime(s)\,ds\right|
\leq\|G^\prime_j\|_\infty|u_i(x)-u_i(y)|.$$ Since
$\|G^\prime_j\|_\infty$ is finite by Remark \ref{rem1}, and since
$u_i$ is locally Lipschitz continuous, it follows that
$G_j(u_i)\in\lipl$ for every $i,j \in\N$. In particular,
$\nabla(G_j(u_i))$ is well-defined a.e. in $\Omega$.

\textit{Claim 2:} For almost every $x\in\Omega$,
\begin{equation}\label{7}
  \nabla(G_j(u_i))(x)=G^\prime_j(u_i(x))\nabla
  u_i(x)\qquad\text{ for all $i,j \in\N$, and}
\end{equation}
$$\int_\Omega\big|\sqrt{Q}\nabla(G_j(u_i))\big|\,dx<\infty.$$
Indeed, consider $x\in\Omega$ such that
\begin{itemize}
\item[i)] $\nabla u_i(x)$, $\nabla(G_j(u_i))(x)$ are defined,
\item[ii)]$|u_i(x)|\neq l_j-k=r_j$ for every $i,j$.
\end{itemize}
Then $G_j(t)$ is $\mathcal{C}^1$ in a neighborhood of $u_i(x)$,
since $G_j(t)$ has corners only at $t=\pm(l_j-k)$ by Remark
\ref{rem1}. By the chain rule, formula \eqref{7} holds at the point
$x$.

Since the set of points for which either (i) or (ii) does not hold has
Lebesgue measure $0$, formula \eqref{7} holds for
a.e. $x\in\Omega$. Thus for every $i$ and $j$,
\bea
 \nonumber\int_\Omega\big|\sqrt{Q}\nabla(G_j(u_i))\big|\,dx&=&\int_\Omega\big|G^\prime_j(u_i(x))\big|\big|\sqrt{Q}\nabla
  u_i(x)\big|\,dx\\\vspace{0,1cm}
 \nonumber&\leq&\|G_j^\prime\|_\infty\int_\Omega\big|\sqrt{Q}\nabla
  u_i(x)\big|\,dx<\infty.
\eea Claim 2 follows.

\textit{Claim 3:} $G_j(u_i)\longrightarrow G_j(u)$ a.e. in $\Omega$
and in $\lp$ for all $j$.\\
Since $u_i\longrightarrow u$ a.e. in $\Omega$ and $G_j$ is
continuous for every $j$, we obviously have that $G_j(u_i)
\longrightarrow G_j(u)$ a.e. in $\Omega$ for every $j$. Then
$$|G_j(u_i)-G_j(u)|^p\longrightarrow0\qquad\text{a.e. in }\Omega,\,\,\text{ for all }j\in\N.$$
We also have for every $x\in\Omega$ that
\bea
\nonumber\left|G_j(u_i)-G_j(u)\right|^p&\leq&2^p\left(|G_j(u_i)|^p+|G_j(u)|^p\right)\\\vspace{0,1cm}
 \nonumber&\leq&2^p\left(|F_j(\bar{u}_i)|^p\,|F^\prime_j(\bar{u}_i)|^{p\mu}+|F_j(\bar{u})|^p\,
    |F^\prime_j(\bar{u})|^{p\mu}\right)\\\vspace{0,1cm}
 \nonumber&\leq&2^p\left(q l_j^{q-1}\right)^{p\mu}\left(|F_j(\bar{u}_i)|^p+
    |F_j(\bar{u})|^p\right)\\\vspace{0,1cm}
 \nonumber&\leq&2^p\left(q l_j^{q-1}\right)^{p\mu}\left(\left(q l_j^{q-1}\bar{u}_i+
    l_j^q\right)^p+\left(q l_j^{q-1}\bar{u}+l_j^q\right)^p\right)\\\vspace{0,1cm}
 \label{8}&\leq&2^{2p}\left(q l_j^{q-1}\right)^{p\mu}
    \left(\left(q l_j^{q-1}\right)^p\bar{u}_i^p+\left(q
    l_j^{q-1}\right)^p\bar{u}^p+2l_j^{q p}\right)\\\vspace{0,1cm}
 \nonumber&\leq&2^{2p}\left(q l_j^{q-1}\right)^{p\mu}
    \max{\!\left(\left(q l_j^{q-1}\right)^p,2l_j^{q
    p}\right)}\left(\bar{u}_i^p+\bar{u}^p+1\right)\\\vspace{0,1cm}
 \nonumber&=&2^{2p}\left(q l_j^{q-1}\right)^{p\mu}
    \max{\!\left(\left(q l_j^{q-1}\right)^p,2l_j^{q
    p}\right)}\\\vspace{0,1cm}
\nonumber&&\hspace{4,5cm}\left((|u_i|+k)^p+(|u|+k)^p+1\right)\\\vspace{0,1cm}
 \nonumber\phantom{\left|G_j(u_i)-G_j(u)\right|^p}&\leq&2^{3p}\left(q l_j^{q-1}\right)^{p\mu}
    \max{\!\left(\left(q l_j^{q-1}\right)^p,2l_j^{q p}\right)}\\\vspace{0,1cm}
 \nonumber&&\hspace{4,5cm}\left(|u_i|^p+|u|^p+(2k^p+1)\right)\\\vspace{0,1cm}
 \nonumber&\leq&2^{3p}\left(q l_j^{q-1}\right)^{p\mu}
    \max{\!\left(\left(q l_j^{q-1}\right)^p,2l_j^{q p}\right)}\\\vspace{0,1cm}
 \nonumber&&\hspace{4,5cm}(2k^p+1)\left(|u_i|^p+|u|^p+1\right)\\\vspace{0,1cm}
 \nonumber&=&C(p,\sigma,q,l_j,k)\left(|u_i|^p+|u|^p+1\right),
\eea and $\left(|u_i|^p+|u|^p+1\right)\in L^1(\Omega)$ for every
$i\in\N$.

Since $\left(|u_i|^p+|u|^p+1\right)\longrightarrow\left(2|u|^p+1
\right)$ for a.e. $x\in\Omega$, and since
\bea
 \nonumber\int_\Omega\left(|u_i|^p+|u|^p+1\right)\,dx&=&\Lnormb{u_i}{p}{\Omega}^p+\Lnormb{u}{p}{\Omega}^p+|\Omega|\\\vspace{0,1cm}
 \nonumber&\rightarrow&2\Lnormb{u}{p}{\Omega}^p+|\Omega|\,\,\,\,=\,\,\,\,\int_\Omega\left(2|u|^p+1\right)\,dx,
\eea the Lebesgue Sequentially Dominated Convergence Theorem gives
$$G_j(u_i)\longrightarrow G_j(u)\qquad\text{in }\lp \,\text{
for all }j.$$ Claim 3 is thus proved.

\textit{Claim 4:} $\sqrt{Q}\nabla(G_j(u_i))\longrightarrow
G^\prime_j(u)\sqrt{Q}\nabla u$ a.e. in $\Omega$
and in $[\lp]^n$, for all $j\in\N$.\\
Consider a point $x\in\Omega$ such that
\begin{itemize}
\item[i)] $\nabla(G_j(u_i))(x)$ and $\sqrt{Q}\nabla u(x)$ are
defined,
\item[ii)] $\nabla(G_j(u_i))(x)=G^\prime_j(u_i(x))\nabla u_i(x)$,
\item[iii)] $|u(x)|\neq r_j=l_j-k\,\,\,\,\,$ and
$\,\,\,\,\,|u_i(x)|\neq r_j$ for every $i,j$,
\item[iv)] $u_i(x)\longrightarrow u(x)$,
\item[v)] $\sqrt{Q}\nabla u_i(x)\longrightarrow\sqrt{Q}\nabla u(x)$.
\end{itemize}
The set of points which do not satisfy one or more of these
conditions has Lebesgue measure $0$ for the following reasons,
respectively:
\begin{itemize}
\item[i)] because $G_j(u_i)$ is locally Lipschitz for every $i,j$ by
claim 1, and $u\in\waQ$,
\item[ii)] by claim 2,
\item[iii)] by our choice of the sequence $r_j$,
\item[iv)] because $u_i\longrightarrow u$ a.e. in $\Omega$,
\item[v)] because $\sqrt{Q}\nabla u_i\longrightarrow\sqrt{Q}\nabla u$ a.e. in $\Omega$.
\end{itemize}
For any $x\in\Omega$ satisfying all these conditions,
$$\sqrt{Q}\nabla(G_j(u_i))(x)=G_j^\prime(u_i(x))\sqrt{Q}\nabla
u_i(x)\longrightarrow G_j^\prime(u(x))\sqrt{Q}\nabla u(x)$$
since $G_j^\prime(t)$ is continuous everywhere except at $t=\pm(l_j-k)$
while $|u(x)|\neq l_j-k$. Thus,
$$\sqrt{Q}\nabla(G_j(u_i))\longrightarrow G_j^\prime(u)\sqrt{Q}\nabla
u\qquad\text{ a.e in }\Omega.$$

On the other hand, a.e. in $\Omega$ and for every $i,j\in\N$,
\bea
 \nonumber\left|\sqrt{Q}\nabla(G_j(u_i))-G^\prime_j(u)\sqrt{Q}\nabla
    u\right|^p&=&\left|G^\prime_j(u_i)\sqrt{Q}\nabla
    u_i-G^\prime_j(u)\sqrt{Q}\nabla u\right|^p\\\vspace{0,1cm}
 \nonumber&\leq&\left(|G^\prime_j(u_i)|\,|\sqrt{Q}\nabla u_i|+|G^\prime_j(u)|\,|\sqrt{Q}\nabla u|\right)^p\\\vspace{0,1cm}
 \nonumber&\leq&2^p\left(\|G^\prime_j\|^p_\infty\,|\sqrt{Q}\nabla u_i|^p+\|G^\prime_j\|_\infty^p\,|\sqrt{Q}\nabla u|^p\right)\\\vspace{0,1cm}
 \nonumber&\leq&2^p\left(\beta q^\mu
l_j^{(\mu+1)(q-1)}\right)^p\left(|\sqrt{Q}\nabla
u_i|^p+|\sqrt{Q}\nabla u|^p\right),
\eea and the functions on the right in the last inequality belong to
$L^1(\Omega)$ for every $i,j$.

Also, $(|\sqrt{Q}\nabla u_i|^p+|\sqrt{Q}\nabla u|^p)\longrightarrow
2|\sqrt{Q}\nabla u|^p\,$ for a.e. $x\in\Omega$, and
\bea
 \nonumber\int_\Omega\left(|\sqrt{Q}\nabla u_i|^p+|\sqrt{Q}\nabla
    u|^p\right)dx&=&\Lnormb{\sqrt{Q}\nabla u_i}{p}{\Omega}^p
+\Lnormb{\sqrt{Q}\nabla u}{p}{\Omega}^p\\\vspace{0,1cm}
 \nonumber&\rightarrow&2\Lnormb{\sqrt{Q}\nabla u}{p}{\Omega}^p
= \int_\Omega2|\sqrt{Q}\nabla u|^p\,dx.
\eea Then by the Lebesgue Sequentially Dominated Convergence Theorem
we have $$\sqrt{Q}\nabla(G_j(u_i))\longrightarrow
G^\prime_j(u)\sqrt{Q}\nabla u\qquad\text{ in }[\lp]^n\,\text{ for
all }j\in\N.$$ Claim 4 is thus proved.

Claims 3 and 4 together prove that $G_j(u)\in\waQ$, that
$\{G_j(u_i)\}_{i\in\N}$ is a sequence of locally Lipschitz functions
in $\Omega$ which represents $G_j(u)$ and that
$$\sqrt{Q}\nabla(G_j(u))=G_j^\prime(u)\sqrt{Q}\nabla u\qquad\text{ a.e. in }\Omega,$$
which is formula \eqref{G-Qgrad}.$\qquad\Box$

\section{\textbf{Proof of Theorem \ref{main_thm}.}}\label{proof}

\emph{Step 1.} We will use the notation
$$\avnorm{w}{\alpha}:=\avnormb{w}{\alpha}{B_r}=\left(\frac{1}{|B_r|}
\int_{B_r}|w|^\alpha\,dx\right)^\frac{1}{\alpha}
=\left(\fint_{B_r}|w|^\alpha\,dx\right)^\frac{1}{\alpha}$$ for any
$\alpha\geq1$, any function $w$ and any $\rho$-ball $B_r=B(y,r)$
with $0<r<r_1(y)$. For $k>0$, define
\begin{eqnarray}
\nonumber\bar{z}&=&|z|+k,\qquad z\in\R,\\
\nonumber\bar{b}(x)&=&b(x)+k^{1-p}e(x),\qquad x\in\Omega,\\
\nonumber\bar{h}(x)&=&h(x)+k^{-p}g(x),\qquad x\in\Omega,\\
\nonumber\bar{d}(x)&=&d(x)+k^{1-p}f(x),\qquad x\in\Omega.
\end{eqnarray}
Then the following new structural inequalities for the coefficients
are easily obtained from \eqref{struct} with $\gamma=\delta=p$:
\begin{eqnarray}
\nonumber\xi\cdot A(x,z,\xi)&\geq&
a^{-1}|\sqrt{Q(x)}\cdot\xi|^p-\bar{h}(x)\bar{z}^p,\\
\label{struct2}\left|\wt{A}(x,z,\xi)\right|&\leq&a|\sqrt{Q(x)}
\cdot\xi|^{p-1}+\bar{b}(x)\bar{z}^{p-1},\\
\nonumber\left|B(x,z,\xi)\right|&\leq&c|\sqrt{Q(x)}
\cdot\xi|^{\psi-1}+\bar{d}(x)\bar{z}^{p-1}.
\end{eqnarray}

In fact, when we deal with a specific solution pair $(u, \nabla
u)$, as is the case in Theorem \ref{main_thm}, we will only need to
assume the analogue of \eqref{struct2} in which $z$ and $\xi$ are
replaced respectively by $u(x)$ and $\nabla u(x)$ for all $x\in
\Omega$.

Now consider the functions $F_j(\bar{t})$ and $G_j(t)$
defined in Theorem \ref{Gtest-thm}. Let $\eta\in\text{Lip}_0(B_r)$
be any of the Lipschitz cutoff functions provided by \eqref{cutoff}
for a $\rho$--ball $B_r$. Then by Theorem \ref{Gtest-thm}, Corollary
\ref{C4} and Proposition \ref{10}, each function $\varphi_j(x):
=\eta(x)^p G_j(u(x))$ is a feasible test function in \eqref{weaksol}.

In order to simplify notation, we will not explicitly show the
dependence of $A,\,\wt{A}$ or $B$ on the variables $x,\,u(x)$ and
$\nabla u(x)$. Also, we will often not show the dependence of any
function of $x$ on $x$.

\emph{Step 2.} We start by deriving  some pointwise estimates which give
lower bounds for $\nabla\varphi_j\cdot A+\varphi_jB$. By the
structural conditions \eqref{struct2},
\begin{eqnarray}
\nonumber\nabla\varphi_j\cdot
    A+\varphi_jB\!&\!=\!&\!\sqrt{Q}\nabla\varphi_j\cdot\wt{A}+\varphi_jB\\
\nonumber\!&\!=\!&p\eta^{p-1}G_j(u)\sqrt{Q}\nabla\eta\cdot\wt{A}+\eta^pG^\prime_j(u)\sqrt{Q}\nabla
    u\cdot\wt{A}+\eta^pG_j(u)B\\
\nonumber\!&\!\geq\!&\!\eta^pG^\prime_j(u)\left[a^{-1}\big|\sqrt{Q}\nabla
    u\big|^p-\bar{h}\bar{u}^p\right]-\eta^p|G_j(u)|
    \left[c\big|\sqrt{Q}\nabla u\big|^{\psi-1}+\bar{d}\bar{u}^{p-1}\right]\\
\nonumber\!&\!\!&\!\hspace{3,2cm}-p\eta^{p-1}|G_j(u)|\big|\sqrt{Q}\nabla\eta\big|\left[a\big|\sqrt{Q}\nabla
    u\big|^{p-1}+\bar{b}\bar{u}^{p-1}\right].
\end{eqnarray}

Then it follows from Remark \ref{rem1} that
\begin{eqnarray}
\nonumber\nabla\varphi_j\cdot
    A+\varphi_jB\!&\!\geq\!&\!a^{-1}\eta^p F_j^\prime(\bar{u})^{\mu+1}
    \big|\sqrt{Q}\nabla u\big|^p\\
\nonumber\!&\!\!&\!\hspace{1cm}-\eta^pF_j(\bar{u}) F_j^\prime(\bar{u})^\mu
    \left[c\big|\sqrt{Q}\nabla u\big|^{\psi-1}+\bar{d}\bar{u}^{p-1}\right]\\
\nonumber\!&\!\!&\!\hspace{1cm}-\beta
    q^{-1} F_j^\prime(\bar{u})^{\mu+1}\eta^p\bar{h}\bar{u}^p- p\eta^{p-1}\bar{b}F_j(\bar{u})F_j^\prime(\bar{u})^\mu
    \big|\sqrt{Q}\nabla\eta\big|\bar{u}^{p-1}\\
\nonumber\!&\!\!&\!\hspace{1cm}-ap\eta^{p-1}F_j(\bar{u})
    F_j^\prime(\bar{u})^\mu\big|\sqrt{Q}\nabla\eta\big|\big|\sqrt{Q}\nabla
    u\big|^{p-1}
\end{eqnarray}
\begin{eqnarray}
\nonumber\hspace{3cm}\!&\!=\!&\!a^{-1}\eta^p
F_j^\prime(\bar{u})^{\mu+1-p}
    \big|F^\prime_j(\bar{u})\sqrt{Q}\nabla  u\big|^p\\
\nonumber\!&\!\!&\!\hspace{1cm}-c\eta^{p+1-\psi}F_j(\bar{u})
    F_j^\prime(\bar{u})^{\mu+1-\psi}\big|\eta
    F_j^\prime(\bar{u})\sqrt{Q}\nabla u\big|^{\psi-1}\\
\nonumber\!&\!\!&\!\hspace{1cm}-p\eta^{p-1}\bar{b}
    F_j^\prime(\bar{u})^{\mu+1-p}
    \big|F_j(\bar{u})\sqrt{Q}\nabla\eta\big|
    \big(F_j^\prime(\bar{u})\bar{u}\big)^{p-1}\\
\nonumber\!&\!\!&\!\hspace{1cm}-ap F_j^\prime(\bar{u})^{\mu+1-p}
    \big|F_j(\bar{u})\sqrt{Q}\nabla\eta\big|\big|\eta
    F^\prime_j(\bar{u})\sqrt{Q}\nabla u\big|^{p-1}\\
\nonumber\!&\!\!&\!\hspace{1cm}-\bar{h}\beta
    q^{-1}\eta^p \big(F^\prime_j(\bar{u})\bar{u}\big)^p
    F_j^\prime(\bar{u})^{\mu+1-p}\\
\nonumber\!&\!\!&\!\hspace{1cm}-\bar{d}F_j(\bar{u})\eta^p
    \big(F_j^\prime(\bar{u})\bar{u}\big)^{p-1}
    F_j^\prime(\bar{u})^{\mu+1-p}\\
\nonumber\!&\!\geq\!&\!a^{-1}\eta^p F_j^\prime(\bar{u})^{\mu+1-p}
    \big|F^\prime_j(\bar{u})\sqrt{Q}\nabla  u\big|^p\\
\nonumber\!&\!\!&\!\hspace{1cm}-c\eta^{p+1-\psi}F_j(\bar{u})
    F_j^\prime(\bar{u})^{\mu+1-\psi}\big|\eta
    F_j^\prime(\bar{u})\sqrt{Q}\nabla u\big|^{\psi-1}\\
\nonumber\!&\!\!&\!\hspace{1cm}-pq^{p-1}\eta^{p-1}\bar{b}
    F_j^\prime(\bar{u})^{\mu+1-p}
    \big|F_j(\bar{u})\sqrt{Q}\nabla\eta\big| F_j(\bar{u})^{p-1}\\
\nonumber\!&\!\!&\!\hspace{1cm}-ap F_j^\prime(\bar{u})^{\mu+1-p}
    \big|F_j(\bar{u})\sqrt{Q}\nabla\eta\big|\big|\eta
    F^\prime_j(\bar{u})\sqrt{Q}\nabla u\big|^{p-1}\\
\nonumber\!&\!\!&\!\hspace{1cm}-\bar{h}\beta
    q^{p-1}\eta^p F_j(\bar{u})^p F_j^\prime(\bar{u})^{\mu+1-p}\\
\nonumber\!&\!\!&\!\hspace{1cm}-\bar{d}q^{p-1}\eta^p F_j(\bar{u})^p
    F_j^\prime(\bar{u})^{\mu+1-p}.
\end{eqnarray}
Although the last two terms are identical apart from the
multiplicative factor $\beta$, we will treat them separately in
order to simplify calculations later in the corollaries following
our main theorem. By Lemma \ref{appendix1} in Appendix 1, we can
replace $\big|\sqrt{Q}\nabla u\big|$ by
$\big|\sqrt{Q}\nabla\bar{u}\big|$ in the previous inequalities.
Setting $v_j: =F_j(\bar{u}),$ we have $\sqrt{Q}\nabla
v_j=F^\prime_j(\bar{u})\sqrt{Q}\nabla\bar{u}$ by Lemma \ref{appendix2}
in Appendix 1. Thus
\begin{eqnarray}
\nonumber\nabla\varphi_j\cdot
    A+\varphi_jB\!&\!\geq\!&\!a^{-1} F_j^\prime(\bar{u})^{\mu+1-p}
    \big|\eta\sqrt{Q}\nabla v_j\big|^p\\
\nonumber\!&\!\!&\!\hspace{0,3cm}-c\eta^{p+1-\psi}v_j
    F_j^\prime(\bar{u})^{\mu+1-\psi}\big|\eta\sqrt{Q}\nabla
    v_j\big|^{\psi-1}\\
\label{1}\!&\!\!&\!\hspace{0,3cm}-pq^{p-1}\eta^{p-1}\bar{b}
    F_j^\prime(\bar{u})^{\mu+1-p}\big|v_j\sqrt{Q}\nabla\eta\big|v_j^{p-1} \\
\nonumber\!&\!\!&\!\hspace{0,3cm}-ap F_j^\prime(\bar{u})^{\mu+1-p}
\big|v_j\sqrt{Q}\nabla\eta\big|\big|\eta\sqrt{Q}\nabla
v_j\big|^{p-1}\\
\nonumber\!&\!\!&\!\hspace{0,3cm}-\bar{h}\beta
    q^{p-1}\eta^pv_j^p F_j^\prime(\bar{u})^{\mu+1-p}-\bar{d}q^{p-1}
\eta^pv_j^p  F_j^\prime(\bar{u})^{\mu+1-p}.
\end{eqnarray}
Now, recalling that $u$ is a weak solution of \eqref{eqdiff} and that
$\eta$ and thus $\varphi_j$ have support in $B_r$, we have
\begin{eqnarray}\label{99999}\int_{B_r}\big[\nabla\varphi_j\cdot A+\varphi_jB\big]\,dx=\int_\Omega\big[\nabla\varphi_j\cdot A+\varphi_jB\big]\,dx=0.
\end{eqnarray}
Integrating \eqref{1} over $\Omega$, dividing by $|B_r|$ and
rearranging terms we then get
\begin{eqnarray}
\nonumber
    a^{-1}\fint_{B_r}\!F_j^\prime(\bar{u})^{\mu+1-p}\big|
    \eta\sqrt{Q}\nabla  v_j\big|^p\,dx\!&\!\leq\!&\!ap\fint_{B_r}\!
    F_j^\prime(\bar{u})^{\mu+1-p}\big|v_j\sqrt{Q}\nabla\eta\big|
    \big|\eta\sqrt{Q}\nabla   v_j\big|^{p-1}dx\\
\nonumber\!&\!\!&\!+pq^{p-1}\!\!\fint_{B_r}\!\eta^{p-1}\bar{b}
    F_j^\prime(\bar{u})^{\mu+1-p}  \big|\sqrt{Q}\nabla\eta\big|v_j^p dx\\
\label{12}\!&\!\!&\!+\fint_{B_r}\!c\eta^{p+1-\psi}v_j
    F_j^\prime(\bar{u})^{\mu+1-\psi}\big|\eta\sqrt{Q}\nabla
v_j\big|^{\psi-1}dx\\
\nonumber\!&\!\!&\!+\beta q^{p-1}\fint_{B_r}\!\bar{h}\eta^pv_j^p
F_j^\prime(\bar{u})^{\mu+1-p}dx\\
\nonumber\!&\!\!&\!+q^{p-1}\fint_{B_r}\!\bar{d}\eta^pv_j^p
    F_j^\prime(\bar{u})^{\mu+1-p}\,dx.
\end{eqnarray}

\emph{Step 3.} Our next aim is to apply Young's inequality to the
first and third terms on the right side in order to absorb all
terms containing $|\eta\sqrt{Q}\nabla v_j|$ into the left side. We begin by
noticing that for any $\theta>0$,
$$\big|v_j\sqrt{Q}\nabla\eta\big|\big|\eta\sqrt{Q}\nabla v_j
\big|^{p-1}\leq\theta\big|\eta\sqrt{Q}\nabla
v_j\big|^p+\frac{1}{\theta^{p-1}}\big|v_j\sqrt{Q}\nabla\eta\big|^p.$$
Hence
\begin{eqnarray}
\nonumber\fint_{B_r}\!F_j^\prime(\bar{u})^{\mu+1-p}\big|v_j\sqrt{Q}
    \nabla\eta\big|\big|\eta\sqrt{Q}\nabla
    v_j\big|^{p-1}dx\!\!&\!\!\leq\!\!&\!\!\theta\fint_{B_r}
    F_j^\prime(\bar{u})^{\mu+1-p}\big|\eta\sqrt{Q}\nabla
    v_j\big|^pdx\\
\label{11}\!\!&\!\!\!\!&\!\!+\frac{1}{\theta^{p-1}}\fint_{B_r}
    F_j^\prime(\bar{u})^{\mu+1-p}\big|v_j\sqrt{Q}\nabla
    \eta\big|^{p-1}dx.
\end{eqnarray}

In order to deal with the third term on the right side of
\eqref{12}, we use Young's inequality with exponents
$\frac{p}{\psi-1}$ and $\frac{p}{p+1-\psi}$. Setting $\nu=
\frac{(\mu+1-p)(\psi-1)}{p}$, we get for any $\theta>0$ that
\begin{equation}\label{13}
 \begin{array}{l}
   \displaystyle\fint_{B_r}\!c\eta^{p+1-\psi}v_j
   F_j^\prime(\bar{u})^{\mu+1-\psi}\big|\eta\sqrt{Q}\nabla
     v_j\big|^{\psi-1}dx\\\vspace{0,1cm}
   \displaystyle\hspace{0,6cm}=\fint_{B_r}\!c\eta^{p+1-\psi}v_j
   F_j^\prime(\bar{u})^{\mu+1-\psi-\nu}
     F_j^\prime(\bar{u})^\nu\big|\eta\sqrt{Q}\nabla
   v_j\big|^{\psi-1}\,dx\\\vspace{0,1cm}
   \displaystyle\hspace{0,6cm}\leq\theta\fint_{B_r}\!\!
   F_j^\prime(\bar{u})^{\mu+1-p}\big|\eta\sqrt{Q}\nabla
     v_j\big|^pdx+\frac{1}{\theta^\frac{\psi-1}{p+1-\psi}}\fint_{B_r}
   \!\!c^\frac{p}{p+1-\psi}\eta^p
     v_j^\frac{p}{p+1-\psi}F_j^\prime(\bar{u})^{\mu-
   \frac{\psi-1}{p+1-\psi}}dx.
 \end{array}
\end{equation}
We explicitly note that $\frac{p}{\psi-1}>1$ and $\nu>0$ by
\eqref{ranges1} and since $\mu=p\sigma-1$; see also
Theorem \ref{Gtest-thm}. Moreover $\mu+1-\psi-\nu>0$.

Combining \eqref{11} and \eqref{13} with \eqref{12}, and choosing
$\theta$ suitably small, we obtain
\begin{eqnarray}
\nonumber
    \fint_{B_r}\!F_j^\prime(\bar{u})^{\mu+1-p}\big|\eta\sqrt{Q}\nabla
    v_j\big|^p\,dx\!&\!\leq\!&\!C\left\{\fint_{B_r}\!
    F_j^\prime(\bar{u})^{\mu+1-p}\big|v_j\sqrt{Q}\nabla\eta
    \big|^p\,dx\right.\\
\nonumber\!&\!\!&\!+q^{\gamma-1}\!\!\fint_{B_r}\!\eta^{p-1}\bar{b}
    F_j^\prime(\bar{u})^{\mu+1-p} \big|\sqrt{Q}\nabla\eta\big|v_j^p\,
    dx\\ \label{14}\!&\!\!&\!+\fint_{B_r}\!c^\frac{p}{p+1-\psi}\eta^p
     v_j^\frac{p}{p+1-\psi} F_j^\prime(\bar{u})^{\mu-\frac{\psi
    -1}{p+1-\psi}}dx\\
\nonumber\!&\!\!&\!+\beta  q^{\gamma-1}\fint_{B_r}\!\bar{h}\eta^p
    v_j^\gamma F_j^\prime(\bar{u})^{\mu+1-p}dx\\
\nonumber\!&\!\!&\!\left.+q^{p-1}\fint_{B_r}\!\bar{d}\eta^pv_j^p
    F_j^\prime(\bar{u})^{\mu+1-p}\,dx\right\}
\end{eqnarray}
for a positive constant $C=C(p,a,\sigma)$.

\emph{Step 4.} Now we would like to pass to the limit as $j\to \infty$
in \eqref{14}. By Theorem \ref{Gtest-thm}, both $\{F_j(\bar{t})\}_{j}$
and $\{F_j^\prime(\bar{t})\}_{j}$ are nondecreasing for every
$\bar{t}$. Then the three sequences $v_j=F_j(\bar{u})$, $F^\prime_j(\bar{u})$
and $\big|\sqrt{Q}\nabla v_j\big|= F^\prime_j(\bar{u})\big|\sqrt{Q}
\nabla\bar{u}\big|$ are nondecreasing. Indeed
\begin{equation*}
v_j\nearrow\bar{u}^q,\qquad F_j^\prime(\bar{u})\nearrow
q\bar{u}^{q-1},\qquad\big|\sqrt{Q}\nabla v_j\big|\nearrow
q\bar{u}^{q-1}|\sqrt{Q}\nabla\bar{u}|
\end{equation*}
a.e. in $\Omega$ as $j$ tends to $\infty$. Passing to the limit in
\eqref{14} and using the monotone convergence theorem then yields
\begin{eqnarray}
\nonumber
    q^{\mu+1}\fint_{B_r}\!\bar{u}^{(q-1)(\mu+1)}\big|\eta\sqrt{Q}\nabla
    \bar{u}\big|^p\,dx\!&\!\leq\!&\!C\left\{q^{\mu+1-p}\fint_{B_r}\!\bar{u}^{(q-1)(\mu+1)+p}\big|\sqrt{Q}\nabla\eta
    \big|^p\,dx\right.\\
\nonumber\!&\!\!&\!+q^\mu\!\!\fint_{B_r}\!\eta^{p-1}\bar{b}\big|\sqrt{Q}\nabla\eta\big|\bar{u}^{(q-1)(\mu+1)+p}\,dx\\
\label{15}\!&\!\!&\!+q^{\mu-\frac{\psi-1}{p+1-\psi}}\fint_{B_r}\!c^\frac{p}{p+1-\psi}\eta^p
     \bar{u}^{(q-1)(\mu+1)+\frac{p}{p+1-\psi}}\,dx\\
\nonumber\!&\!\!&\!+\beta
    q^\mu\fint_{B_r}\!\bar{h}\eta^p\bar{u}^{(q-1)(\mu+1)+p}\,dx\\
\nonumber\!&\!\!&\!\left.+q^\mu\fint_{B_r}\!\bar{d}\eta^p\bar{u}^{(q-1)(\mu+1)+p}\,dx\right\}\\
\nonumber\!&\!\!&\! = C\{I + II + III + IV + V\},\nonumber
\end{eqnarray}
where the integrals may not be finite.

\emph{Step 5.} We will estimate II, III, IV and V separately.  Define
\begin{eqnarray}\label{ne1}
(i)\;\;\;\; Y &=& (\mu+1)(q-1) + p,\nonumber\\
(ii) \;\;\;\;\; t&=&\fracc{p}{p+1-\psi}\textrm{ and }\\
(iii)\;\;\; T&=&\mu - \fracc{\psi-1}{p+1-\psi}. \nonumber
\end{eqnarray}

\noindent We begin with term II:

\begin{eqnarray}\label{ne2}
II &=& q^\mu\fint_{B_r}\eta^{p-1}\bar{b}\big|\sqrt{Q}\nabla\eta\big|
\bar{u}^{Y}\;dx\nonumber\\
&=& q^\mu\fint_{B_r}\;\bar{b}\big|\sqrt{Q}\nabla\eta\big|
\bar{u}^{\frac{Y}{p}}(\eta\bar{u}^{\frac{Y}{p}})^{p-1}\;dx\nonumber\\
&\leq& q^\mu\;\avnorm{\bar{b}|\sqrt{Q}\nabla\eta|
\bar{u}^\frac{Y}{p}}{\frac{p\sigma}{p\sigma-p+1}} \,
\avnorm{\eta\bar{u}^\frac{Y}{p}}{p\sigma}^{p-1}
\eea
by H\"older's inequality with $\frac{p-1}{p\sigma} +
\frac{p\sigma-p+1}{p\sigma} = 1$.  Use H\"older's inequality again on the
first factor with $\frac{\sigma}{p(\sigma-1) + 1}+\frac{(\sigma
-1)(p-1)}{p(\sigma -1) + 1} = 1$ and apply Lemma \ref{appendix3} to
the second factor to obtain
\bea
II&\leq&Cq^\mu\; r^{p-1}\;\avnorm{\bar{b}}{\frac{p\sigma}{(\sigma
-1)(p-1)}}\, \avnorm{\bar{u}^\frac{Y}{p}\sqrt{Q}\nabla\eta}{p}\nonumber\\
&&\;\;\times \Big[\avnorm{\bar{u}^\frac{Y}{p}\sqrt{Q}
\nabla\eta}{p}^{p-1}+\Big(\fracc{Y}{p}\Big)^{p-1}
\avnorm{\bar{u}^{\frac{Y}{p}-1}\eta\sqrt{Q}\nabla\bar{u}}{p}^{p-1}
+\fracc{1}{r^{p-1}}\avnorm{\eta\bar{u}^\frac{Y}{p}}{p}^{p-1}\;\Big].\nonumber
\end{eqnarray}
\noindent  Noting that
\begin{eqnarray}
p'\sigma' = \fracc{p\sigma}{(p-1)(\sigma-1)}\nonumber
\end{eqnarray}
and setting
\begin{eqnarray}
\bar{B}=r^{p-1}\avnorm{\bar{b}}{p'\sigma'}\;\; \textrm{      and }\nonumber
\eea
\begin{eqnarray}\label{ne3}
x=\avnorm{\eta\bar{u}^\frac{Y}{p}}{p},\;
y=\avnorm{\bar{u}^\frac{Y}{p}\sqrt{Q}\nabla\eta}{p},\;
z=\avnorm{\bar{u}^{\frac{Y}{p}-1}\eta\sqrt{Q}\nabla\bar{u}}{p},
\eea
we have
\begin{eqnarray}\label{neII}
II &\leq& Cq^\mu\bar{B}y^p +
Cq^\mu\bar{B}\Big(\fracc{Y}{p}\Big)^{p-1}yz^{p-1} +
Cq^\mu\bar{B}y\Big(\fracc{x}{r}\Big)^{p-1}\nonumber\\
&=& Cq^\mu\Big[ \bar{B}y^p + \bar{B}y\Big(\frac{x}{r}\Big)^{p-1}\Big] + Cq^\mu\bar{B}\Big(\frac{Y}{p}\Big)^{p-1}yz^{p-1}.
\end{eqnarray}
We now use Young's inequality on the second and third terms of
(\ref{neII}).  Fix $s_1\in (0,1)$ to be chosen precisely later in the
proof.  Then, since $p'=\frac{p}{p-1}$ and $q\sim Y$,
\begin{eqnarray}
II&\leq&
Cq^\mu\Big[ \bar{B}y^p + \bar{B}^py^p +\Big(\frac{x}{r}\Big)^p\Big] + Cq^{p\mu}\frac{\bar{B}^p}{s_1^p}\Big(\frac{Y}{p}\Big)^{p(p-1)}y^p + s_1^\frac{p}{p-1}z^p\nonumber\\
&\leq& Cq^\mu\Big[ \bar{B}y^p + \bar{B}^py^p +\Big(\frac{x}{r}\Big)^p\Big] + Cq^{p(\mu+p-1)}\frac{\bar{B}^p}{s_1^p}y^p + s_1^\frac{p}{p-1}z^p\nonumber\\
&\leq& Cq^{p(\mu+p-1)}\Big(1+\fracc{\bar{B}}{s_1}\Big)^p\Big(y^p +
\fracc{x^p}{r^p}\Big) + s_1^{\frac{p}{p-1}}z^p\nonumber\\
&\leq&
Cq^{p(\mu+p-1)}\left(1+\left(\fracc{\bar{B}}{s_1}\right)^p\right)\Big(y^p
+ \fracc{x^p}{r^p}\Big) + s_1^{\frac{p}{p-1}}z^p.\nonumber
\end{eqnarray}
We next estimate III.  Fix $\epsilon_1\in(0,1]$ as provided in the
hypothesis.  Then
\begin{eqnarray}\label{ne4}
 III &=& q^T\fint_{B_r}\;c^t\eta^p \bar{u}^{Y+t-p}\;dx\nonumber\\
&=& q^T\fint_{B_r}\;c^t\eta \bar{u}^{\frac{Y}{p}+t-p}(\eta
\bar{u}^\frac{Y}{p})^{p-1}\;dx\nonumber\\
&\leq& q^T\avnorm{c^t\eta \bar{u}^{\frac{Y}{p}+t-p}}{\frac{p
\sigma}{p\sigma-p+1}}\,\avnorm{\eta\bar{u}^\frac{Y}{p}}{p\sigma}^{p-1}\\
\label{ne4.1}&=&q^T\avnorm{c^t\bar{u}^{t-p}(\eta
\bar{u}^\frac{Y}{p})^{\epsilon_1}(\eta\bar{u}^\frac{Y}{p})^{1-
\epsilon_1}}{\frac{p\sigma}{p\sigma-p+1}}\,
\avnorm{\eta\bar{u}^\frac{Y}{p}}{p\sigma}^{p-1}\, ,
\eea
where H\"older's inequality with $\frac{p-1}{p\sigma} +
\frac{p\sigma-p+1}{p\sigma} = 1$ was used to obtain (\ref{ne4}).  Due
to the restrictions on $p,\sigma,\epsilon_1$ we have
\bea (i)&& \frac{p\sigma - p+1}{\sigma\epsilon_1} \geq 1\textrm{
and}\nonumber\\
(ii)&& \frac{\sigma(p-\epsilon_1)+1-p}{1-\epsilon_1}\geq 1.\nonumber
\eea
Therefore, by H\"older's inequality applied to the triple product in
(\ref{ne4.1}),

\begin{eqnarray}
III&\leq&q^T\;\avnorm{c^t\bar{u}^{t-p}}{\frac{p\sigma}{(p-
      \epsilon_1)(\sigma-1)}}\,
      \avnorm{\eta\bar{u}^\frac{Y}{p}}{p}^{\epsilon_1}\,
      \avnorm{\eta\bar{u}^\frac{Y}{p}}{p\sigma}^{p-\epsilon_1}\,.\nonumber
\end{eqnarray}
Now use Lemma \ref{appendix3} on the last factor to obtain
\begin{eqnarray} \label{ne4.2}
III&\leq&Cq^Tr^{p-\epsilon_1}\;\avnorm{c^t\bar{u}^{t-
p}}{\frac{p\sigma}{(p-\epsilon_1)(\sigma-1)}} \,\avnorm{\eta
\bar{u}^\frac{Y}{p}}{p}^{\epsilon_1}\nonumber\\  &&\times
\Big[\avnorm{\bar{u}^\frac{Y}{p}\sqrt{Q}\nabla\eta}{p}^{p-
\epsilon_1}+\Big(\fracc{Y}{p}\Big)^{p-\epsilon_1}
\avnorm{\bar{u}^{\frac{Y}{p}-1}\eta\sqrt{Q}\nabla\bar{u}}{p}^{p-
\epsilon_1}+\fracc{1}{r^{p-\epsilon_1}}\avnorm{\eta
\bar{u}^\frac{Y}{p}}{p}^{p-\epsilon_1}\;\Big].\nonumber
\eea

Setting
\begin{equation}\label{ne5}
\bar{C} =r^p\; \avnorm
{c^t\bar{u}^{t-p}}{\frac{p\sigma}{(p-\epsilon_1)(\sigma-1)}}
\end{equation}
and making the substitutions given by (\ref{ne3}) yields
\begin{eqnarray}\label{neIII}
III&\leq&Cq^T\bar{C}\left(\frac{x}{r}\right)^{\epsilon_1}y^{p-\epsilon_1}+
Cq^T\bar{C}\Big(\fracc{Y}{p}\Big)^{p-\epsilon_1}\left(\frac{x}{r}\right)^{\epsilon_1}z^{p-\epsilon_1}
+ Cq^T\bar{C}\Big(\fracc{x}{r}\Big)^p\nonumber.
\end{eqnarray}
Fix $s_2\in (0,1)$ to be chosen later, and note that
$\Big(\fracc{p}{\epsilon_1}\Big)' = \fracc{p}{p-\epsilon_1}$.  Then,
using Young's inequality in (\ref{neIII}) in a similar manner as in
(\ref{neII}), we have \bea III&\leq&
Cq^\frac{p(T+p-\epsilon_1)}{\epsilon_1}\Big(1+\fracc{\bar{C}}{s_2}\Big)^\frac{p}{\epsilon_1}\Big(y^p+\fracc{1}{r^p}x^p\Big)
+ s_2^\frac{p}{p-\epsilon_1}z^p\nonumber\\
&\leq&Cq^\frac{p(T+p-\epsilon_1)}{\epsilon_1}\left(1+\left(\fracc{\bar{C}}{s_2}\right)^\frac{p}{\epsilon_1}\right)
\Big(y^p+\fracc{1}{r^p}x^p\Big)+s_2^\frac{p}{p-\epsilon_1}z^p.\nonumber
\end{eqnarray}

Terms IV and V are estimated in the same way as $III$.
Fix $\epsilon_2,\epsilon_3\in (0,1]$. For $s_3,s_4\in (0,1)$ to
be chosen later, and with
\begin{eqnarray}\label{ne7}
(i)\;\;\;\; \bar{H} &=& r^p\avnorm{\bar{h}}{\frac{p\sigma}{(p-\epsilon_2)(\sigma-1)}}\,\textrm{ and }\\
(ii)\;\;\;\bar{D} &=&
r^p\avnorm{\bar{d}}{\frac{p\sigma}{(p-\epsilon_3)(\sigma-1)}},\nonumber
\end{eqnarray}
we have
\begin{eqnarray}\label{neIV}
IV &\leq&
Cq^\frac{p(\mu+1+p-\epsilon_2)}{\epsilon_2}\left(1+\left(\fracc{\bar{H}}{s_3}\right)^\frac{p}{\epsilon_2}\right)\Big(y^p
+ \fracc{1}{r^p}x^p\Big) +s_3^\frac{p}{p-\epsilon_2}z^p
\end{eqnarray}
and
\begin{eqnarray}\label{neV}
V &\leq& Cq^\frac{p(\mu
+p-\epsilon_3)}{\epsilon_3}\left(1+\left(\fracc{\bar{D}}{s_4}\right)^\frac{p}{\epsilon_3}\right)\Big(y^p
+ \fracc{1}{r^p}x^p\Big) + s_4^\frac{p}{p-\epsilon_3}z^p.
\end{eqnarray}

\emph{Step 6.} Set $\epsilon = \min
\{\epsilon_1,\;\epsilon_2,\;\epsilon_3\}$ and choose
$s_1,s_2,s_3,s_4$ depending only on
$p,\epsilon_1,\epsilon_2,\epsilon_3$ so that
\begin{equation*}
s_1^{\frac{p}{p-1}}+s_2^\frac{p}{p-\epsilon_1}+
s_3^\frac{p}{p-\epsilon_2} +s_4^\frac{p}{p-\epsilon_3} \leq
\fracc{1}{2C}
\end{equation*}
where $C$ is as in \eqref{15}.  Then there is a $c_* = c_*(\epsilon,
p,\sigma) >0$ so that \eqref{15} becomes
\begin{eqnarray}
z^p &\leq& Cq^{c_*p}\left( 1+ \left(\fracc{\bar{B}}{s_1}\right)^p +
\left(\fracc{\bar{C}}{s_2}\right)^\frac{p}{\epsilon_1}
+\left(\fracc{\bar{H}}{s_3}\right)^\frac{p}{\epsilon_2} +
\left(\fracc{\bar{D}}{s_4}\right)^\frac{p}{\epsilon_3}\right)\Big(y^p
+ \fracc{x^p}{r^p}\Big)\nonumber\\
&\leq& Cq^{c_*p}\left( 1+ \bar{B} +{\bar{C}}^\frac{1}{\epsilon_1}
+{\bar{H}}^\frac{1}{\epsilon_2} +
{\bar{D}}^\frac{1}{\epsilon_3}\right)^p\Big(y^p
+ \fracc{x^p}{r^p}\Big),\nonumber
\end{eqnarray}
where $C$ now also depends on $s_1,s_2,s_3$ and $s_4$.  Taking
$p^{th}$ roots and setting $\bar{Z}=\Big(1+ \bar{B} +
\bar{C}^\frac{1}{\epsilon_1} +\bar{H}^\frac{1}{\epsilon_2} +
\bar{D}^\frac{1}{\epsilon_3}\Big)$ as in \eqref{Zetabar}, it follows
that
\begin{equation}\label{ne16}
\avnorm{\bar{u}^{\frac{Y}{p}-1}\eta\sqrt{Q}\nabla\bar{u}}{p} \leq
Cq^{c_*}\bar{Z}\Big(\avnorm{\bar{u}^\frac{Y}{p}\sqrt{Q}\nabla\eta}{p}
+ \fracc{1}{r}\avnorm{\eta\bar{u}^\frac{Y}{p}}{p}\Big).
\end{equation}
Using Lemma \ref{appendix3} and estimate
\eqref{ne16} we have, with $b_*=c_*+1$,
\begin{eqnarray}
\avnorm{\eta\bar{u}^\frac{Y}{p}}{p\sigma} &\leq& Cr\avnorm {\sqrt{Q}\nabla\Big(\eta\bar{u}^\frac{Y}{p}\Big)}{p} +
    C\avnorm{\eta\bar{u}^\frac{Y}{p}}{p}\nonumber\\
\label{ne17}&\leq&Cr\Big\{\avnorm{\bar{u}^\frac{Y}{p}\sqrt{Q}\nabla\eta}{p}+
    \Big(\frac{Y}{p}\Big)\avnorm{\bar{u}^{\frac{Y}{p}-1}
    \eta\sqrt{Q}\nabla\bar{u}}{p}\Big\} + C\avnorm{\eta\bar{u}^\frac{Y}{p}}{p}\\
&\leq&Cq^{b_*}\bar{Z}\Big\{\avnorm{\eta\bar{u}^\frac{Y}{p}}{p}+r\avnorm{\bar{u}^\frac{Y}{p}\sqrt{Q}\nabla\eta}{p}\Big\}.\nonumber
\end{eqnarray}

\emph{Step 7.} We now use the accumulating sequence of Lipschitz
cutoff functions in \eqref{cutoff}. For any
$\Omega^\prime\subset\Omega,$ let $\chi_{_{\Omega^\prime}}:\Omega
\rightarrow\R$ be defined by
$$\chi_{_{\Omega^\prime}}(x)=\begin{cases}
1\qquad\text{if }x\in \Omega^\prime\\
0\qquad\text{if }x\notin \Omega^\prime.
\end{cases}$$ For each $j$, let $S_j=supp\;\eta_j$ and recall that
$\eta_j = 1$ on $S_{j+1}$. Since $s^*>p\sigma'$, there exists $s$ so
that $1\leq s<\sigma$ and $s'p=s^*$.  Then for each $j$,
\begin{eqnarray}\label{ne18}
\avnorm{\bar{u}^\frac{Y}{p}\chi_{_{S_{j+1}}}}{p\sigma} &\leq&
Cq^{b_*}\bar{Z}\Big\{\avnorm{\bar{u}^\frac{Y}{p}\chi_{_{S_j}}}{p}
 +r \avnormb{\sqrt{Q}\nabla\eta_j}{s^*}{B_r}\avnorm{\bar{u}^\frac{Y}{p}\chi_{_{S_j}}}{sp}\Big\}\nonumber\\
&\leq&Cq^{b_*}\bar{Z}N^j\avnorm{\bar{u}^\frac{Y}{p}\chi_{_{S_j}}}{sp}.
\end{eqnarray}
Rewriting \eqref{ne18} so that $\bar{u}$ appears to power $1$ inside
each norm, we see that
\begin{equation}\label{ne19}
\avnorm{\bar{u}\chi_{_{S_{j+1}}}}{\sigma Y} \leq
\Big(C\bar{Z}N^jq^{b_*}\Big)^\frac{p}{Y}\avnorm
{\bar{u}\chi_{_{S_j}}}{sY}.
\end{equation}
Note that $\sigma Y > sY$ since $1\leq s<\sigma$.  Thus, $\bar{u} \in
L^{sY}(S_j)$ implies the stronger inclusion $\bar{u}\in L^{\sigma
Y}(S_{j+1})$.  We will use this fact and a M\"oser iteration to
obtain the conclusion of the Theorem \ref{main_thm}.  Set $\mathcal{X} =
\frac{\sigma}{s} > 1$ and let $q_0=1$.  For each
$j\in\mathbb{N}$, choose $q_j> 1$ so that $Y_j=(\mu +1)(q_j-1)+p$
and $Y_j = p\mathcal{X}^j$.  With $Y_0 = p,$ we have
\begin{equation}\label{ne23}
Y_j = p{\mathcal X}^j \quad\mbox{for $j\geq 0$}.
\end{equation} Choosing $Y=Y_j$ in \eqref{ne19} gives
\begin{equation}\label{ne20}
\avnorm{\bar{u}\chi_{_{S_{j+1}}}}{sp\mathcal{X}^{j+1}} \leq
(C\bar{Z})^{\mathcal{X}^{-j}}N^{j\mathcal{X}^{-j}}q_j^{b_*\mathcal{X}^{-j}}\avnorm{\bar{u}\chi_{_{S_j}}}{sp\mathcal{X}^j}.
\end{equation}

\noindent Let $\Psi_1 = \sum_{j=1}^{\infty} \mathcal{X}^{-j}$ and
$\Psi_2 = \sum_{j=1}^{\infty} j\mathcal{X}^{-j}$, recalling that
$\mathcal{X}>1$. Then, since $q_j\backsim \mathcal{X}^j$ and
$B_{\tau r}\subset S_{j}$ for each $j\in\mathbb{N}$, we obtain
\begin{equation}\label{ne21}
\avnorm{\bar{u}\chi_{_{B_{\tau r}}}}{sp\mathcal{X}^{j+1}}
\,\leq\,(C\bar{Z})^{\Psi_1}(N\mathcal{X}^{b_*})^{\Psi_2}
\avnorm{\bar{u}\chi_{_{S_1}}}{sp\mathcal{X}}.
\end{equation}

Let $\eta_0\in C^\infty_0(B_r)$ be a nonnegative cutoff function
so that $S_1\subset \{x:\eta_0(x) = 1\}$ and $\eta_0\leq1$
in $B_r$. Then since $Y_0=p$ and $q_0 =1$, \eqref{ne17} and
\eqref{ne21} imply that
\begin{equation}\label{ne21.1}
\avnorm{\bar{u}\chi_{_{B_{\tau
r}}}}{sp\mathcal{X}^{j+1}}\,\leq\,C\bar{Z}^{\Psi_1+1}(N
\mathcal{X}^{b_*})^{\Psi_2}\avnormb{\bar{u}}{sp}{B_r}.
\end{equation}
Since this holds for every $j\in\mathbb{N}$, it follows that
\begin{equation}\label{ne22}
 \norm \bar{u}\norm_{L^\infty(B_{\tau
r})}\leq C\bar{Z}^{\Psi_0}\avnormb{\bar{u}}{sp}{B_r}
\end{equation}
where $\Psi_0 = \disp\sum_{j=0}^\infty \mathcal{X}^{-j}$ and $C$ are
independent of $u,r,b,c,d,e,f,g,h$. This completes the proof of
Theorem \ref{main_thm}.
$\qquad\Box$


\section{\textbf{Appendix 1}}

In Appendix 1, we will prove some facts used in the proof of Theorem
\ref{main_thm} that are related to the chain rule and the
iteration process. See also \cite{SW2} for results related to the
chain rule.

\begin{lem}\label{appendix2}
Let $(u, \nabla u) \in\waQ$ and $\phi\in\mathcal{C}^1(\R)$ with
$\phi^\prime\in L^\infty(\R)$. Then $\phi(u)\in\waQ$ and
$$\sqrt{Q}\nabla\big(\phi(u)\big)=\phi^\prime(u)\sqrt{Q}\nabla u
\qquad\text{a.e. in }\Omega.$$
\end{lem}
\textbf{Proof:} Let $\{u_j\}_{j\in\N}\subset\waQ \cap
\text{Lip}_{\text{loc}}(\Omega)$  be a
representative sequence for $(u, \nabla u)\in\waQ$. Then as usual, up
to subsequences,
\begin{eqnarray}
\label{25}&&u_j\rightarrow u\quad\text{a.e. in }\Omega\text{ and in
$L^p(\Omega)$, and}\\
\nonumber&&\sqrt{Q}\nabla u_j\rightarrow\sqrt{Q}\nabla
   u\quad\text{a.e. in }\Omega\text{ and in }\big(L^p(\Omega)\big)^n.
\end{eqnarray}

\emph{Claim 1:} $\phi(u_j)\in\text{Lip}_{\text{loc}}(\Omega)$ for
every $j$. Indeed, by the fundamental theorem of calculus, for
every $x,y\in\Omega$,
$$\big|\phi\big(u_j(x)\big)-\phi\big(u_j(y)\big)\big|=\left|
\int_{u_j(y)}^{u_j(x)}\phi^\prime(t)\,dt\right|\leq
\|\phi^\prime\|_\infty|u_j(x)-u_j(y)|,$$ and the claim follows from
$u_j\in\text{Lip}_{\text{loc}}(\Omega)$.

\emph{Claim 2:} $\phi(u_j)\rightarrow\phi(u)$ in $L^p(\Omega)$ and
a.e. in $\Omega$. In fact, since $u_j\rightarrow u$ a.e. in $\Omega$ and
$\phi$ is continuous, then $\phi(u_j)\to \phi(u)$ a.e. in
$\Omega$. Also,
$$\big|\phi(t)\big|=\left|\phi(0)+\int_0^t\phi^\prime(s)\,ds
\right|\leq|\phi(0)|+\|\phi^\prime\|_\infty|t|:=B_0+A_0|t|.$$
Then
\begin{eqnarray}
\nonumber\big|\phi(u_j)-\phi(u)\big|^p&\leq&2^p\left(\big|\phi(u_j)\big|^p
+\big|\phi(u)\big|^p\right)\\
\nonumber&\leq&4^p\left(2B_0^p+A_0^p|u_j|^p+A_0^p|u|^p\right)\\
\nonumber&\leq&c(1+|u_j|^p+|u|^p).
\end{eqnarray}
Since $1+|u_j|^p+|u|^p\rightarrow1+2|u|^p$ a.e. in $\Omega$ and in
$L^1(\Omega)$ by relation \eqref{25}, by Lebesgue's Sequentially
Dominated Convergence Theorem implies that
$$\int_\Omega\big|\phi(u_j)-\phi(u)\big|^p\,dx\rightarrow0.$$

\emph{Claim 3:}
$\sqrt{Q}\nabla\big(\phi(u_j)\big)\rightarrow\phi^\prime(u)\sqrt{Q}\nabla
u$ a.e. in $\Omega$ and in $\big(L^p(\Omega)\big)^n$. Indeed, since
$u_j$ and $\phi(u_j)$ are locally Lipschitz, their gradients exist
a.e. in $\Omega$, and
$$\sqrt{Q}\nabla\big(\phi(u_j)\big)=\phi^\prime(u_j)\sqrt{Q}\nabla
u_j\quad\text{a.e. in }\Omega.$$ Since $\phi^\prime$ is continuous,
\eqref{25} gives
$\sqrt{Q}\nabla\big(\phi(u_j)\big)\rightarrow\phi^\prime(u)\sqrt{Q}\nabla
u$ a.e. in $\Omega$. Moreover,
\begin{eqnarray}
\nonumber\big|\sqrt{Q}\nabla\big(\phi(u_j)\big)-\phi^\prime(u)\sqrt{Q}\nabla
    u\big|^p&\leq&2^p\left(\big|\sqrt{Q}\nabla\big(\phi(u_j)\big)\big|^p+\|\phi^\prime\|_\infty^p\big|\sqrt{Q}\nabla
    u\big|^p\right)\\
\nonumber&\leq&2^p\|\phi^\prime\|_\infty^p\left(\big|\sqrt{Q}\nabla
    u_j\big|^p+\big|\sqrt{Q}\nabla u\big|^p\right),
\end{eqnarray}
and $\big|\sqrt{Q}\nabla u_j\big|^p+\big|\sqrt{Q}\nabla
u\big|^p\rightarrow2\big|\sqrt{Q}\nabla u\big|^p$ a.e. in $\Omega$
and in $L^1(\Omega)$ by \eqref{25}. Again by Lebesgue's
Sequentially Dominated Convergence Theorem,
$$\int_\Omega\big|\sqrt{Q}\nabla\big(\phi(u_j)\big)-\phi^\prime(u)
\sqrt{Q}\nabla u\big|^p\,dx\rightarrow0.$$
Therefore, the sequence of locally Lipschitz functions
$\{\phi(u_j)\}_{j\in\N}\subset\waQ$ is Cauchy in $\waQ$ and defines
an element of $\waQ$ having $\phi(u)$ as its $L^p$--part and
$\phi^\prime(u)\sqrt{Q}\nabla u$ as its gradient--part. This
completes the proof.$\qquad\Box$

\begin{lem}\label{appendix1}
Let $u\in\waQ$ and $k\in\R$. Then there exists $\bar{u}\in\waQ$
whose $L^p$--part is $|u|+k$ and whose gradient--part $\nabla\bar{u}$
satisfies
\begin{equation}\label{26}
\sqrt{Q(x)}\nabla\bar{u}(x)=\begin{cases}
                               \,\sqrt{Q(x)}\nabla u(x)\,\,\,\qquad\text{if }u(x)\geq0,\\
                               -\sqrt{Q(x)}\nabla u(x)\qquad\text{if }u(x)<0
                            \end{cases}\qquad\text{a.e. in }\Omega.
\end{equation}
\end{lem}
\begin{rem}
Choosing $k=0$, it follows that if $u\in \waQ$ then $|u|\in\waQ$, and
the gradient--part of $|u|$ satisfies \eqref{26}.
\end{rem}
\begin{rem}
For $u, k$ and $\bar{u}$ as in Lemma \ref{appendix1},
$|\sqrt{Q}\nabla\bar{u}|=|\sqrt{Q}\nabla u|$ a.e. in $\Omega$.
\end{rem}

\textbf{Proof of Lemma \ref{appendix1}.}  For any $\theta>0$, define
$\phi_\theta:\R\rightarrow\R$ by $\phi_\theta(t)
=(t^2+\theta^2)^\frac{1}{2}.$  Then
\begin{itemize}
\item $\phi_\theta\in\mathcal{C}^1(\R)$ with
$\phi^\prime_\theta(t)=t/(t^2+\theta^2)^\frac{1}{2}$ and
$\phi^\prime_\theta\in L^\infty(\R)$,
\item $0\leq\phi_\theta(t)\leq|t|+\theta$ and
$|\phi_\theta^\prime(t)|\leq1$ for $t\in\R$,
\item as $\theta\rightarrow 0,$  $\phi_\theta(t)\rightarrow|t|$ and
$\displaystyle\phi_\theta^\prime(t)\rightarrow\text{sign}(t):=\begin{cases}
1\qquad\text{ if }t>0,\\
0\qquad\text{ if }t=0,\\
-1\qquad\!\text{ if }t<0.
\end{cases}$
\end{itemize}

Now let $u\in\waQ$ and $\lambda\geq0$ be such that
$\big|\big\{x\in\Omega:|u(x)|=\lambda\big\}\big|=0$. Let
$\{u_j\}_{j\in\N}\subset\waQ\cap \text{Lip}_{\text{loc}}(\Omega)$ be a
representative sequence for $(u,\nabla u)$. We may assume as usual that
\begin{eqnarray}
\nonumber&&u_j\rightarrow u\quad\text{a.e. in }\Omega\text{ and in
$L^p(\Omega)$, and}\\
\nonumber&&\sqrt{Q}\nabla u_j\rightarrow\sqrt{Q}\nabla
   u\quad\text{a.e. in }\Omega\text{ and in }\big(L^p(\Omega)\big)^n.
\end{eqnarray}
Keeping $\lambda$ fixed, let $\varphi_{j,\lambda}(x)
=\phi_\frac{1}{j}\big(u_j(x)+\lambda\big)$ for $x\in\Omega$. As
shown in the proof of Lemma \ref{appendix2}, $\varphi_{j,\lambda}
\in\text{Lip}_{\text{loc}}(\Omega)$ for every $j$ and $\varphi_{j,
\lambda}\rightarrow|u+\lambda|$ a.e. in $\Omega$ and in $L^p(\Omega)$
as $j \to \infty$. Moreover, since $u\neq-\lambda$ a.e. in $\Omega$,
\begin{equation*}
\sqrt{Q}\nabla\varphi_{j,\lambda}(x)\rightarrow V_\lambda(x):=
  \begin{cases}
     \sqrt{Q}\nabla u(x)\,\,\,\,\qquad\text{if }u(x)\geq-\lambda,\\
     -\sqrt{Q}\nabla u(x)\qquad\text{if }u(x)<-\lambda
  \end{cases}
\end{equation*}
a.e. in $\Omega$ and in $\big(L^p(\Omega)\big)^n$ as $j \to
\infty$. Then $\{\varphi_{j,\lambda}\}_{j\in\N}\subset\waQ$ is a
Cauchy sequence and thus defines an element $\varphi_\lambda\in\waQ$
having $|u+\lambda|$ as its $L^p$--part and $V_\lambda$ as its
gradient--part.

In case $\big|\big\{x\in\Omega:|u(x)|=0\big\}\big|=0,$ we choose
$\lambda=0$ in the preceding argument and conclude the proof of the
lemma. In case $\big|\big\{x\in\Omega:|u(x)|=0\big\}\big|>0,$ choose a
sequence $\lambda_m\searrow0$ such that (see Corollary
\ref{goodlevelsets2})
\begin{equation}\label{27}
\big|\big\{x\in\Omega:|u(x)|=\lambda_m\big\}\big|=0\qquad\text{ for all
}m.
\end{equation}
Then $\varphi_{\lambda_m}=|u+\lambda_m|\rightarrow|u|$ a.e. in
$\Omega$, and since $$\int_{\Omega}\big||u+\lambda_m|-|u|\big|^p\,dx
\leq \lambda_m^p|\Omega|,$$ we also have that
$\varphi_{\lambda_m}\rightarrow|u|$ in $L^p(\Omega)$.

Let us show that for a.e. $x\in\Omega$,
\begin{equation*}
V_{\lambda_m}(x)\rightarrow V(x):=
  \begin{cases}
    \sqrt{Q}\nabla u(x)\,\,\,\,\qquad\text{if }u(x)\geq0,\\
    -\sqrt{Q}\nabla u(x)\qquad\text{if }u(x)<0.
  \end{cases}
\end{equation*}
Indeed, if $u(x)\geq0$ then $u(x)\geq-\lambda_m$, and hence
$V_{\lambda_m}(x)=\sqrt{Q}\nabla u(x)$ for all $m\in\N$. On the
other hand, if $u(x)<0$ then $u(x)<-\lambda_m$ for all large $m$, and then
$V_{\lambda_m}(x)=-\sqrt{Q}\nabla u(x)$, again for all large $m$.

Since for all $m$,
$$|V_{\lambda_m}-V|^p\leq2^p\big(|V_{\lambda_m}|^p+
|V|^p\big)=2^{p+1}|\sqrt{Q}\nabla u|^p,$$
Lebesgue's Dominated Convergence Theorem yields that
$V_{\lambda_m}\rightarrow V$ in $\big(L^p(\Omega)\big)^n$. Then
$\{\varphi_{\lambda_m}\}_{m\in\N}$ is a Cauchy sequence in $\waQ$, and
it converges to an element of $\waQ$ having $|u|$ as its $L^p$--part
and $V$ as its gradient--part. We denote this element by $|u|$.
Finally, since $|u|\in\waQ$ and $k\in\waQ$ with $\sqrt{Q}\nabla k =0$
for every $k\in\R$, we obtain that $\bar{u}=|u|+k\in\waQ$ and
\eqref{26} holds. The proof of the lemma is now complete.$\qquad\Box$

\begin{lem}\label{appendix3}
Let $k\geq0$, $u\in\waQ$ be a weak solution of \eqref{eqdiff}, and
$\bar{u}=|u|+k$.  Let $\eta\in\text{Lip}_0(\Omega)$ and
$\text{supp}(\eta)\subset B$ for a $\rho$-ball $B=B(y,r)$ with
$r<r_1(y)$. Let \eqref{Cond1} and Sobolev's inequality
\eqref{Sobolev} be true. Suppose that \eqref{E3} holds for $B$ and
$\eta$, and also that, for $t\ge 1$ as in \eqref{E3}, condition
\eqref{3.6-0} holds for $B$ with $t'$ given by $1/t+1/t'=1$. If
$q\geq1$, $p>1$, $\sigma>1$ and $\theta>0$, then
\begin{eqnarray}
\nonumber\left(\fint_{B}\big|\eta\bar{u}^\frac{(q-1)\theta+p}{p}\big|^{p\sigma}\,dx\right)^\frac{1}{p\sigma}
     \!\!\!&\!\!\leq\!\!&\!Cr\!\left[\!\frac{(q-1)\theta+p}{p}\left(\fint_{B}\!\big|\eta\bar{u}^\frac{(q-1)\theta}{p}\sqrt{Q}\nabla\bar{u}\big|^p
     \,dx\right)^\frac{1}{p}\right.\\
\label{16}\!\!&\!\!\!\!&\!\!+\left.\left(\fint_{B}\!\big|\bar{u}^\frac{(q-1)\theta+p}{p}\sqrt{Q}\nabla\eta\big|^pdx
     \right)^{\!\frac{1}{p}}\right]\!+\!C\!
    \left(\fint_{B}\!\big|\eta\bar{u}^\frac{(q-1)\theta+p}{p}\big|^pdx\!\right)^{\!\frac{1}{p}},
\end{eqnarray}
where the integrals may not be finite.
\end{lem}
\textbf{Proof:} For any $l>0,$ let
\begin{equation*}
H_l(t) =
 \begin{cases}
  \frac{(q-1)\theta+p}{p}l^\frac{(q-1)\theta}{p}t+\frac{(q-1)\theta}{p}l^\frac{(q-1)\theta+p}{p}\hspace{1,3cm}\text{if }t<-l,\\
  |t|^\frac{(q-1)\theta+p}{p}\text{sign}(t)\hspace{4cm}\text{if }|t|\leq l,\\
  \frac{(q-1)\theta+p}{p}l^\frac{(q-1)\theta}{p}t-\frac{(q-1)\theta}{p}l^\frac{(q-1)\theta+p}{p}\hspace{1,4cm}\text{if }t>l.\\
 \end{cases}
\end{equation*}
Then $H_l\in\mathcal{C}^1(\R)$ with
\begin{equation*}
H_l^\prime(t)=
 \begin{cases}
  \frac{(q-1)\theta+p}{p}l^\frac{(q-1)\theta}{p}\hspace{1,3cm}\text{if }t<-l,\\
  \frac{(q-1)\theta+p}{p}|t|^\frac{(q-1)\theta}{p}\hspace{1,1cm}\text{if }|t|\leq l,\\
  \frac{(q-1)\theta+p}{p}l^\frac{(q-1)\theta}{p}\hspace{1,4cm}\text{if }t>l,\\
 \end{cases}
\end{equation*}
and $H^\prime_l\in L^\infty(\R)$ with
$\|H^\prime_l\|_\infty\leq\frac{(q-1)\theta+p}{p}l^\frac{(q-1)\theta}{p}$.
Notice that $H_l^\prime(t)$ is nondecreasing in $l$ for every
$t\in\R$, while $H_l(t)$ is nondecreasing in $l$ only for $t\geq0$.

By Lemmas \ref{appendix1} and \ref{appendix2}, $H_l(\bar{u})\in\waQ$ with
$\sqrt{Q}\nabla\big(H_l(\bar{u})\big)=H^\prime_l(\bar{u})\sqrt{Q}\nabla\bar{u}$
a.e. in $\Omega$. Then, by Proposition \ref{C4} and the assumptions on
$\eta$, $\eta H_l(\bar{u})\in\wbQ$ with support in $B$ and
$$\sqrt{Q}\nabla\big(\eta H_l(\bar{u})\big)=\eta
H^\prime_l(\bar{u})\sqrt{Q}\nabla\bar{u}+H_l(\bar{u})\sqrt{Q}\nabla\eta.$$
By Sobolev's inequality \eqref{Sobolev},
\begin{equation*}
\begin{array}{l}
\displaystyle\left(\fint_{B}\!\big|\eta
H_l(\bar{u})\big|^{p\sigma}dx\right)^\frac{1}{p\sigma}\leq\\
\hspace{0,3cm}\displaystyle Cr\!\left[\left(\fint_{B}\!\big|\eta
H^\prime_l(\bar{u})\sqrt{Q}\nabla\bar{u}\big|^pdx\right)^\frac{1}{p}
\!\!+\!\left(\fint_{B}\!\big|H_l(\bar{u})\sqrt{Q}\nabla\eta\big|^pdx\right)^\frac{1}{p}\right]\!+\!
C\!\left(\fint_{B}\!\big|\eta
H_l(\bar{u})\big|^pdx\!\right)^{\frac{1}{p}}.
\end{array}
\end{equation*}
Since $\bar{u}\geq0$ in $\Omega$, both $H_l(\bar{u})$ and
$H^\prime_l(\bar{u})$ are nondecreasing in $l$ and
$$H_l(\bar{u})\nearrow\bar{u}^\frac{(q-1)\theta+p}{p},\qquad
H^\prime_l(\bar{u})\nearrow\frac{(q-1)\theta+p}{p}\bar{u}^\frac{(q-1)\theta}{p}$$
a.e. in $\Omega$ as $l\to \infty$. Passing to the limit in the
previous inequality and using the monotone convergence theorem, we get
\begin{eqnarray}
\nonumber\left(\fint_{B}\!\big|\eta\bar{u}^\frac{(q-1)\theta+p}{p}\big|^{p\sigma}dx\right)^\frac{1}{p\sigma}
    \!\!&\!\!\leq\!\!&\!\!Cr\!\left[\frac{(q-1)\theta+p}{p}\left(\fint_{B}\!\big|\eta
    \bar{u}^\frac{(q-1)\theta}{p}\sqrt{Q}\nabla\bar{u}\big|^pdx\right)^\frac{1}{p}\right.\\
\nonumber\!\!&\!\!\!\!&\!\!\hspace{0,2cm}+\left.\left(\fint_{B}\!\big|\bar{u}^\frac{(q-1)\theta+p}{p}
    \sqrt{Q}\nabla\eta\big|^pdx\!\right)^\frac{1}{p}\!\right]\!+\!C\!
    \left(\fint_{B}\!\big|\eta\bar{u}^\frac{(q-1)\theta+p}{p}\big|^pdx\!\right)^\frac{1}{p},
\end{eqnarray}
where the integrals may not be finite. This completes the
proof. $\qquad\Box$

\section{\textbf{Appendix 2}}

In Appendix 2, we will prove the following three theorems related to
the structural assumptions about equation \eqref{eqdiff}

\begin{thm}\label{32} Consider the differential equation \eqref{eqdiff}:
\begin{equation*}
\text{div}\big(A(x,u,\nabla u)\big)=B(x,u,\nabla u).
\end{equation*}
Suppose that the structural assumptions \eqref{struct} hold relative
to a symmetric nonnegative definite matrix $Q(x)$. If $H(x)$ is
another symmetric nonnegative definite matrix and
\begin{equation}\label{equivalent}
\frac{1}{C}\langle Q(x)\xi,\xi\rangle\leq\langle
H(x)\xi,\xi\rangle\leq C\langle Q(x)\xi,\xi\rangle\qquad\text{for
all }\xi\in\R^n,\,\text{ a.e. }x\in\Omega,
\end{equation}
then there is a vector ${\hat A}(x,z,\xi)$ such that
\begin{eqnarray}
\nonumber A(x,z,\xi)&=&\sqrt{H(x)}{\hat A}(x,z,\xi),\\
\nonumber\xi \cdot A(x,z,\xi) &\geq& \big(C^\frac{p}{2}
a\big)^{-1}\Big|\sqrt{H(x)}\cdot \xi\Big|^p - h|z|^\gamma - g,\\
\label{struct2bis} \Big|\hat{A}(x,z,\xi)\Big| &\leq&
\big(C^\frac{p}{2} a\big)\Big|\sqrt{H(x)}\cdot\xi\Big|^{p-1} +
\big(C^\frac{1}{2}b\big)|z|^{\gamma -1} + \big(C^\frac{1}{2}e\big),\\
\nonumber\Big|B(x,z,\xi)\Big|&\leq&\big(C^\frac{\psi-1}{2}\,c\big)\Big|\sqrt{H(x)}\cdot\xi\Big|^{\psi-1}+d|z|^{\delta-1}+f
\end{eqnarray}
for $\xi\in\R^n,z\in\R$ and a.e $x\in\Omega$. Here, $C$ is the same
constant as in \eqref{equivalent}, and $a, b, c, d, e, f, g, h$ are as
in \eqref{struct}.
\end{thm}

\vspace{0,1cm}

Next we will show that many linear equations satisfy \eqref{struct}.
Consider the linear equation
\begin{equation}\label{34}
\text{div}\big(Q(x)\nabla u\big)+\mathbf{HR}u+\mathbf{S^\prime
G}u+Fu=f+\mathbf{T^\prime g}\qquad\text{in }\Omega,
\end{equation}
where $Q(x)$ is symmetric and nonnegative definite,
$\mathbf{R}=\{R_i\}_{i=1}^n$, $\mathbf{S}=\{S_i\}_{i=1}^n$,
$\mathbf{T}=\{T_i\}_{i=1}^n$ are collections of vector fields
subunit with respect to $Q(x)$, and where the operator coefficients
$\mathbf{H}=\{H_i\}_{i=1}^n$, $\mathbf{G}=\{G_i\}_{i=1}^n$ and $F$ as
well as the inhomogeneous data $\mathbf{g}=\{g_i\}_{i=1}^n$ and $f$ are
measurable. See also \cite{SW1}. We will prove the following fact
about such equations.

\begin{thm}\label{35} The linear equation \eqref{34} satisfies the structural
conditions \eqref{struct} with $p=\gamma=\psi=\delta=2$ relative
to the matrix $Q(x)$.
\end{thm}

Finally, we will prove the next result concerning conditions
\eqref{struct} and \eqref{struct3}.

\begin{thm}\label{36}
For the differential equation \eqref{eqdiff}, the structural
assumptions \eqref{struct} are satisfied if and only if
\eqref{struct3} is satisfied.
\end{thm}

\vspace{0,1cm}

For the proofs, we will need some technical results which we collect
in the following lemmas. We state the first two without proofs.

\textbf{Notation:} For any $k\in\N$, we will denote the
identity $k\times k$ matrix by $I_k$ and the zero $k\times k$ matrix
by $0_k$. Also, $\langle\cdot,\cdot\rangle_{\R^k}$ and
$|\cdot|_{\R^k}$ will denote respectively the inner product and the norm in
$\R^k$. When we work in $\R^n$, i.e., when $k=n$, we will usually omit
the subscript $\R^k$.

Also, let $\text{Mat}(n,\R)$ be the set of $n\times n$ real
matrices, $\mathcal{O}(n)$ be the set of $n\times n$ real
orthogonal matrices, and $\mathcal{S}^{n\times n}=\text{Symm}(n,\R)$
be the set of symmetric $n\times n$ real matrices.
For any $Q\in\mathcal{S}^{n\times n}$, we will write $Q\geq0$
if $Q$ is nonnegative definite. Since $Q$ is symmetric and hence
diagonalizable, the condition $Q\geq 0$ is the same as assuming that
$Q$ has nonnegative eigenvalues. Finally, if
$Q,H\in\mathcal{S}^{n\times n}$, we will write $Q\geq H$ if
$Q-H\geq0$, i.e. if $\langle Q\xi,\xi\rangle\geq\langle
H\xi,\xi\rangle$ for all $\xi\in\R^n$.

\begin{lem}\label{37}
(i) If $Q,H\in\mathcal{S}^{n\times n}$ with $Q, H \geq 0$,
then $Q+H\geq 0.$

(ii) If $Q\in\mathcal{S}^{n\times n}$, $Q\geq 0$ and
$M\in\text{Mat}(n,\R)$, then $M^TQM\geq 0$, where $M^T$ denotes the
transpose of $M$.

(iii) If $Q\in\mathcal{S}^{n\times n}$, $Q\geq 0$ and
$\text{det}\,Q\neq0$, then $Q^{-1}\geq 0.$
\end{lem}

\begin{lem}\label{38} If $Q\in\mathcal{S}^{n\times n}$ with $Q\geq
0$, there is a unique matrix $\sqrt{Q} \in \mathcal{S}^{n\times n}$
such that $\sqrt{Q} \geq 0$ and $\sqrt{Q}\sqrt{Q}=Q$. Moreover

i) $\lambda\geq0$ is an eigenvalue for $Q$ with eigenvector
$v$ if and only if $\sqrt\lambda$ is an eigenvalue for $\sqrt{Q}$ with
eigenvector $v$.

ii) $\sqrt{Q}$ is invertible if and only if $Q$ is invertible.
\end{lem}

\begin{pro}\label{39} Let $H,Q\in\mathcal{S}^{n\times n}$ with $Q, H
\geq 0$ and suppose that there is a constant $C>0$ such that
\begin{equation}\label{40}
\frac{1}{C}\langle Q\xi,\xi\rangle\leq\langle H\xi,\xi\rangle\leq
C\langle Q\xi,\xi\rangle\qquad\text{ for all }\,\xi\in\R^n.
\end{equation}
Then there is an invertible matrix $M\in\text{Mat}(n,\R)$ such that
\begin{enumerate}
\item  $Q=M^THM$,
\item $\sqrt{Q}=\sqrt{H}M=M^T\sqrt{H}$,
\item $\frac{1}{\sqrt{C}}|\xi|\leq|M^T\xi|\leq
\sqrt{C}|\xi|\,\,\,\,\text{ for all }\xi\in\R^n$.
\end{enumerate}
\end{pro}

\textbf{Proof of Proposition \ref{39}:} Step 1. We claim that
$$\text{Ker}Q=\text{Ker}H,$$ i.e., $Q\xi=0$ if and only if $H\xi=0$
for any $\xi\in\R^n$. Indeed, suppose $Q\xi=0$. Then $\langle
H\xi,\xi\rangle=0$ by \eqref{40} and
$$|\sqrt{H}\xi|^2=\langle\sqrt{H}\xi,\sqrt{H}\xi\rangle=
\langle\sqrt{H}\sqrt{H}\xi,\xi\rangle=\langle H\xi,\xi\rangle=0,
$$
i.e., $\sqrt{H}\xi=0$ and $\xi$ is an eigenvector of $\sqrt{H}$
for the eigenvalue $0$. This in turn implies by Lemma \ref{38} that
$\xi$ is an eigenvector of $H$ for the eigenvalue $0$. Thus $\xi\in
\text{Ker}H$. The same argument shows that if $H\xi=0$ then
$\xi\in\text{Ker}Q$, and thus the claim is proved.

Step 2. Assume that one of the two matrices is invertible, i.e., has
empty kernel. Then by Step 1 the other matrix is also invertible, and so
by Lemma \ref{38}, both $\sqrt{H}$ and $\sqrt{Q}$ are invertible. In
this case, we may define
$$M=\big(\sqrt{H}\big)^{-1} \sqrt{Q}.$$
Then
$M^T=\sqrt{Q}^T\Big(\big(\sqrt{H}\big)^{-1}\Big)^T=\sqrt{Q}\big(\sqrt{H}\big)^{-1}$
and hence
\begin{itemize}
\item[i)] $M^THM=\sqrt{Q}\big(\sqrt{H}\big)^{-1}\sqrt{H}\sqrt{H}\big(\sqrt{H}\big)^{-1}\sqrt{Q}=Q$,
\item[ii)] $\sqrt{H}M=\sqrt{H}\big(\sqrt{H}\big)^{-1}\sqrt{Q}=\sqrt{Q}$,
\item[iii)] $M^T\sqrt{H}=\sqrt{Q}\big(\sqrt{H}\big)^{-1}\sqrt{H}=\sqrt{Q}$.
\end{itemize}
Thus in this case parts (1) and (2) of the proposition are satisfied.
Next note that \eqref{40} implies that both $H-\frac{1}{C}Q, CQ-H \geq 0.$
Then Lemma \ref{37}, part 2, implies that both of the following are
also nonnegative definite:
\begin{eqnarray}
\nonumber&&I_n-\frac{1}{C}\big(\sqrt{H}\big)^{-1}Q\big(\sqrt{H}\big)^{-1}=
\big(\sqrt{H}\big)^{-1}\Big(H-\frac{1}{C}Q\Big)\big(\sqrt{H}\big)^{-1},\\
\nonumber&&C\big(\sqrt{H}\big)^{-1}Q\big(\sqrt{H}\big)^{-1}-I_n=
\big(\sqrt{H}\big)^{-1}\Big(CQ-H\Big)\big(\sqrt{H}\big)^{-1}.
\end{eqnarray}
Equivalently, for every $\xi\in\R^n$,
\begin{eqnarray}
\nonumber&&|\xi|^2\geq\frac{1}{C}\langle\big(\sqrt{H}\big)^{-1}Q\big(\sqrt{H}\big)^{-1}\xi,\xi\rangle,\\
\nonumber&&C\langle\big(\sqrt{H}\big)^{-1}Q\big(\sqrt{H}\big)^{-1}\xi,\xi\rangle\geq|\xi|^2.
\end{eqnarray}
But then
$$|M^T\xi|^2=\langle M^T\xi,M^T\xi\rangle=\langle MM^T\xi,\xi\rangle=
\langle\big(\sqrt{H}\big)^{-1}Q\big(\sqrt{H}\big)^{-1}\xi,\xi\rangle\geq\frac{1}{C}|\xi|^2$$
and $|M^T\xi|^2\leq C|\xi|^2$. Thus we finally obtain that
$$\frac{1}{\sqrt{C}}|\xi|\leq|M^T\xi|\leq\sqrt{C}|\xi|,$$
which is part (3) of the statement. This proves the desired result in
case both $Q, H$ are invertible.

Step 3. It remains only to consider the case when neither $H$ nor $Q$ is
invertible since $\text{Ker}H=\text{Ker}Q$ by Step 1.  Since both matrices
are symmetric, each is diagonalizable. Moreover, eigenvectors
related to different eigenvalues are orthogonal in $\R^n$.

Consider the subspaces $V:=\text{Ker}H\subset\R^n$ and $V^\bot$. Then
$\R^n=V\oplus V^\bot$, and letting $k=\text{dim}V$, we have
$k\geq1$. If $k=n$ then $Q=H=0_n$ and the conclusion of Proposition \ref{39}
is obvious with $M=I_n$. Thus we may assume $k\leq n-1$.

Now choose an orthonormal basis $\{v_1,\ldots,v_k\}$ in $V$ and
another one $\{v_{k+1},\ldots,v_n\}$ in $V^\bot$. Then
$\mathcal{B}:=\{v_1,\ldots,v_k,v_{k+1},\ldots,v_n\}$ is an
orthonormal basis in $\R^n$. Let $\mathcal{B}^\prime:=\{e_1,\ldots,
e_n\}$ be the canonical basis in $\R^n$, and let $O$ be the matrix
which expresses the change of basis between $\mathcal{B}$ and
$\mathcal{B}^\prime$.  Then $O\in\mathcal{O}(n)$ and
$$Q=O^T\left(\begin{array}{cc}
               0_k & 0 \\
                 0 & Q_1\\
             \end{array}\right)O,\qquad\qquad
H=O^T\left(\begin{array}{cc}
               0_k & 0 \\
                 0 & H_1\\
             \end{array}\right)O.$$
Here $Q_1\in\text{Mat}(n-k,\R)$ is the invertible matrix associated to
the bijective linear map $T_Q:V^\bot\rightarrow V^\bot$ defined by
$T_Q(x) =Qx$, expressed with respect to the basis
$\{v_{k+1},\ldots,v_n\}$ of $V^\bot$. Similarly
$H_1\in\text{Mat}(n-k,\R)$ is the invertible matrix associated to the
bijective linear map $T_H:V^\bot\rightarrow V^\bot$ defined by
$T_H(x)=Hx$, also expressed with respect to the basis
$\{v_{k+1},\ldots,v_n\}$ of $V^\bot$.

Since $Q,H$ are symmetric and nonnegative definite, so are $Q_1,H_1$.
Then we can apply the result from Step 2 to find an invertible
$M_1\in\text{Mat}(n-k,\R)$ such that
\begin{itemize}
\item[i)] $M_1^TH_1M_1=Q_1$,
\item[ii)] $\sqrt{H_1}M_1=\sqrt{Q_1}=M_1^T\sqrt{H_1}$,
\item[iii)] $\frac{1}{\sqrt{C}}|w|_{\R^{n-k}}\leq|M_1w|_{\R^{n-k}}\leq
\sqrt{C}|w|_{\R^{n-k}}$ for all $w\in\R^{n-k}$.
\end{itemize}
Now define $$M=O^T\left(\begin{array}{cc}
                           I_k & 0 \\
                             0 & M_1 \\
                         \end{array}\right)O.$$
Then
\begin{eqnarray}
\nonumber M^THM&=&O^T\left(\begin{array}{cc}
                           I_k & 0 \\
                             0 & M^T_1 \\
                         \end{array}\right)O
                  O^T\left(\begin{array}{cc}
                           0_k & 0 \\
                             0 & H_1\\
                         \end{array}\right)O
                  O^T\left(\begin{array}{cc}
                           I_k & 0 \\
                             0 & M_1 \\
                         \end{array}\right)O\\
\nonumber&=&O^T\left(\begin{array}{cc}
                        0_k & 0 \\
                          0 & M_1^TH_1M_1\\
                     \end{array}\right)O\\
\nonumber&=&O^T\left(\begin{array}{cc}
                        0_k & 0 \\
                          0 & Q_1\\
                     \end{array}\right)O=Q.
\end{eqnarray}
Moreover, $\sqrt{Q}=O^T\left(\begin{array}{cc}
                                           0_k & 0 \\
                                           0 & \sqrt{Q_1} \\
                                         \end{array}\right)O$ and
$\sqrt{H}=O^T\left(\begin{array}{cc}
                                           0_k & 0 \\
                                           0 & \sqrt{H_1} \\
                                         \end{array}\right)O$, so
that $$\sqrt{H}M=O^T\left(\begin{array}{cc}
                            0_k & 0 \\
                            0 & \sqrt{H_1}M_1 \\
                          \end{array}\right)O=
                 O^T\left(\begin{array}{cc}
                            0_k & 0 \\
                            0 & \sqrt{Q_1} \\
                          \end{array}\right)O=\sqrt{Q}.$$
In the same way, $$M^T\sqrt{H}=O^T\left(\begin{array}{cc}
                                         0_k & 0 \\
                                         0 & M_1^T\sqrt{H_1} \\
                                       \end{array}\right)O=
                              O^T\left(\begin{array}{cc}
                                         0_k & 0 \\
                                         0 & \sqrt{Q_1} \\
                                       \end{array}\right)O=\sqrt{Q}.$$
Finally let $\xi\in\R^n$ and $\eta=O\xi$. Write $\eta=(v,w)$ with
$w\in\R^{n-k}$ and $v\in\R^k$. Then
$M^T\xi=O^T\left(\begin{array}{cc}
                    I_k & 0 \\
                    0 & M_1^T \\
                 \end{array}\right)\eta=O^T\left(\begin{array}{c}
                                                    v \\
                                                    M_1^Tw \\
                                                 \end{array}\right)$.
Since $O\in\mathcal{O}(n)$, we have
\begin{eqnarray}
\nonumber|M^T\xi|^2_{\R^n}&=&|v|^2_{\R^k}+|M_1^Tw|^2_{\R^{n-k}}\\
\nonumber&\leq&|v|^2_{\R^k}+C|w|^2_{\R^{n-k}}\\
\nonumber&\leq& C\big(|v|^2_{\R^k}+|w|^2_{\R^{n-k}}\big)\\
\nonumber&=&C|\eta|^2_{\R^n}=C|O\xi|^2_{\R^n}=C|\xi|^2_{\R^n}.
\end{eqnarray}
Similarly, $|M^T\xi|^2_{\R^n}\geq\frac{1}{C}|\xi|^2_{\R^n}$. Thus
$$\frac{1}{\sqrt{C}}|\xi|_{\R^n}\leq|M^T\xi|_{\R^n}\leq\sqrt{C}|\xi|_{\R^n},$$
and the proof of Proposition \ref{39} is complete. $\qquad\Box$

\begin{cor}\label{42} Let $Q(x)$ and $H(x)$ be symmetric nonnegative
definite matrices depending on $x\in\Omega$, and suppose there is a
constant $C>0$ so that
\begin{equation}\label{41}
\frac{1}{C}\langle Q(x)\xi,\xi\rangle\leq\langle
H(x)\xi,\xi\rangle\leq C\langle Q(x)\xi,\xi\rangle
\end{equation}
for all $\xi\in\R^n$ and a.e. $x\in\Omega$. Then for
a.e. $x\in\Omega$, there is an invertible matrix $M(x)$ such that
\begin{enumerate}
\item  $Q(x)=M^T(x)H(x)M(x)$,
\item $\sqrt{Q(x)}=\sqrt{H(x)}M(x)=M^T(x)\sqrt{H(x)}$,
\item $\frac{1}{\sqrt{C}}|\xi|\leq|M^T(x)\xi|\leq
\sqrt{C}|\xi|\,\,\,\,\text{ for all }\xi\in\R^n$.
\end{enumerate}
\end{cor}
\textbf{Proof:} This follows immediately by applying Proposition
\ref{39} at each point $x\in\Omega$ where \eqref{41}
holds. $\qquad\Box$

\vspace{0,1cm}

\textbf{Proof of Theorem \ref{32}:} Let $Q$ and $H$ satisfy the
hypothesis of Theorem \ref{32}. By Corollary \ref{42}, for
a.e. $x\in\Omega$, there is an invertible matrix $M(x)$ satisfying
properties (1), (2) and (3) relative to $Q$ and $H$. For any such
$x$, define
$$\hat{A}(x,z,\xi) =M(x)\tilde{A}(x,z,\xi).$$
Then by property (2) in Corollary \ref{42},
$$A(x,z,\xi)=\sqrt{Q(x)}\tilde{A}(x,z,\xi)=\sqrt{H(x)}M(x)
\tilde{A}(x,z,\xi)=\sqrt{H(x)}\hat{A}(x,z,\xi).$$
Therefore, by properties (2) and (3),
\begin{eqnarray}
\nonumber \xi \cdot A(x,z,\xi) &\geq& a^{-1}\Big|\sqrt{Q(x)}\cdot
\xi\Big|^p - h|z|^\gamma - g\\
\nonumber&\geq&a^{-1}C^{-\frac{p}{2}}\Big|\sqrt{H(x)}\cdot
\xi\Big|^p - h|z|^\gamma - g,\\
\nonumber\Big|B(x,z,\xi)\Big|&\leq&c\Big|\sqrt{Q(x)}\cdot\xi
\Big|^{\psi-1}+d|z|^{\delta-1}+f\\
\nonumber&\leq&C^\frac{\psi-1}{2}\,c\Big|\sqrt{H(x)}\cdot\xi
\Big|^{\psi-1}+d|z|^{\delta-1}+f,
\end{eqnarray}
where $a, b, c, d, e, f, g, h$ are as in \eqref{struct}.
In order to prove the third part of \eqref{struct2bis}, we first note
that
$$
|\eta|=\sup_{\zeta\in\R^n,\,|\zeta|=1}|\langle\zeta,\eta\rangle| \quad
 \mbox{for any $\eta\in\R^n$.}
$$
Then by property (3) in Corollary \ref{42},
\begin{eqnarray}
\nonumber|M(x)\eta|&=&\sup_{\zeta\in\R^n,\,|\zeta|=1}|\langle\zeta,M(x)\eta\rangle|\\
\nonumber&=& \sup_{\zeta\in\R^n,\,|\zeta|=1}|\langle
M^T(x)\zeta,\eta\rangle|\\
\nonumber&\leq&\sup_{\zeta\in\R^n,\,|\zeta|=1}|M^T(x)\zeta||\eta|\,\,\leq\,\,\sqrt{C}|\eta|.
\end{eqnarray}
Hence
\begin{eqnarray}
\nonumber\Big|\hat{A}(x,z,\xi)\Big|&=&\Big|M(x)\tilde{A}(x,z,\xi)\Big|\\
\nonumber&\leq&\sqrt{C}\Big|\tilde{A}(x,z,\xi)\Big|\\
\nonumber&\leq&\sqrt{C}\Big[ a\Big|\sqrt{Q(x)}\cdot\xi\Big|^{p-1} + b|z|^{\gamma -1} + e\Big]\\
\nonumber&\leq&C^\frac{p}{2}a\Big|\sqrt{H(x)}\cdot\xi\Big|^{p-1} +
C^\frac{1}{2}b|z|^{\gamma -1} + C^\frac{1}{2}e,
\end{eqnarray}
which completes the proof.$\qquad\Box$

\begin{pro}\label{43} For $x\in \Omega$, consider a symmetric
nonnegative definite matrix $Q(x)$ and a vector field
$T(x) =\sum_{j=1}^nt_j(x)\frac{\partial}{\partial
x_j}=\big(t_1(x),\ldots,t_n(x)\big)$ which is subunit with respect
to $Q(x)$, i.e.,
\begin{equation}\label{44}
\Big(\sum_{i=1}^nt_i(x)\xi_i\Big)^2\leq\langle
Q(x)\xi,\xi\rangle\qquad\text{for a.e $x\in\Omega$ and all
$\xi\in\R^n$.}
\end{equation}
Then there exists a vector $V(x)$ such that
\begin{enumerate}
\item $T(x)=\sqrt{Q(x)}V(x)$ for a.e. $x\in\Omega$,
\item $|V(x)|\leq1$ for a.e. $x\in\Omega$.
\end{enumerate}
\end{pro}
\textbf{Proof:} Consider any point $x_0\in\Omega$ at which \eqref{44}
holds with $x=x_0$ for every $\xi\in\R^n$. Denote
$T=T(x_0)=\big(t_1(x_0),\ldots,t_n(x_0)\big)$, $Q=Q(x_0)$ and $K=
\text{Ker}Q=\{\xi\in\R^n:Q\xi=0\}$. Write $\R^n=K\oplus K^\bot$, and
accordingly write $T=T_1+T_2$ with $T_1\in K$ and $T_2\in
K^\bot$. Then by \eqref{44} at $x_0$,
$$|\langle T,\xi\rangle|^2\leq\langle Q\xi,\xi\rangle\qquad \text{ for
all }\xi\in\R^n.$$
Choosing $\xi=T_1$ gives
$$|T_1|^2=|\langle T,T_1\rangle|^2\leq\langle QT_1,T_1\rangle=\langle 0,T_1\rangle=0,$$
hence $T=T_2\in K^\bot$.

We may assume that $K\subsetneqq\R^n$, since otherwise $Q=0$, $T=0$
and then the conclusion of the proposition holds at $x_0$ by choosing
$V(x_0)=0$. Now note that there is an orthogonal matrix
$O\in\mathcal{O}(n)$ such that $$Q=O^T\left(\begin{array}{cc}
                                              Q_1 & 0 \\
                                                0 & 0_k \\
                                            \end{array}\right)O,$$ where $k=\text{dim} K\geq0$ and all the
eigenvalues of $Q_1$ are strictly positive. Then
$$\sqrt{Q}=O^T\left(\begin{array}{cc}
                        \sqrt{Q_1} & 0 \\
                                 0 & 0_k \\
                    \end{array}\right)O.$$
Also, $Q_1$ is an invertible symmetric matrix which
corresponds to the invertible linear operator $L_{Q,K^\bot}$ defined
on $K^\bot$ by $L_{Q,K^{\bot}}(\xi)=Q\xi$. Hence we may define
$$N=O^T\left(\begin{array}{cc}
                 \big(\sqrt{Q_1}\big)^{-1} & 0 \\
                                         0 & 0_k \\
             \end{array}\right)O.$$
The matrix $N$ is symmetric and
$$\sqrt{Q}N=O^T\left(\begin{array}{cc}
                       I_{n-k} & 0 \\
                             0 & 0_k \\
                     \end{array}\right)O$$
corresponds to the canonical projection of $\R^n$ onto $K^\bot$. Since
$T\in K^\bot$, we have $T=\sqrt{Q}NT.$ Now set $NT=V$. Then
$T=\sqrt{Q}V$ and
\begin{equation*}
\begin{array}{l}
|V|^2\,=\,\sup_{|\xi|=1}|\langle
V,\xi\rangle|^2\,=\,\sup_{|\xi|=1}|\langle NT,\xi\rangle|^2
\,=\,\sup_{|\xi|=1}|\langle
T,N\xi\rangle|^2\vspace{0,2cm}\\
\hspace{2,5cm}\leq\,\sup_{|\xi|=1}\langle QN\xi,N\xi\rangle
\,=\,\sup_{|\xi|=1}\langle
\sqrt{Q}N\xi,\sqrt{Q}N\xi\rangle\,=\,1,
\end{array}
\end{equation*}
where the inequality follows from the fact that $T$ is subunit.
Thus the desired result holds at $x_0$ and the proof of Proposition
\ref{43} is complete.  $\qquad\Box$

\textbf{Proof of Theorem \ref{35}:} Rewrite \eqref{34} in the form
\begin{eqnarray}
\nonumber&&\text{div}\big(Q(x)\nabla u\big)+\sum_{i=1}^nS^\prime_i
G_iu-\sum_{i=1}^nT_i^\prime g_i=f-\sum_{i=1}^nH_iR_iu-Fu.
\end{eqnarray}
Since the vector fields $R_i$, $S_i$, $T_i$ are all subunit with
respect to $Q(x)$, Proposition \ref{43} shows that they can be
expressed as
\begin{equation}\label{45}
R_i(x)=\sqrt{Q(x)}\check{R}_i(x),\qquad
S_i(x)=\sqrt{Q(x)}\check{S}_i(x),\qquad
T_i(x)=\sqrt{Q(x)}\check{T}_i(x),
\end{equation}
where
\begin{equation}\label{46}
|\check{R}_i(x)|\leq1,\qquad|\check{S}_i(x)|\leq1,\qquad|\check{T}_i(x)|\leq1
\end{equation}
for every $i$ and a.e. $x\in\Omega$. Now write $R_i$, $S_i$, $T_i$,
$\check{R}_i$, $\check{S}_i$, $\check{T}_i$ as
\begin{equation*}
\begin{array}{ll}
\displaystyle R_i(x)=\sum_{j=1}^n R_{ij}(x)\frac{\partial}{\partial
x_j}=\big(R_{i1}(x),\ldots,R_{in}(x)\big),& \check{R}_i(x)=
\big(\check{R}_{i1}(x),\ldots,\check{R}_{in}(x)\big),
\end{array}
\end{equation*}
and similarly,
\begin{equation*}
\begin{array}{ll}
S_i(x)=\big(S_{i1}(x),\ldots,S_{in}(x)\big),& \check{S}_i(x)=
\big(\check{S}_{i1}(x),\ldots,\check{S}_{in}(x)\big), \\ 
T_i(x)=\big(T_{i1}(x),\ldots,T_{in}(x)\big),&\check{T}_i(x)=
\big(\check{T}_{i1}(x),\ldots,\check{T}_{in}(x)\big). 
\end{array}
\end{equation*}
For every $i,j=1,\ldots,n$ and a.e $x\in\Omega$, \eqref{45} gives
$$S_{ij}=\big(S_i\big)_j=\big(\sqrt{Q}\check{S}_i\big)_j=
\sum_{k=1}^n\big(\sqrt{Q}\big)_{jk}\big(\check{S}_i\big)_k= 
\sum_{k=1}^n\big(\sqrt{Q}\big)_{jk}\check{S}_{ik},$$ and in a
similar way,
$$R_{ij}=\sum_{k=1}^n\big(\sqrt{Q}\big)_{jk}\check{R}_{ik},\qquad
T_{ij}=\sum_{k=1}^n\big(\sqrt{Q}\big)_{jk}\check{T}_{ik},
$$
where in the notation we have suppressed dependence on $x$. Letting
$\check{R}$, $\check{S}$, $\check{T}$ denote respectively the matrices
$\big[\check{R}_{ij}\big]$, $\big[\check{S}_{ij}\big]$,
$\big[\check{T}_{ij}\big]$, we obtain for a.e. $x\in \Omega$ that
\begin{eqnarray}
\nonumber&&\mathbf{S^\prime G}u=\sum_{i=1}^nS_i^\prime G_i
u=-\text{div}\bigg(\Big(\sum_{i=1}^nS_{i1}G_iu,\ldots,\sum_{i=1}^n
S_{in}G_iu\Big)\bigg)\\   
\nonumber&&\hspace{1cm}=-\text{div}\bigg(\Big(\sum_{i,k=1}^n
\big(\sqrt{Q}\big)_{1k}\check{S}_{ik}G_iu,\ldots, 
\sum_{i,k=1}^n\big(\sqrt{Q}\big)_{nk}\check{S}_{ik}G_iu\Big)\bigg)\\
\nonumber&&\hspace{1cm}=-\text{div}\Big(\sqrt{Q}\check{S}^T\mathbf{G}u\Big).
\end{eqnarray}
In the same way,
 $$
\mathbf{T^\prime g}=-\text{div}\Big(\sqrt{Q}\check{T}^T\mathbf{g}\Big)
$$
a.e. in $\Omega$. On the other hand,
$$\mathbf{HR}u=\sum_{i,j=1}^nH_iR_{ij}\frac{\partial u}{\partial x_j}=
\sum_{i,j,k=1}^nH_i\big(\sqrt{Q}\big)_{jk}\check{R}_{ik}\frac{\partial
u}{\partial x_j}=\langle\mathbf{H},\check{R}\sqrt{Q}\nabla
u\rangle.$$

Then we can rewrite \eqref{34} as follows for a.e. $x\in\Omega$:
\begin{equation}\label{newform}
\text{div}\Big(Q(x)\nabla
u-\sqrt{Q}\check{S}^T\mathbf{G}u+\sqrt{Q}\check{T}^T\mathbf{g}\Big)
=f-\langle\mathbf{H},\check{R}\sqrt{Q}\nabla u\rangle-Fu.
\end{equation}
To compare this form with \eqref{eqdiff} and with the structural
conditions \eqref{struct} in case all of $p, \gamma, \psi, \delta$ are
equal to $2$, let
\begin{eqnarray}
\nonumber A(x,z,\xi)&=&Q(x)\xi-\sqrt{Q(x)}\check{S}^T(x)\mathbf{G}(x)z
+\sqrt{Q(x)}\check{T}^T(x)\mathbf{g}(x),\\
\nonumber \widetilde{A}(x,z,\xi)&=&\sqrt{Q(x)}\xi-\check{S}^T(x)
\mathbf{G}(x)z+\check{T}^T(x)\mathbf{g}(x),\\
\nonumber
B(x,z,\xi)&=&f(x)-\langle\mathbf{H}(x),\check{R}(x)\sqrt{Q(x)}
\xi\rangle-F(x)z.
\end{eqnarray}
Then $A(x,z,\xi)=\sqrt{Q(x)}\widetilde{A}(x,z,\xi)$ and
\eqref{newform} takes the form \eqref{eqdiff}. By \eqref{46}, for a.e
$x\in\Omega$ and every $\eta\in\R^n$,
$$|\check{R}\eta|^2=\sum_{i=1}^n\Big(\sum_{j=1}^n\check{R}_{ij}
\eta_j\Big)^2\leq\sum_{i=1}^n\bigg(\Big(\sum_{j=1}^n\check{R}_{ij}^2\Big)
\Big(\sum_{j=1}^n\eta_j^2\Big)\bigg)=n|\eta|^2,$$ and in the same way,
$$|\check{S}^T\eta|^2\leq n|\eta|^2,\qquad|\check{T}^T\eta|^2\leq n|\eta|^2.$$
Then for a.e. $x\in\Omega$ and every $z\in\R$ and $\xi\in\R^n$,
\begin{eqnarray}
\nonumber \xi\cdot A(x,z,\xi)&=&|\sqrt{Q(x)}\xi|^2-\langle\xi,
\sqrt{Q(x)}\check{S}^T(x)\mathbf{G}(x)z\rangle+
\langle\xi,\sqrt{Q(x)}\check{T}^T(x)\mathbf{g}(x)\rangle\\
\nonumber&\geq&|\sqrt{Q(x)}\xi|^2-|\sqrt{Q(x)}\xi||\check{S}^T(x)
\mathbf{G}(x)z|- |\sqrt{Q(x)}\xi||\check{T}^T(x)\mathbf{g}(x)|\\
\nonumber&\geq&|\sqrt{Q(x)}\xi|^2-\frac{1}{4}|\sqrt{Q(x)}\xi|^2
-4|\check{S}^T(x)\mathbf{G}(x)z|^2\\
\nonumber&&\hspace{4.5cm}-\frac{1}{4}|\sqrt{Q(x)}\xi|^2
-4|\check{T}^T(x)\mathbf{g}(x)|^2\\
\nonumber&\geq&\frac{1}{2}|\sqrt{Q(x)}\xi|^2-4n|\mathbf{G}(x)|^2
|z|^2-4n|\mathbf{g}(x)|^2.
\end{eqnarray}
Moreover,
\begin{eqnarray}
\nonumber |\widetilde{A}(x,z,\xi)|&\leq&|\sqrt{Q(x)}\xi|
+|\check{S}^T(x)\mathbf{G}(x)z|+|\check{T}^T(x)\mathbf{g}(x)|\\
\nonumber&\leq&|\sqrt{Q(x)}\xi|+\sqrt{n}|\mathbf{G}(x)||z|+
\sqrt{n}|\mathbf{g}(x)|,\\
\nonumber |B(x,z,\xi)|&\leq&|f(x)|+|\langle\mathbf{H}(x),\check{R}(x)
\sqrt{Q(x)} \xi\rangle|+|F(x)z|\\
\nonumber&\leq&|\mathbf{H}(x)||\check{R}(x)\sqrt{Q(x)}\xi|+|F(x)||z|+|f(x)|\\
\nonumber&\leq&\sqrt{n}|\mathbf{H}(x)||\sqrt{Q(x)}\xi|+|F(x)||z|+|f(x)|.
\end{eqnarray}
Thus the structural conditions \eqref{struct} hold with $p=\gamma =
\psi = \delta = 2$ and with
\begin{equation}\label{48}
\begin{array}{c}
a=2,\qquad d=|F(x)|,\qquad
g=4n|\mathbf{g}(x)|^2,\qquad h(x)=4n|\mathbf{G}(x)|^2\\
b=\sqrt{n}|\mathbf{G}(x)|,\qquad e=\sqrt{n}|\mathbf{g}(x)|,\qquad
c(x)=\sqrt{n}|\mathbf{H}(x)|,\qquad f=|f(x)|.
\end{array}
\end{equation}
This completes the proof of Theorem \ref{35}. $\qquad\Box$

\begin{rem}
By Theorems \ref{32} and \ref{35}, the linear equation
\eqref{34} satisfies the structural assumptions \eqref{struct} with
$p=\gamma=\psi=\delta=2$ not only with respect to $Q(x)$ but also with
respect to any other symmetric matrix $H(x)\geq 0$ such
that for a.e. $x\in\Omega$ and every $\xi\in\R^n$,
$$\frac{1}{C}\langle H(x)\xi,\xi\rangle\leq \langle
Q(x)\xi,\xi\rangle\leq C\langle H(x)\xi,\xi\rangle.$$
\end{rem}

\textbf{Proof of Theorem \ref{36}:} Step 1. It is easy to see that
if \eqref{struct} is satisfied, then \eqref{struct3} is also satisfied
with $\tilde{a}(x,z,\xi):=|\tilde{A}(x,z,\xi)|$ for $\xi\in\R^n$,
$z\in\R$ and a.e. $x\in\Omega$. Indeed, if \eqref{struct} holds, then
for every $\eta,\xi\in\R^n$, $z\in\R$ and a.e. $x\in\Omega$,
$$\big|\eta\cdot A(x,z,\xi)\big|=\big|\eta\cdot\sqrt{Q(x)}
\tilde{A}(x,z,\xi)\big|\leq\big|\sqrt{Q(x)}\eta\big|\big|
\tilde{A}(x,z,\xi)\big| =\big|\sqrt{Q(x)}\eta\big|\tilde{a}(x,z,\xi).
$$
Moreover,
$$\tilde{a}(x,z,\xi)=\big|\tilde{A}(x,z,\xi)\big|\leq a\big|
\sqrt{Q(x)}\xi\big|^{p-1} + b|z|^{\gamma -1} + e,$$
and thus \eqref{struct3} holds.

Step 2. We will now prove that \eqref{struct3} implies
\eqref{struct}. Fix any $x\in\Omega$ such that \eqref{struct3} is
satisfied for all $\xi,\eta\in\R^n$ and all $z\in\R$. \emph{Claim:}
$A(x,z,\xi)\in\big(\text{ker}Q(x)\big)^\bot$ for $\xi\in\R^n$,
$z\in\R$.  Indeed, define $K=\text{ker}Q(x)$ and recall from Lemma
\ref{38} that since $Q(x)$ is symmetric and nonnegative, then also
$K=\text{ker}\sqrt{Q(x)}$. Consider the decomposition $\R^n=K\oplus
K^\bot$ and write $A(x,z,\xi)=A_1+A_2$ with $A_1\in K$ and $A_2\in
K^\bot$. From the first inequality in \eqref{struct3} with $\eta=A_1$,
we get
\begin{equation*}
\nonumber|A_1|^2=A_1\cdot A_1=A_1\cdot (A_1+A_2)=A_1\cdot
A(x,z,\xi)\leq |\sqrt{Q(x)}A_1| {\tilde a}(x,z,\xi)=0.
\end{equation*}
Hence $A(x,z,\xi)=A_2\in K^\bot$, which proves the claim.

Now suppose that $k:=\text{dim}K<n$ and choose an orthogonal matrix
$O\in\mathcal{O}(n)$ such that
$$Q(x)=O^T\left(\begin{array}{cc}
                   Q_1 & 0 \\
                     0 & 0_k \\
                \end{array}\right)O,$$
with $Q_1$ symmetric, nonnegative and invertible.
Then $$\sqrt{Q(x)}=O^T\left(\begin{array}{cc}
                                \sqrt{Q_1} & 0 \\
                                         0 & 0_k \\
                            \end{array}\right)O.$$
Next define $$N(x)=O^T\left(\begin{array}{cc}
                               \big(\sqrt{Q_1}\big)^{-1} & 0 \\
                                                       0 & 1_k \\
                            \end{array}\right)O,$$
so that
$$\sqrt{Q(x)}N(x)=O^T\left(\begin{array}{cc}
                             1_{n-k} & 0 \\
                                   0 & 0_k \\
                           \end{array}\right)O,$$
i.e., the linear mapping $L:\R^n\rightarrow\R^n$ defined by
$L\eta =\sqrt{Q(x)}N(x)\eta$ is the canonical projection of $\R^n$
onto $K^\bot$. Since $A(x,z,\xi)\in K^\bot$, then
$$A(x,z,\xi)=\sqrt{Q(x)}N(x)A(x,z,\xi)=\sqrt{Q(x)}\tilde{A}(x,z,\xi)$$
where $\tilde{A}(x,z,\xi):=N(x)A(x,z,\xi)$. Hence
\begin{eqnarray}
\nonumber\big|\tilde{A}(x,z,\xi)\big|&=&\sup_{|\eta|=1}\big|\langle
\eta,\tilde{A}(x,z,\xi)\rangle\big|\\
\nonumber&=&\sup_{|\eta|=1}\big|\langle
\eta,N(x)A(x,z,\xi)\rangle\big|\,\,\,=\,\,\,
\sup_{|\eta|=1}\big|\langle N(x)\eta,A(x,z,\xi)\rangle\big|.
\end{eqnarray}
The last term can be estmated by using \eqref{struct3} to obtain
\begin{eqnarray}
\nonumber\big|\langle
N(x)\eta,A(x,z,\xi)\rangle\big|&\leq&\big|\sqrt{Q(x)}N(x)
\eta\big| \tilde{a}(x,z,\xi)\\
\nonumber &\leq&\big|\sqrt{Q(x)}N(x)\big||\eta|\tilde{a}(x,z,\xi)
\,\,\,=\,\,\,|\eta|\tilde{a}(x,z,\xi).
\end{eqnarray}
Therefore,
$$\big|\tilde{A}(x,z,\xi)\big|\leq\tilde{a}(x,z,\xi)\leq a\Big|
\sqrt{Q(x)}\cdot\xi\Big|^{p-1} + b|z|^{\gamma -1} + e.$$
The proof Theorem \ref{36} is now complete.$\qquad\Box$

\vspace{0,5cm}

\begin{center}
\textbf{Acknowledgement:}
\end{center}
The present work began while the first two authors were at the
Department of Mathematics of Rutgers University. They wish to
express their thanks for the stimulating environment and hospitality
that was provided. The first author also wishes to thank the
Nonlinear Analysis Center at Rutgers for having hosted his visit.

\vspace{0,5cm}

\bibliographystyle{plain}

\end{document}